\documentclass[reqno, huge]{amsart}
\usepackage[dvipsnames, table]{xcolor}
\usepackage{amsthm}
\theoremstyle{plain}
\usepackage{amssymb}
\usepackage{marvosym}
\usepackage{bm}
\usepackage{mathrsfs}
\usepackage{enumerate}  
\usepackage{mathtools}
\usepackage{tikz-cd}
\usepackage{graphicx}
\usepackage{float}
\usepackage[titletoc,title]{appendix}
\usepackage{marginnote}
\usepackage[disable]{todonotes}
\usepackage{etoolbox}
\usepackage[all]{xy}
\makeatletter
\patchcmd{\Ginclude@eps}{"#1"}{#1}{}{}
\makeatother
\usepackage[outdir=./]{epstopdf}
\usepackage[numbers]{natbib}
\usepackage[utf8]{inputenc}
\usepackage[english]{babel}
\usepackage{changepage}
\usepackage{makecell}
\usepackage{enumitem}

\definecolor{lightblue}{HTML}{1F88CD}
\definecolor{lightgrey}{HTML}{727272}
\definecolor{lightblue2}{HTML}{009EC1}
\definecolor{mypink}{HTML}{FD00B0}
\definecolor{lightred}{HTML}{ff4d4d}

\usepackage{hyperref}
\hypersetup{
	colorlinks=true,
    linkcolor={OliveGreen},
    citecolor={lightred},
	urlcolor={black}
}

\newtheorem*{theorem*}{Theorem}
\newtheorem{theorem}{Theorem}[section]
\newtheorem{corollary}[theorem]{Corollary}
\newtheorem{lemma}[theorem]{Lemma}
\newtheorem{conjecture}[theorem]{Conjecture}

\newtheorem{proposition}[theorem]{Proposition}
\theoremstyle{definition}

\theoremstyle{definition}
\newtheorem{definition}[theorem]{Definition}
\theoremstyle{definition}
\newtheorem{remark}[theorem]{Remark}
\theoremstyle{definition}

\theoremstyle{definition}

\theoremstyle{definition}

\theoremstyle{definition}
\newtheorem{question}[theorem]{Question}
\theoremstyle{definition}

\theoremstyle{definition}
\newtheorem{question!}[theorem]{Question!}
\theoremstyle{definition}

\makeatletter
\newcommand*\sbt{\mathpalette\sbt@{.75}}
\newcommand*\sbt@[2]{\mathbin{\vcenter{\hbox{\scalebox{#2}{$\m@th#1\bullet$}}}}}
\makeatother

\newcommand{\ra}{\rightarrow}
\newcommand{\xra}{\xrightarrow}

\newcommand{\sst}{\subset}

\newcommand{\bR}{\bm{\mathrm{R}}}
\newcommand{\bL}{\bm{\mathrm{L}}}

\newcommand{\D}{\mathrm{D}}

\newcommand{\ZZ}{\mathbb{Z}}

\newcommand{\CC}{\mathbb{C}}
\newcommand{\PP}{\mathbb{P}}

\newcommand{\ch}{\mathrm{ch}}

\newcommand{\HH}{\mathrm{HH}}

\newcommand{\pr}{\mathrm{pr}}
\newcommand{\ev}{\mathrm{ev}}

\newcommand{\GL}{\widetilde{\mathrm{GL}}^+(2,\mathbb{R})}

\renewcommand{\Re}{\operatorname{Re}}
\renewcommand{\Im}{\operatorname{Im}}

\DeclareMathOperator{\identity}{id}

\DeclareMathOperator{\rk}{rk}

\DeclareMathOperator{\Coh}{\mathrm{Coh}}

\DeclareMathOperator{\Ext}{Ext}

\DeclareMathOperator{\Hom}{Hom}
\DeclareMathOperator{\RHom}{RHom}

\DeclareMathOperator{\ext}{ext}

\DeclareMathOperator{\Pic}{Pic}

\DeclareMathOperator{\CaCl}{CaCl}

\DeclareMathOperator{\cone}{cone}
\DeclareMathOperator{\Stab}{Stab}
\DeclareMathOperator{\Gr}{Gr}

\newcommand{\cX}{\mathcal{X}}

\newcommand{\cC}{\mathcal{C}}
\newcommand{\cA}{\mathcal{A}}
\newcommand{\cE}{\mathcal{E}}

\newcommand{\cH}{\mathcal{H}}
\newcommand{\cB}{\mathcal{B}}

\newcommand{\cI}{\mathcal{I}}
\newcommand{\cT}{\mathcal{T}}
\newcommand{\cQ}{\mathcal{Q}}
\newcommand{\Ku}{\mathcal{K}u}
\newcommand{\cP}{\mathcal{P}}
\newcommand{\cD}{\mathcal{D}}

\newcommand{\cN}{\mathcal{N}}
\newcommand{\cM}{\mathcal{M}}
\DeclareMathOperator{\cF}{\mathcal{F}}

\DeclareMathOperator{\E}{\mathcal{E}}
\DeclareMathOperator{\oh}{\mathcal{O}}
\newcommand{\RHomb}{\mathrm{RHom}^{\bullet}}

\usetikzlibrary{decorations.pathmorphing}

\title[Categorical Torelli for GM threefolds]{Categorical Torelli theorems for Gushel--Mukai threefolds} 


\address{Department of Mathematics, Maison du Nombre, 6 avenue de la Fonte, L-4364 Esch-sur-Alzette, Luxembourg}
\email{augustinas.jacovskis@uni.lu}

\address{Max Planck Institute for Mathematics, Vivatsgasse 7, 53111 Bonn, Germany}
\email{xlin@mpim-bonn.mpg.de, lin-x18@tsinghua.org.cn}

\address{School of Mathematical Sciences, Zhejiang University, Hangzhou, Zhejiang Province 310030, P. R. China}
\email{jasonlzy0617@gmail.com}

\address{Simons Laufer Mathematical Sciences
Institute, Berkeley, CA 94720, USA}
\address{Institut de Mathématiqes de Toulouse, UMR 5219, Université de Toulouse, Université Paul Sabatier, 118 route de
Narbonne, 31062 Toulouse Cedex 9, France}
\email{Shizhuozhang@msri.org,shizhuo.zhang@math.univ-toulouse.fr}

\author{Augustinas Jacovskis, Xun Lin, Zhiyu Liu and Shizhuo Zhang}
\date{\today}

\begin{document}

\begin{abstract}
We show that a general ordinary Gushel--Mukai (GM) threefold $X$ can be reconstructed from its Kuznetsov component $\Ku(X)$ together with an extra piece of data coming from tautological subbundle of the Grassmannian $\mathrm{Gr}(2,5)$. We also prove that $\Ku(X)$ determines the birational isomorphism class of $X$, while $\Ku(X')$ determines the isomorphism class of a special GM threefold $X'$ if it is general. As an application, we prove a conjecture of Kuznetsov--Perry in dimension three under a mild assumption. Finally, we use $\Ku(X)$ to restate a conjecture of Debarre--Iliev--Manivel regarding fibers of the period map for ordinary GM threefolds.
\end{abstract}

\maketitle

{
\hypersetup{linkcolor=black}
\setcounter{tocdepth}{1}
\tableofcontents
}

\section{Introduction}

In recent times, derived categories have played an important role in algebraic geometry; in many cases, much of the geometric information of a variety/scheme $X$ is encoded by its bounded derived category of coherent sheaves $\D^b(X)$. In this setting, one of the most fundamental questions that can be asked is whether $\D^b(X)$ recovers $X$ up to isomorphism, in other words, whether a \emph{derived Torelli theorem} holds for $X$. For varieties with ample or anti-ample canonical bundle (which include Fano varieties and varieties of general type), this question was answered affirmatively by Bondal--Orlov in \cite{bondal2001reconstruction}.

\subsection{Kuznetsov components and categorical Torelli theorems}

Therefore, for the class of varieties above, it is natural to ask whether they are also determined up to isomorphism by \emph{less} information than the whole derived category $\D^b(X)$. A natural candidate for this is a subcategory $\Ku(X)$ of $\D^b(X)$ called the \emph{Kuznetsov component}.
This subcategory has been studied extensively by Kuznetsov and others (e.g. \cite{kuznetsov2003derived, kuznetsov2009derived, kuznetsov2018derived}) for many Fano varieties, including Gushel--Mukai (GM) varieties. 

The question of whether $\Ku(X)$ determines $X$ up to isomorphism has been studied for certain cases in the setting of Fano threefolds. In \cite{bernardara2012categorical}, the authors show that the Kuznetsov component completely determines cubic threefolds up to isomorphism, in other words, a \emph{categorical Torelli theorem} holds for cubic threefolds $Y$. The same result was also verified in \cite{pertusi2020some}. On the other hand, for many Fano varieties, the Kuznetsov component $\Ku(X)$ does not determine the isomorphism class, but only the birational isomorphism class of $X$. This is known as a \emph{birational categorical Torelli theorem}. For instance, Kuznetsov components determine the birational isomorphism class of every index $1$ prime Fano threefolds of even genus $g\geq 8$. For GM threefolds -- the focus of our paper -- by \cite{kuznetsov2019categorical} it is known that there are birational GM threefolds with equivalent Kuznetsov components. So there are two natural questions to ask in this setting:

\begin{question}\label{intro questions about Kuznetsov torelli}\leavevmode 
\begin{enumerate}
    \item Does $\Ku(X)$ determine the birational equivalence class of $X$?
    \item What extra data along with $\Ku(X)$ do we need to identify a particular GM threefold $X$ from its birational equivalence class?
\end{enumerate}
\end{question}

\subsection{Main Results}

\subsubsection{(Refined) categorical Torelli for Gushel--Mukai threefolds}

In the present paper, we deal with the case of index $1$ prime Fano threefolds of degree $10$ and genus $6$, also known as Gushel--Mukai threefolds (GM threefolds for short), which are split into two types: ordinary GM threefolds which arise as a quadric section of a linear section of the Grassmannian $\mathrm{Gr}(2,5)$, and special GM threefolds which arise as double covers of a codimension three linear section of $\mathrm{Gr}(2,5)$, branched over a degree ten K3 surface. By \cite{kuznetsov2018derived}, we have a semiorthogonal decomposition
\[\D^b(X)=\langle \Ku(X), \cE, \oh_X\rangle,\]
where $\cE$ is the pull-back of the tautological subbundle on $\Gr(2,5)$ along the natural map $X\to \Gr(2,5)$.

Our first main theorem is concerned with ordinary GM threefolds and answers Question \ref{intro questions about Kuznetsov torelli} (2):

\begin{theorem}[{Theorem \ref{refined categorical torelli theorem}}] \label{intro gluing data theorem}
Let $X$ be a general ordinary GM threefold and $\pi \colon \D^b(X) \ra \Ku(X)$ be the right adjoint to the inclusion $\Ku(X) \sst \D^b(X)$. Then the data of $\Ku(X)$ along with the object $\pi(\cE)$ is enough to determine $X$ up to isomorphism.
\end{theorem}


On the other hand, for special GM threefolds which are general (``general special" for short), we show that a categorical Torelli theorem holds:

\begin{theorem}[{Theorem \ref{theorem_torelli_sGM}}] \label{intro sgm torelli}
Let $X$ and $X'$ be general special GM threefolds, and assume that there is an equivalence of categories $\Ku(X) \simeq \Ku(X')$. Then $X$ and $X'$ are isomorphic.
\end{theorem}

\subsubsection{Birational categorical Torelli for Gushel--Mukai threefolds}

Next, returning to the setting of ordinary GM threefolds, we show that a birational categorical Torelli theorem holds for general ordinary GM threefolds, which answers Question~\ref{intro questions about Kuznetsov torelli} (1).

\begin{theorem}[{Theorem \ref{birational OGM torelli theorem}}] \label{intro birational ogm torelli theorem}
Let $X$ and $X'$ be general ordinary GM threefolds, and suppose that there is an equivalence of categories $\Ku(X)\simeq \Ku(X')$. Then $X$ is birationally equivalent to $X'$.
\end{theorem}

In \cite{kuznetsov2019categorical}, the authors studied GM varieties of arbitrary dimension and proved the Duality Conjecture \cite[Conjecture 3.7]{kuznetsov2018derived} for them, i.e. they showed that the period partner or period dual of a GM variety $X$ shares the same Kuznetsov component $\Ku(X)$ as $X$. Combining earlier results \cite[Theorem 4.20]{debarre2015gushel} on the birational equivalence of these varieties, this gives strong evidence for the following conjecture:

\begin{conjecture}[{\cite[Conjecture 1.7]{kuznetsov2019categorical}}] \label{conjecture_KP2019}
If $X$ and $X'$ are GM varieties of the same dimension such that there is an equivalence $\Ku(X)\simeq \Ku(X')$, then $X$ and $X'$ are birationally equivalent. 
\end{conjecture}

Thus our result Theorem \ref{intro birational ogm torelli theorem} actually proves Conjecture \ref{conjecture_KP2019} under the assumption that $X$ and $X'$ are both of dimension $3$, ordinary and general.

Moreover, by a careful study of Bridgeland moduli spaces of stable objects in the Kuznetsov components $\cA_X$ for not only smooth ordinary GM threefolds but also special GM threefolds $X$, we can prove that the Kuznetsov component of a general ordinary GM threefold can not be equivalent to the one of a general special GM threefold. Therefore, combined with Theorem \ref{intro birational ogm torelli theorem} and \ref{intro sgm torelli}, we have the following improved version of Theorem \ref{intro birational ogm torelli theorem}, which allows threefolds to be either ordinary or special:

\begin{theorem}[{Theorem~\ref{theorem_KP_conjecture} and Corollary~\ref{prop_general_preserves}}]\label{KP_conjecture_general}
If $X$ and $X'$ are general ordinary or general special GM threefolds such that there is an equivalence $\Ku(X)\simeq \Ku(X')$, then $X$ and $X'$ are birationally equivalent. 
\end{theorem}

\subsubsection{The Debarre--Iliev--Manivel Conjecture}

In \cite{debarre2012period}, the authors conjecture that the general fiber of the classical period map from the moduli space of ordinary GM threefolds to the moduli space of $10$ dimensional principally polarised abelian varieties is birational to the disjoint union of the minimal model $\cC_m(X)$ of the Fano surface of conics and a moduli space of stable sheaves $M_G^X(2,1,5)$, both quotiented by involutions, which we call the \emph{Debarre--Iliev--Manivel Conjecture} (cf.~Conjecture \ref{DIM conjecture}). 
Within the moduli space of smooth ordinary GM threefolds, we define the fiber of the ``categorical period map" through $[X]$ as the isomorphism classes of all ordinary GM threefolds $X'$ whose Kuznetsov components satisfy $\Ku(X')\simeq\Ku(X)$. Then the following categorical analogue of the \emph{Debarre--Iliev--Manivel conjecture} follows from Theorem~\ref{intro birational ogm torelli theorem} and results on Bridgeland moduli spaces with respect to the two $(-1)$-classes in the numerical Grothendieck group of $\cA_X$. 

\begin{theorem}[Theorem \ref{categorical period map}] \label{intro categorical period map}
A general fiber of the ``categorical period map" through an ordinary GM threefold $X$ is the union of $\mathcal{C}_m(X)/\iota$ and $M_G^X(2,1,5)/\iota'$ where $\iota, \iota'$ are geometrically meaningful involutions.
\end{theorem}

As an application, the \emph{Debarre--Iliev--Manivel} Conjecture \ref{DIM conjecture} can be restated in an equivalent form as follows:

\begin{conjecture}
\label{DIM_equi_conjecture}
Let $X$ be a general ordinary GM threefold. The intermediate Jacobian $J(X)$ determines the Kuznetsov component $\Ku(X)$.
\end{conjecture}

\begin{remark}
In \cite{debarre2012period}, the authors actually conjecture that a general fiber of the period map is \emph{birational} to the disjoint union of two surfaces, parametrizing conic transforms and conic transforms of a line transform of $X$, which is birational to the disjoint union of $\cC_m(X)$ and $M_G^X(2,1,5)$, both quotiented by involutions. In Corollary \ref{line_trans} we show that this birational equivalence is indeed an \emph{isomorphism}.
\end{remark}

\subsubsection{Uniqueness of Serre-invariant stability conditions}
One of the key steps when we identify Bridgeland moduli spaces via an equivalence of Kuznetsov components in the proofs of Theorems \ref{intro gluing data theorem} and \ref{intro birational ogm torelli theorem}. A stability condition $\sigma$ on the Kuznetsov component $\Ku(X)$ of a prime Fano threefold $X$ is \emph{Serre-invariant} if $S_{\Ku(X)}\cdot\sigma=\sigma\cdot g$ for some $g\in\widetilde{\mathrm{GL}}^+(2,\mathbb{R})$ (see Section \ref{section_double_tilted_stability_condition}). Serre-invariance is one of the fundamental tools in studying relationship of classical Gieseker moduli spaces and Bridgeland moduli spaces for Kuznetsov components (cf.~\cite{altavilla2019moduli,pertusi2020some,Zhang2020Kuzconjecture,liu2021note,FeyzbakhshPertusi2021stab}). A natural question is whether any two Serre-invariant stability conditions are in the same $\widetilde{\mathrm{GL}}^+(2,\mathbb{R})$-orbit. In the present paper, we answer this question affirmatively. 

\begin{theorem}[Theorem \ref{all_in_one_orbit}] \label{uniquness_Serre_invariant_introduction}
Let $X$ be a prime Fano threefold of index $1$ of genus $g\geq 6$, or a del Pezzo threefold of degree $d\geq 2$.
Then all Serre-invariant stability conditions on $\Ku(X)$ are in the same $\widetilde{\mathrm{GL}}^+(2,\mathbb{R})$-orbit.
\end{theorem}


\medskip

\subsection{Methods}

For convenience, we work with the alternative Kuznetsov component $\cA_X$, defined by the semiorthogonal decomposition $\D^b(X) = \langle \cA_X , \oh_X , \cE^\vee \rangle$ and there is an equivalence $\Xi\colon\Ku(X) \simeq \cA_X$. We prove the above theorems \ref{intro gluing data theorem}, \ref{intro birational ogm torelli theorem}, \ref{KP_conjecture_general} and \ref{intro categorical period map} by considering the moduli spaces of Bridgeland stable objects in the alternative Kuznetsov component $\cA_X$ with respect to $(-1)$-classes in the numerical Grothendieck group of $\cA_X$, i.e. a vector $v$ with $\chi(v,v)=-1$ where $\chi$ is the Euler form. Up to sign, there are two $(-1)$-classes in the numerical Grothendieck group of $\cA_X$, call them $-x$ and $y-2x$. 

First, we show that the moduli space with the class $-x$ is isomorphic to the minimal model $\cC_m(X)$ of the Fano surface of conics (Theorem~\ref{theorem_irreduciblecomponent_conics}). Indeed, we first show that the unique exceptional curve contracted in $\cC(X)$ is the rational curve of conics whose ideal sheaf $I_C$ is not in $\cA_X$ and that the image is the smooth point represented by $\pi(\cE)$ (Proposition~\ref{prop_contraction_P1_OGM}), so $\mathcal{C}_m(X)$ forms an irreducible component of the moduli space $\mathcal{M}_{\sigma}(\cA_X,-x)$ of stable objects in $\cA_X$ with respect to $-x$.

Then we show this component actually occupies the whole moduli space (Proposition \ref{moduli -x bijective}), which is the most difficult and technical part of the article and we only briefly sketch the argument here. We start with a stable object $F\in\cA_X$ of the class $-x$. It suffices to show that $F$ is isomorphic to the projection of ideal sheaf $I_C$ of a conic $C\subset X$. First, we assume that $F$ is semistable in the double tilted heart $\Coh^0_{\alpha, \beta}(X)$ (cf.~Section~\ref{section_double_tilted_stability_condition}). Then by a wall-crossing argument, we prove that $F[-1]$ is a slope-semistable sheaf of rank one. Since its class is $[F]=-[I_C]$, we get $F\cong I_C[1]$. Next, we assume that $F$ is not semistable in the double tilted heart $\Coh^0_{\alpha, \beta}(X)$. Our main tools are inequalities in \cite{li2016stability} and Theorem \ref{Naoki BG}, which allow us to bound the rank and first two Chern characters $\ch_1,\ch_2$ of the destabilizing objects and their cohomology objects. Since $F\in \cA_X$, by using the Euler characteristics $\chi(\mathcal{O}_X,-)$ and $\chi(\cE^{\vee},-)$ we can obtain a bound on $\ch_3$. Then we deduce that the Harder--Narasimhan factors of $F$ are the expected ones (Proposition \ref{F iso prIC}). As a result,  $\mathcal{M}_{\sigma}(\cA_X,-x)\cong\mathcal{C}_m(X)$. Similarly, we identify the moduli space $M_G^X(2,1,5)$ of Gieseker semistable sheaves of rank $2$, $c_1=1, c_2=5$ and $c_3=0$ on $X$ with the Bridgeland moduli space $\cM_\sigma(\cA_X, y-2x)$ in Theorem~\ref{irreducible component pr' theorem}.

As we have seen, $\cC(X)$ is exactly the blow-up of $\cC_m(X)\cong \mathcal{M}_{\sigma}(\cA_X,-x)$ at the point $\Xi(\pi(\cE))$, hence the data $(\mathcal{K}u(X), \pi(\cE))$ determines $\mathcal{C}(X)$. A classical result of Logachev \cite{logachev2012fano} states that $X$ can be determined up to isomorphism from $\cC(X)$. Thus Theorem \ref{intro gluing data theorem} is proved.

We prove Theorem \ref{intro sgm torelli} via another method. By considering the equivariant Kuznetsov components $\Ku(X)^{\mu_2}$, first discussed in \cite{KP2017}, and exploiting the fact that $X$ is the double cover of a degree $5$ index $2$ prime Fano threefold $Y$, branched over a quadric hypersurface $\mathcal{B} \sst Y$. In this case, the equivariant Kuznetsov component is equivalent to $\D^b(\mathcal{B})$ where $\mathcal{B}$ is a K3 surface. Therefore, a number of results concerning the Fourier--Mukai partners of K3 surfaces can be used to deduce that $\Ku(X)^{\mu_2} \simeq \Ku(X')^{\mu_2}$ implies $\mathcal{B} \cong \mathcal{B}'$. Then the fact that the del Pezzo threefold $Y$ of degree $5$ is rigid can be used to deduce that indeed, $X \cong X'$.

To prove Theorem \ref{intro birational ogm torelli theorem}, we invoke a few more results from \cite{debarre2012period}. More precisely, an equivalence of categories $\Phi\colon \cA_X \simeq \cA_{X'}$ identifies the moduli space $\mathcal{M}_{\sigma}(\cA_X,-x)$ with either $\cM_\sigma(\cA_X,-x)$ or $\cM_\sigma(\cA_{X'}, y-2x)$. The former case gives an isomorphism of minimal surfaces $\cC_m(X) \cong \cC_m(X')$. Blowing $\cC_m(X)$ up at the smooth point associated to $\pi(\cE)$ gives $\cC(X)$, and blowing up $\cC_m(X')$ at the image of $\pi(\cE)$ under $\Phi$ gives $\cC(X'_C)$, where $X'_C$ is certain birational transformation of $X'$, associated with a conic $C\subset X'$. Then by Logachev's Reconstruction Theorem for $\cC(X)$, $X$ is isomorphic to $X'_C$ which is birational to $X'$. For the latter case, we start with the isomorphism $\cC_m(X) \cong M_G^{X'}(2,1,5)$. In fact, $M_G^{X'}(2,1,5)$ is birational to $\cC(X'_L)$, where $X'_L$ is another birational transformation of $X'$, associated with a line $L \sst X'$. Since $\cC(X'_L)$ is a surface of general type, we get $\cC_m(X) \cong \cC_m(X_L')$. Then by the same argument as in the previous case, $X$ is isomorphic to some birational transformation of $X'$.

Finally, the proof of Theorem \ref{KP_conjecture_general} is similar to that of Theorem~\ref{intro birational ogm torelli theorem}. Firstly, we identify the Bridgeland moduli spaces $\mathcal{M}_{\sigma}(\cA_{X'}, -x)$ and $\mathcal{M}_{\sigma}(\cA_{X'}, y-2x)$ on a special GM threefold $X'$ with $\mathcal{C}_m(X')$ and $M^{X'}_G(2,1,5)$ respectively (Theorem~\ref{theorem_irreduciblecomponent_conics} and Theorem~\ref{irreducible component pr' theorem}), where $\mathcal{C}_m(X')$ is the contraction of the Fano surface $\mathcal{C}(X')$ of conics on $X'$ along one of the components to a singular point. Then if $X$ is ordinary, the equivalence $\Phi\colon\cA_{X}\simeq\cA_{X'}$ would identify those moduli spaces on a general ordinary GM threefold $X$ with those on a special GM threefold $X'$; we show that this is impossible by analyzing their singularities. Then Theorem~\ref{KP_conjecture_general} reduces to Theorem~\ref{intro birational ogm torelli theorem} and Theorem~\ref{intro sgm torelli}.

\subsection{Related Work}
\subsubsection{Categorical Torelli theorems}
There is a very nice survey article \cite{pertusi2022categorical} on recent results and remaining open questions on this topic. In \cite{bernardara2012categorical} and \cite{pertusi2020some}, the authors prove categorical Torelli theorems for cubic threefolds. In \cite{altavilla2019moduli} and \cite{bernardara2013semi}, the authors prove categorical Torelli theorems for general quartic double solids. In \cite{li2021refined} and \cite{li2022refined}, the authors prove a refined categorical Torelli theorem for Enriques surfaces. In \cite{jz2021brillnoether}, the authors generalize Theorem~\ref{refined categorical torelli theorem} to all prime Fano threefolds of genus $g\geq 6$. In \cite{GLZ2021conics}, the authors prove a birational categorical Torelli theorem for general non-Hodge-special Gushel--Mukai fourfolds. 

\subsubsection{Identifying classical moduli spaces as Bridgeland moduli spaces for Kuznetsov components}
In the present article, we realize the Fano surface of conics and a certain Gieseker moduli space of semistable sheaves as Bridgeland moduli spaces of stable objects in Kuznetsov components of GM threefolds. In \cite{pertusi2020some}, the authors realize the Fano surface of lines $\Sigma(Y_d)$ (for $d\geq 2$) as a Bridgeland moduli space of stable objects in the Kuznetsov component $\Ku(Y_d)$. In \cite{liu2021note}, the authors realize the moduli space of rank two instanton sheaves on a del Pezzo threefold $Y_d$ (for $d\geq 3$) and the compactification of the moduli space of ACM sheaves on $X_{4d+2}$ (for $d\geq 3$) as Bridgeland moduli spaces of stable objects in $\Ku(Y_d)$ and $\Ku(X_{4d+2})$, respectively. In \cite{FeyzbakhshPertusi2021stab}, the authors realize the moduli space of Ulrich bundles of arbitrary rank on a cubic threefold $Y_3$ as an open locus of a Bridgeland moduli space of stable objects in $\Ku(Y_3)$. 

\subsubsection{Serre-invariant stability conditions}
In \cite{Pertusi2021serreinv} and \cite{pertusi2020some}, the authors prove that stability conditions on Kuznetsov components of every del Pezzo threefold $Y_d$ of degree $d\geq 1$ and every index $1$ prime Fano threefold of genus $g\geq 6$ are Serre-invariant.  In \cite{FeyzbakhshPertusi2021stab}, the authors prove the uniqueness of Serre-invariant stability conditions for a general triangulated category satisfying a list of very natural assumptions, which include Kuznetsov components of a series of prime Fano threefolds.

\subsection{Notation and conventions}

\begin{itemize}
\item We work over the field $k=\CC$. All triangulated categories and abelian categories are assumed to be $k$-linear.

\item We use $\hom$ and $\ext^{i}$ to represent the dimension of the vector spaces $\Hom$ and~$\Ext^{i}$.

   \item The numerical $K$ group of a triangulated category $\cD$ is denoted by $\cN(\cD)$, which is the Grothendieck group $K_0(\cD)$ modulo the kernel of the Euler form $\chi(E, F) = \sum_i (-1)^i \ext^i(E, F)$
   
   \item We denote the bounded derived category of a smooth projective variety $X$ by $\D^b(X)$. The derived dual functor is denoted by $\mathbb{D}:=\mathrm{R}\mathcal{H}om_X(-, \oh_X)$.

    
    \item We denote the phase and slope with respect to a weak stability condition $\sigma$ by $\phi_{\sigma}$ and $\mu_{\sigma}$, respectively. The maximal and minimal slopes (phases) of the Harder--Narasimhan factors of a given object $F$ will be denoted by $\mu_{\sigma}^+(F)$ ($\phi^+_{\sigma}(F)$) and $\mu_{\sigma}^-(F)$ ($\phi^-_{\sigma}(F)$), respectively.
    
    \item $\cH^i_{\cA}$ means the $i$-th cohomology with respect to the heart $\cA$. When the $\cA$ subscript is dropped, we take the heart to be $\Coh(X)$.
    \item The symbol $\simeq$ denotes an equivalence of categories and a birational equivalence of varieties. The symbol $\cong$ denotes an isomorphism of varieties.
    \item Let $X$ be a GM threefold. Then a conic means a closed subscheme $C \sst X$ with Hilbert polynomial $p_C(t)=1+2t$, and a line means a closed subscheme $L \sst X$ with Hilbert polynomial $p_L(t)=1+t$.
\end{itemize}

\subsection{Organization of the paper}
In Section \ref{sod section}, we collect basic facts about semiorthogonal decompositions. In Section \ref{GM threefolds and their derived categories section}, we introduce Gushel--Mukai threefolds and their Kuznetsov components. In Section \ref{stability section}, we introduce the definition of weak stability conditions on $\D^b(X)$, and the induced stability conditions on the alternative Kuznetsov components $\cA_X$ of GM threefolds. In Section \ref{gluing data machine section}, we introduce a distinguished object $\pi(\cE)\in\mathcal{K}u(X)$ and its alternative Kuznetsov component analogue $\Xi(\pi(\cE))\in\cA_X$ and prove its stability. In Section \ref{hilbert schemes section} we discuss the geometry of the Fano surface of conics of a GM threefold. In Section \ref{refined categorical Torelli section}, we construct the Bridgeland moduli space of $\sigma$-stable objects with class $-x$ in $\cA_X$. In Section \ref{M_G(2,1,5) section}, we construct the Bridgeland moduli space of $\sigma$-stable objects with respect to the other $(-1)$-class $y-2x$ in $\cA_X$. In Section \ref{torelli_sec}, we prove several birational/refined categorical Torelli theorems (Theorems~\ref{intro gluing data theorem}, \ref{intro sgm torelli} and \ref{intro birational ogm torelli theorem}) and Conjecture~\ref{conjecture_KP2019} in dimension three with mild assumptions. In Section \ref{DIM conjecture section}, we describe the general fiber of the ``categorical period map" for ordinary GM threefolds \ref{intro categorical period map}, and restate the \emph{Debarre--Iliev--Manivel conjecture} in terms of Conjecture~\ref{DIM conjecture equivalent prop}. Finally, we study Serre-invariant stability conditions on Kuznetsov components and show that they are contained in one $\widetilde{\mathrm{GL}}^+(2,\mathbb{R})$ orbit in Appendix \ref{appendix}.

\subsection*{Acknowledgements}
Firstly, it is our pleasure to thank Arend Bayer for very useful discussions on the topics of this project. We would like to thank Sasha Kuznetsov for answering many of our questions on Gushel--Mukai threefolds. We thank Atanas Iliev, Laurent Manivel,  Daniele Faenzi, Dmitry Logachev, Will Donovan, Bernhard Keller, Alexey Elagin, Xiaolei Zhao, Chunyi Li, Laura Pertusi, Song Yang, Alex Perry, Pieter Belmans, Qingyuan Jiang, Enrico Fatighenti, Naoki Koseki, Bingyu Xia, Yong Hu and Luigi Martinelli for helpful conversations on several related topics. We would like to thank Daniele Faenzi for sending us the preprint \cite{Faenzi2021genusten} and Soheyla Feyzbakhsh for sending us the preprint \cite{FeyzbakhshPertusi2021stab}. We thank Pieter Belmans for useful comments on the first draft of our article. The last author would like to thank Tingyu Sun for constant support and encouragement. We would like to thank the anonymous referee for their careful reading of
our paper, and for their very useful and insightful comments.

The first and last authors are supported by ERC Consolidator Grant WallCrossAG, no. 819864. The first author was also supported by the Luxembourg National Research Fund (FNR–17113194).

\section{Semiorthogonal decompositions} \label{sod section}

In this section, we collect some useful facts about semiorthogonal decompositions. Background on triangulated categories and derived categories of coherent sheaves can be found in \cite{huybrechts2006fourier}, for example. From now on, let $\D^b(X)$ denote the bounded derived category of coherent sheaves on a smooth projective variety $X$, and for $E, F \in \D^b(X)$, define 
\[ \RHom^\bullet(E, F) = \bigoplus_{i \in \ZZ} \Hom(E, F[i])[-i] .  \]

\subsection{Exceptional collections and semiorthogonal decompositions}

\begin{definition}
Let $\cD$ be a triangulated category and $E \in \cD$. We say that $E$ is an \emph{exceptional object} if $\RHom^\bullet(E, E) = k$. Now let $(E_1, \dots, E_m)$ be a collection of exceptional objects in $\cD$. We say it is an \emph{exceptional collection} if $\RHom^\bullet(E_i, E_j) = 0$ for $i > j$.
\end{definition}

\begin{definition}
Let $\cD$ be a triangulated category and $\cC$ be a triangulated subcategory of $\cD$. We define the \emph{right orthogonal complement} of $\cC$ in $\cD$ as the full triangulated subcategory
\[ \cC^\bot = \{ X \in \cD \mid \Hom(Y, X) =0  \text{ for all } Y \in \cC \}.  \]
The \emph{left orthogonal complement} is defined similarly, as 
\[ {}^\bot \cC = \{ X \in \cD \mid \Hom(X, Y) =0  \text{ for all } Y \in \cC \}.  \]
\end{definition}

\begin{definition}
Let $\cD$ be a triangulated category. We say a triangulated subcategory $\cC \sst \cD$ is \emph{admissible} if the inclusion functor $i \colon \cC \hookrightarrow \cD$ has left adjoint $i^*$ and right adjoint $i^!$.
\end{definition}

\begin{definition}
Let $\cD$ be a triangulated category, and $( \cC_1, \dots, \cC_m  )$ be a collection of full admissible subcategories of $\cD$. We say that $\cD = \langle \cC_1, \dots, \cC_m \rangle$ is a \emph{semiorthogonal decomposition} of $\cD$ if $\cC_j \sst \cC_i^\bot $ for all $i > j$, and the subcategories $(\cC_1, \dots, \cC_m )$ generate $\cD$, i.e. the category resulting from taking all shifts and cones of objects in the categories $(\cC_1, \dots, \cC_m )$ is equivalent to $\cD$.
\end{definition}


Let $S_{\cD}$ be the Serre functor of $\cD$, then we have the following standard result, see e.g.~\cite[Section 3]{bayer2017stability}:

\begin{proposition}[{\cite[Section 3]{bayer2017stability}}] \label{serre s.o.d. proposition}
If $\cD = \langle \cD_1 , \cD_2 \rangle$ is a semiorthogonal decomposition, then $\cD = \langle S_{\cD}(\cD_2), \cD_1 \rangle = \langle \cD_2, S^{-1}_{\cD}(\cD_1) \rangle$ are also semiorthogonal decompositions.
\end{proposition}


\subsection{Mutations} \label{mutations subsection}

Let $\cC \sst \cD$ be an admissible triangulated subcategory. Then the \emph{left mutation functor} $\bL_{\cC}$ through $\cC$ is defined as the functor lying in the canonical functorial exact triangle 
\[  i i^! \ra \identity \ra \bL_{\cC}   \]
and the \emph{right mutation functor} $\bR_{\cC}$ through $\cC$ is defined similarly, by the triangle 
\[ \bR_{\cC} \ra \identity \ra i i^*  .  \]
When $E \in \D^b(X)$ is an exceptional object, and $F \in \D^b(X)$ is any object, the left mutation $\bL_E F$ fits into the triangle 
\[ E \otimes \RHom^\bullet(E, F) \ra F \ra \bL_E F , \]
and the right mutation $\bR_E F$ fits into the triangle
\[ \bR_E F \ra F \ra E \otimes \RHom^\bullet(F, E)^\vee  . \]

\begin{proposition}[{\cite[Lemma 2.6]{kuz:fractional-CY}}]\label{prop-serre-functor}
Let $\cD = \langle \cA, \cB \rangle $ be a semiorthogonal decomposition. Then 
\[ S_{\cB} = \bR_{\cA} \circ S_{\cD}  \, \, \, \text{ and } \, \, \, S_{\cA}^{-1} = \bL_{\cB} \circ S_{\cD}^{-1} . \]
\end{proposition}

\begin{lemma}[{\cite[Lemma 2.7]{kuznetsov2010derivedcubicfourfolds}}] \label{change_sod}
Let $\cD=\langle \cC_1, \cC_2,...,\cC_n \rangle$ be a semiorthogonal decomposition with all
components being admissible. Then for each $1\leq k\leq n-1$, there is a semiorthogonal decomposition
\[\cD=\langle \cC_1,..., \cC_{k-1}, \bL_{\cC_k}\cC_{k+1},\cC_k, \cC_{k+2} ...,\cC_n \rangle\]
and for each $2\leq k\leq n$ there is a semiorthogonal decomposition
\[\cD=\langle \cC_1,..., \cC_{k-2}, \cC_k,  \bL_{\cC_k}\cC_{k-1},\cC_{k+1} ...,\cC_n \rangle.\]
\end{lemma}

\section{Gushel--Mukai threefolds and their derived categories} \label{GM threefolds and their derived categories section}

Let $X$ be a prime Fano threefold of index $1$ and degree $H^3=10$, where $H$ is the ample generator of $\CaCl(X)$. Then $X$ is either a quadric section of a linear section of codimension $2$ of the Grassmannian $\mathrm{Gr}(2,5)$, in which case it is called an ordinary Gushel--Mukai (GM) threefold, or $X$ is a double cover of a degree $5$ and index $2$ Fano threefold $Y$ ramified in a quadric hypersurface, in which case it is called a special GM threefold. In the latter case, it has a natural involution $\tau \colon X\rightarrow X$ induced by the double cover $\pi \colon X\rightarrow Y$. By \cite{mukai:fano-3fold,bayer:vector-bundle-on-fano-threefold}, there exists a stable vector bundle $\mathcal{E}$ of rank $2$ with $c_1(\mathcal{E})=-H, c_2(\mathcal{E})=4L$ and $c_3(\cE)=0$, where $L$ is the class of a line on $X$. In addition, $\mathcal{E}$ is exceptional and $H^{\bullet}(X,\mathcal{E})=0$. In fact, $\cE$ is the pullback of the tautological bundle on the Grassmannian $\Gr(2,5)$. By \cite[Proposition 4.1]{debarre2012period}, $\cE$ is the unique stable sheaf with $c_1(\mathcal{E})=-H, c_2(\mathcal{E})=4L$ and $c_3(\cE)=0$.

Furthermore, there is a standard short exact sequence
\begin{equation}\label{tau-seq}
    0 \ra \cE \ra \oh_X^{\oplus 5} \ra \cQ \ra 0
\end{equation}
where $\cQ$ is the pull-back of the tautological quotient bundle on $\Gr(2,5)$ along the natural map $X\to \Gr(2,5)$. Since $\rk(\cE)=2$, we have $\cE(H)\cong \cE^{\vee}$.

\begin{definition} \label{kuznetsov components definition}
Let $X$ be a GM threefold.
\begin{itemize}
    \item The \emph{Kuznetsov component} of $X$ is defined as $\Ku(X) := \langle \cE,  \oh_X \rangle^\bot$. In particular, it fits into the semiorthogonal decomposition $\D^b(X) = \langle \Ku(X), \cE , \oh_X \rangle$;
    \item The \emph{alternative Kuznetsov component} of $X$ is defined as $\cA_X :=  \langle \oh_X , \cE^\vee \rangle^\bot$. In particular, it fits into the semiorthogonal decomposition $\D^b(X) = \langle \cA_X, \oh_X , \cE^\vee \rangle$.
\end{itemize}
\end{definition}

\begin{remark}\label{rmk-tau}
By \cite[Proposition 2.6]{kuznetsov2018derived}, there is a natural involutive autoequivalence functor $\tau_{\cA}:=S_{\cA_X}[-2]$ of $\cA_X$. When $X$ is special, it is induced by the natural involution $\tau$ on $X$ as $\tau_{\cA}=\tau^*|_{\cA_X}$.
\end{remark}

\begin{definition} \label{natural projection functor A_X}
The left adjoint to the inclusion $\cA_X \hookrightarrow \D^b(X)$ is given by $\pr := \bL_{\oh_X} \bL_{\cE^{\vee}} \colon \D^b(X) \ra \cA_X$. We call this the \emph{projection functor}.
\end{definition}

The analogous natural projection functor can be defined for $\Ku(X)$, and we denote it by $\pr':=\bL_{\cE}\bL_{\oh_X}$.

\subsection{Kuznetsov components} \label{kuznetsov components subsection}

Let $K_0(\cD)$ denote the Grothendieck group of a triangulated category $\cD$. We have the bilinear Euler form
\[ \chi(E, F) = \sum_{i \in \ZZ} (-1)^i \ext^i(E, F) \]
for $[E], [F] \in K_0(\cD)$. By the Hirzebruch--Riemann--Roch formula, it takes the following form on GM threefolds. We have \cite[p. 5]{kuznetsov2009derived} $\chi(u,v)=\chi_0(u^* \cap v)$ where $u \mapsto u^*$ is an involution of $\oplus_{i=0}^3 H^i(X, \mathbb{Q})$ given by multiplication with $(-1)^i$ on $H^{2i}(X, \mathbb{Q})$, and $\chi_0$ is given by
\[ \chi_0(x+yH +zL + wP) = x + \frac{17}{6} y + \frac{1}{2} z + w  , \] 
where $L$ is the class of lines and $P$ is the class of points. The \emph{numerical Grothendieck group} of $\cD$ is $\cN(\cD) = K_0(\cD) / \ker \chi $. 

\begin{lemma}[{\cite[p. 5]{kuznetsov2009derived}}]
The numerical Grothendieck group $\cN(\Ku(X))$ of the Kuznetsov component is a rank $2$ integral lattice generated by the basis elements $v = 1-3L + \frac{1}{2}P$ and $w= H - 6L + \frac{1}{6}P$. Using this basis, $\chi$ is given by the matrix
\[ \begin{pmatrix}  
-2 & -3 \\
-3 & -5
\end{pmatrix} . \]
\end{lemma}

\subsection{Alternative Kuznetsov components}


As in \cite[Proposition 3.9]{kuznetsov2009derived}, the following lemma follows from a straightforward computation.

\begin{lemma}
\label{lemma_Grothendieckgroup_kuz}
The numerical Grothendieck group of $\mathcal{A}_{X}$ is a rank 2 integral lattice with basis vectors $x=1-2L$ and $y=H-4L-\frac{5}{6}P$, and the Euler form with respect to the basis is 
\[
\begin{pmatrix}
-1 & -2 \\
-2 & -5
\end{pmatrix}.
\] 
\end{lemma}


\begin{remark}
It is straightforward to check that the $(-1)$-classes of $\mathcal{N}(\cA_X)$ are $x=1-2L$ and $2x-y=2-H+\frac{5}{6}P$, up to sign.
\end{remark}



Indeed, the Kuznetsov components from Subsection \ref{kuznetsov components subsection} and the alternative Kuznetsov components from this section are equivalent:

\begin{lemma} \label{equivalence of original and alternative kuznetsov component}
The original and alternative Kuznetsov components are equivalent. More precisely, there is an equivalence of categories $\Xi \colon \Ku(X) \xrightarrow{\sim} \cA_X$ given by $E \mapsto \bL_{\oh_X}(E\otimes\mathcal{O}_X(H))$, with inverse given by $F \mapsto (\bR_{\oh_X} F) \otimes\mathcal{O}_X(-H)$.
\end{lemma}

\begin{proof}
Using Lemma \ref{change_sod} and noting that $\cE\otimes \oh_X(H)\cong \cE^{\vee}$, we manipulate the semiorthogonal decomposition as follows:
\begin{align*}
    \D^b(X) &= \langle\Ku(X),\cE,\mathcal{O}_X\rangle \\
     &\simeq \langle\Ku(X)\otimes\mathcal{O}_X(H),\cE^{\vee},\mathcal{O}_X(H)\rangle \\
     &\simeq \langle\mathcal{O}_X,\Ku(X)\otimes\mathcal{O}_X(H),\cE^{\vee}\rangle \\
     &\simeq \langle \bL_{\oh_X}(\Ku(X)\otimes\mathcal{O}_X(H)),\oh_X,\cE^{\vee}\rangle .
\end{align*}
Now comparing with the definition of $\cA_X$, we get $\cA_X \simeq \bL_{\oh_X}(\Ku(X)\otimes\mathcal{O}_X(H))$ and the desired result follows. The reverse direction is similar.
\end{proof}

\section{Bridgeland stability conditions} \label{stability section}

In this section, we recall (weak) Bridgeland stability conditions on $\D^b(X)$, and the notions of tilt stability, double-tilt stability, and stability conditions induced on Kuznetsov components from weak stability conditions on $\D^b(X)$. We follow \cite[\S~2]{bayer2017stability}.

\subsection{Weak stability conditions}

Let $\cD$ be a triangulated category, and $K_0(\cD)$ its Grothendieck group. Fix a surjective morphism $v \colon K_0(\cD) \ra \Lambda$ to a finite rank lattice.

\begin{definition}
A \emph{stability condition (resp.~weak stability condition)} on $\cD$ is a pair $\sigma = (\cA, Z)$ where $\cA$ is the heart of a bounded t-structure on $\cD$, and $Z \colon \Lambda \ra \CC$ is a group homomorphism such that the following conditions hold:
\begin{enumerate}
    \item The composition $Z \circ v \colon K_0(\cA) \cong K_0(\cD) \ra \CC$ satisfies: for any $E\neq 0 \in \cD$ we have $\Im Z(E) \geq 0$ and if $\Im Z(E) = 0$ then $\Re Z(E) < 0$ (resp.~$\Re Z(E) \leq 0$). From now on, we write $Z(E)$ rather than $Z(v(E))$.
\end{enumerate}
We define a \emph{slope function} $\mu_\sigma$ for $\sigma$ using $Z$. For any $E \in \cA$, set
\[
\mu_\sigma(E) := \begin{cases}  - \frac{\Re Z(E)}{\Im Z(E)}, & \Im Z(E) > 0 \\
+ \infty , & \text{else}.
\end{cases}
\]
We say an object $0 \neq E \in \cA$ is $\sigma$-(semi)stable if $\mu_\sigma(F) < \mu_\sigma(E/F)$ (respectively $\mu_\sigma(F) \leq \mu_\sigma(E/F)$) for all proper subobjects $F \sst E$. 
\begin{enumerate}[resume]
    \item Any object $E \in \cA$ has a Harder--Narasimhan filtration in terms of $\sigma$-semistability defined above.
    \item There exists a quadratic form $Q$ on $\Lambda \otimes \mathbb{R}$ such that $Q|_{\ker Z}$ is negative definite, and $Q(E) \geq 0$ for all $\sigma$-semistable objects $E \in \cA$. This is known as the \emph{support property}.
\end{enumerate}
\end{definition}

\begin{definition}
Let $\sigma=(\cA, Z)$ be a stability condition on $\cD$. The \emph{phase} of a $\sigma$-semistable object $E\in \cA$ is
\[\phi(E) := \frac{1}{\pi} \mathrm{arg}(Z(E)) \in (0,1].\]
Specially, if $Z(E) = 0$ then $\phi(E) = 1$. If $F = E[n]$, then we define 
\[\phi(F) := \phi(E) + n\]

A \emph{slicing} $\cP$ of $\cD$ consists of full additive subcategories $\cP(\phi) \subset \cD$ for each $\phi \in \mathbb{R}$ satisfying
\begin{enumerate}[resume]
\item for $\phi\in (0,1]$, the subcategory $\cP(\phi)$ is given by the zero object and all $\sigma$-semistable objects whose phase is $\phi$;
\item for $\phi + n$ with $\phi\in (0,1]$ and $n\in \mathbb{Z}$, we set $\cP(\phi + n) := \cP(\phi)[n]$.
\end{enumerate}
\end{definition}

We will use both notations $\sigma = (\cA,Z)$ and $\sigma = (\cP,Z)$ for a stability condition $\sigma$ with heart $\cA = \cP((0,1])$ where $\cP$ is the slicing of $\sigma$.

We say $\sigma$ is a \emph{numerical stability condition} on $\cD$ if the surjective morphism $v\colon K_0(\cD)\to \Lambda$ factors through the natural surjection $K_0(\cD)\twoheadrightarrow \cN(\cD)$ (assuming $\cN(\cD)$ is well-defined).

Next, we recall two natural group actions on the set of stability conditions $\Stab(\cD)$.
\begin{enumerate}
    \item An element $\tilde{g} = (g, G)$ in the universal covering $\widetilde{\mathrm{GL}}^+(2,\mathbb{R})$ of the group $\mathrm{GL}^+(2,\mathbb{R})$ consists of an increasing function $g \colon \mathbb{R} \rightarrow \mathbb{R}$ such that $g(\phi+1) = g(\phi) + 1$ and a  matrix $G\in \mathrm{GL}^+(2,\mathbb{R})$ with $\det(G)>0$. It acts on the right on the stability manifold by $\sigma \cdot \tilde{g} := (G^{-1}\circ Z,\cP(g(\phi)))$ for any $\sigma = (\cP,Z)\in \Stab(\cD)$ (see \cite[Lemma 8.2]{bridgeland}).
    
    \item Let $\mathrm{Aut}_\Lambda(\cD)$ be the group of exact autoequivalences of $\cD$, whose action $\Phi_*$ on $K_0(\cD)$ is compatible with $v\colon K_0(\cD)\to \Lambda$. For $\Phi\in \mathrm{Aut}_\Lambda(\cD)$ and $\sigma=(\cP, Z) \in \mathrm{Stab}(\cD)$, we define a left action of the group of linear exact autoequivalences $\mathrm{Aut}_\Lambda(\cD)$ by $\Phi\cdot \sigma = (\Phi(\cP), Z\circ \Phi_*^{-1})$, where $\Phi_*$ is the automorphism of $K_0(\cD)$ induced by $\Phi$.
\end{enumerate}

\subsection{Tilt-stability}

Let $(X, H)$ be a polarised  smooth projective variety of dimension $n$ and $\sigma_H = (\Coh(X), Z_{H})$ be the standard weak stability condition on $\Coh(X)$ defined as
\[Z_H(E):=-H^{n-1}\ch_1(E)+\mathfrak{i}H^n\rk(E).\]
Its $\sigma_H$-stability coincides with classical $\mu_H$-stability (slope stability). Now for a fixed real number $\beta$, consider the following subcategories\footnote{The angle brackets here mean extension closure.} of $\Coh(X)$:
\begin{align*} 
\cT^{\beta}&= \langle  E \in \Coh(X) \mid E \text{ is } \sigma_H\text{-semistable with } \mu_{\sigma_H}(E) > \beta \rangle \\
\cF^{\beta}&= \langle E \in \Coh(X) \mid E \text{ is } \sigma_H\text{-semistable with } \mu_{\sigma_H}(E) \leq \beta  \rangle.
\end{align*}
Then it is a result of \cite{happel1996tilting} that the tilted heart $\Coh^{\beta}(X) := \langle \cT^{\beta}, \cF^{\beta}[1] \rangle$ is the heart of a bounded t-structure on $\Coh(X)$.

\begin{proposition}[{\cite{bayer2011bridgeland, bayer2016space}}] \label{BG_origin}
Let $\alpha>0$ and $\beta \in \mathbb{R}$. Then the pair $\sigma_{\alpha, \beta} = (\Coh^\beta(X) , Z_{\alpha, \beta})$ defines a weak stability condition on $\D^b(X)$, where
\[ Z_{\alpha, \beta}(E) = \frac{1}{2} \alpha^2 H^n \ch_0^\beta(E) - H^{n-2} \ch_2^\beta(E) + \mathfrak{i} H^{n-1} \ch_1^\beta(E). \]
The quadratic form $Q$ is given by the discriminant
\[ \Delta_H(E) = (H^{n-1} \ch_1(E))^2 - 2 H^n \ch_0(E) H^{n-2} \ch_2(E) .  \]
We denote the slope function by $\mu_{\alpha, \beta}:=\mu_{\sigma_{\alpha, \beta}}$.

\end{proposition}

The weak stability conditions $\sigma_{\alpha, \beta}$ constructed above are also known as \emph{tilt-stability} and the heart $\Coh^{\beta}(X)$ are called the \emph{tilted heart}.

Now pick a weak stability condition $\sigma_{\alpha, \beta}$. We define
\begin{align*} 
\cT^{0}_{\alpha, \beta}&= \langle  E \in \Coh^{\beta}(X) \mid E \text{ is } \sigma_{\alpha, \beta}\text{-semistable with } \mu_{\alpha, \beta}(E) > 0 \rangle \\
\cF^{0}_{\alpha, \beta}&= \langle E \in \Coh^{\beta}(X) \mid E \text{ is } \sigma_{\alpha, \beta}\text{-semistable with } \mu_{\alpha, \beta}(E) \leq 0  \rangle.
\end{align*}

Moreover, we ``rotate" the stability function $Z_{\alpha, \beta}$ by setting 
\[ Z_{\alpha, \beta}^0 := \frac{1}{\mathfrak{i}} Z_{\alpha, \beta}.\]
Then we have the following result:

\begin{proposition}[{\cite[Proposition 2.15]{bayer2017stability}}]
The pair $\sigma^0_{\alpha, \beta}=(\Coh^0_{\alpha, \beta}(X)=\langle \cT^{0}_{\alpha, \beta}, \cF^{0}_{\alpha, \beta}[1]  \rangle, Z_{\alpha, \beta}^0)$ defines a weak stability condition on $\D^b(X)$. We denote the slope function by $\mu^0_{\alpha, \beta}:=\mu_{\sigma^0_{\alpha, \beta}}$.

\end{proposition}

We now state a useful lemma that relates 2-Gieseker-stability (see \cite[Definition 4.3]{bayer2020desingularization}) and tilt-stability. 

\begin{lemma}[{\cite[Lemma 2.7]{bayer2016space}, \cite[Proposition 4.8, 4.9]{bayer2020desingularization}}] \label{bms lemma 2.7}
Let $E\in \D^b(X)$.

\begin{enumerate}
    \item  Let $\beta<\mu(E)$. Then $E\in \Coh^{\beta}(X)$ is $\sigma_{\alpha, \beta}$-(semi)stable for $\alpha\gg 0$ if and only if $E\in \Coh(X)$ and $E$ is 2-Gieseker-(semi)stable.
    
    \item If $E\in \Coh^{\beta}(X)$ is $\sigma_{\alpha, \beta}$-semistable for $\beta\geq \mu(E)$ and $\alpha\gg 0$, then $\cH^{-1}(E)$ is a torsion free $\mu$-semistable sheaf and $\cH^0(E)$ is supported in dimension not greater than one. If $\beta>\mu(E)$ and $\alpha>0$, then $\cH^{-1}(E)$ is also reflexive.
\end{enumerate}

\end{lemma}

\subsection{Stronger BG inequalities}

In this subsection, we state stronger Bogomolov--Gieseker (BG) style inequalities, which hold for tilt-semistable objects. These will be useful later on for ruling out potential walls for tilt-stability of objects in $\D^b(X)$. The first is a stronger version of Proposition \ref{BG_origin},  which was proved by Chunyi Li in \cite[Proposition 3.2]{li2016stability} for Fano threefolds of Picard number one.

\begin{lemma}[Stronger BG I] \label{Li BG}
Let $X$ be an index $1$ prime Fano threefold with degree $d$, and $E\in \D^b(X)$ a $\sigma_{\alpha, \beta}$-stable object where $\alpha >0$. Let $k:=\lfloor \mu(E)\rfloor$. Then we have:
\[\frac{H \cdot \ch_2(E) }{H^3 \cdot \ch_0(E)} \leq \max \left\{ k\mu_H(E)-\frac{k^2}{2}, \frac{1}{2} \mu_H(E)^2 - \frac{3}{4d}, (k+1)\mu_H(E)-\frac{(k+1)^2}{2} \right\} .\]

Moreover, if the equality holds, then $E$ has rank one or two.
\end{lemma}

The second is due to Naoki Koseki and Chunyi Li. It is based on \cite[Lemma 4.2, Theorem 4.3]{koseki2020bogomolov}. Chunyi Li also sent us a similar inequality from his upcoming paper \cite{li2021strongBG}.

\begin{theorem}[Stronger BG II] \label{Naoki BG}
Let $X$ be an index $1$ Fano threefold of degree $d$, and $E \in \Coh^0(X)$ be a $\sigma_{\alpha, 0}$-semistable object for some $\alpha >0$ with $\vert \mu_H(E)\vert \in [0,1]$ and $\rk (E) \geq 2$. Then 
\[  \frac{H \cdot \ch_2(E) }{H^3 \cdot \ch_0(E)} \leq \max \left\{ \frac{1}{2} \mu_H(E)^2 - \frac{3}{4d}, \mu_H(E)^2 - \frac{1}{2}\vert\mu_H(E)\vert \right\} . \]
\end{theorem}

Before we prove Theorem \ref{Naoki BG}, we first state an easy lemma.

\begin{lemma} \label{Lem:K3 surface inequality}
Let $S$ be a K3 surface of degree $d$ and $H_S$ the ample polarisation. Let $E$ be a $\mu_{H_S}$-semistable sheaf in $\D^b(S)$ with $\rk (E) \geq 2$. Then 
\[ \frac{\ch_2(E)}{H_S^2 \cdot \rk(E)} \leq \frac{1}{2} \mu_{H_S}(E)^2 - \frac{3}{4d} . \]
\end{lemma}

\begin{proof}
Let $v(E)$ be the Mukai vector of $E$. We have 
\begin{align*}
    v(E)^2 &= H_S^2 \cdot  \ch_1(E)^2 - 2 \rk(E)^2 -2 \rk(E) \cdot \ch_2(E) \\
    &\geq -2 \geq -\frac{1}{2}\rk(E)^2 .
\end{align*}
Dividing through by $\rk(E)^2$ and rearranging, we get
\begin{align*}
    \frac{\ch_2(E)}{\rk(E)} &\leq \frac{1}{2}\mu_{H_S}(E)^2 H_S^2 - \frac{3}{4}
\end{align*}
as required.
\end{proof}

\begin{proof}[Proof of Theorem \ref{Naoki BG}]
Let $f \colon [0,1] \ra \mathbb{R}$ be defined as
\[ f(t) := \max \left\{ \frac{1}{2}t^2 - \frac{3}{4d}, t^2 - \frac{1}{2}t \right\} \]
Note that $f$ is star-shaped \cite[Definition 3.2]{koseki2020bogomolov} and satisfies $f(0)=0$ and $f(1)=1/2$ as well as 
\[ t^2 - \frac{1}{2}t \leq f(t) \leq \frac{1}{2} t^2  \]
for all $t \in [0,1]$. We now follow the strategy of proof in \cite[Theorem 4.3]{koseki2020bogomolov}. Assume for a contradiction that there is an $E \in \D^b(X)$ such that the inequality in the statement of Theorem \ref{Naoki BG} is not true. Then conditions (a) and (b) in \cite[Lemma 3.3]{koseki2020bogomolov} are satisfied for $f$. Then by \emph{loc. cit.}, the restriction $E|_{S_d}$ where $S_d$ is a general hyperplane section of $X_d$ is $\mu_{H_{S_d}}$-semistable. Also note that $\mu_{H_{S_d}}(E|_{S_d}) = \mu_H(E)$ and $S_d$ is a smooth K3 surface. But then by assumption
\[ \frac{\ch_2(E|_{S_d})}{H_{S_d}^2 \cdot \rk(E|_{S_d})} > \frac{1}{2} \mu_{H_{S_d}}(E)^2 - \frac{3}{4d} \]
which contradicts Proposition \ref{Lem:K3 surface inequality}, so the assumption is false and the result follows.
\end{proof}

\subsection{Stability conditions on the Kuznetsov component of a GM threefold}\label{section_double_tilted_stability_condition}

Proposition 5.1 in \cite{bayer2017stability} gives a criterion for checking when weak stability conditions on a triangulated category can be used to induce stability conditions on a subcategory. Each of the criteria of this proposition can be checked for $\cA_X \sst \D^b(X)$ to give stability conditions on $\cA_X$. 

More precisely, let $\cA(\alpha, \beta) = \Coh^0_{\alpha, \beta}(X) \cap \cA_X$ and $Z(\alpha, \beta) = Z_{\alpha, \beta}^0|_{\cA_X}$. Furthermore, if we take suitable $(\alpha, \beta)$, 
by {\cite[Theorem 6.9]{bayer2017stability}} and \cite[Proposition 3.2]{Pertusi2021serreinv} we have:

\begin{theorem}
\label{Theorem_stabilitycondition_exists}
Let $X$ be a GM threefold. Then  $\sigma(\alpha,\beta)$ is a stability condition on $\mathcal{A}_X$ for all $(\alpha, \beta)\in V$, where
\[V:=\{(\alpha, \beta)\colon -\frac{1}{10}<\beta<0, 0<\alpha<-\beta\}.\]
\end{theorem}

Now we introduce a special class of stability condition, which will play a central role in our paper.

\begin{definition}
Let $\sigma$ be a stability condition on a triangulated category $\cD$. It is called \emph{Serre-invariant} if $S_{\cD} \cdot \sigma=\sigma \cdot g$ for some $g\in\widetilde{\mathrm{GL}}^+(2,\mathbb{R})$, where $S_{\cD}$ is the Serre functor of $\cD$.
\end{definition}

We recall a recent result proved in \cite{Pertusi2021serreinv}. 
\begin{theorem}\label{Serre-invariant-GM}
Let $X$ be a GM threefold and $\sigma$ (resp.~$\sigma'$) be a stability condition on $\mathcal{K}u(X)$ (resp.~$\mathcal{A}_X$) defined by \cite{bayer2017stability}. Then $\sigma$ (resp.~$\sigma'$) is Serre-invariant. 
\end{theorem}


\begin{proposition}\label{ext23-stable}
Let $X$ be a GM threefold and $E$ a non-zero object in $\cA_X$ such that $\mathrm{ext}^1(E, E)\leq 3$ and $-\chi(E, E)$ is not a perfect square. Then $E$ is $\sigma$-stable for every Serre-invariant stability condition $\sigma$ on $\cA_X$. 
\end{proposition}

\begin{proof}
The proof is the same as in \cite[Lemma 9.12]{Zhang2020Kuzconjecture}. We omit the details. 
\end{proof}

\section{Projection of \texorpdfstring{$\cE$}{the vector bundle} into \texorpdfstring{$\Ku(X)$}{the Kuznetsov component}} \label{gluing data machine section}

In this section, we consider the object that results from projecting the vector bundle $\cE$ into $\Ku(X)$, and its stability in $\Ku(X)$. We start with a lemma.

\begin{lemma}\label{cohmo-bdl}
Let $X$ be a GM threefold.

\begin{enumerate}
    \item $\RHom^\bullet(\cQ(-H), \cE)=\RHom^\bullet(\cE, \cQ^{\vee})=k^2$ when $X$ is ordinary.
    
    \item $\RHom^\bullet(\cQ(-H), \cE)=\RHom^\bullet(\cE, \cQ^{\vee})=k^3\oplus k[-1]$ when $X$ is special.
    
    \item $\RHom^\bullet(\cE, \cQ(-H))=k[-2]$.

    \item $\RHomb(\cQ^{\vee}, \cE)=\RHomb(\cE^{\vee}, \cQ)=k[-2]$.
    
\end{enumerate}

\end{lemma}

\begin{proof}
When $X$ is ordinary, (1) and (2) follow from the Koszul resolution of $X\subset \mathrm{Gr}(2,5)$ and the Borel--Weil--Bott Theorem. When $X$ is special with the double cover $\pi\colon X\to Y$, note that $\pi_*\oh_X=\oh_Y\oplus \oh_Y(-1)$. Then the (1) and (2) follow from the projection formula and \cite[Lemma 2.14, Proposition 2.15]{sanna2014rational}. And applying $\Hom(-,\cE)$ to \eqref{tau-seq} and using Serre duality, we get $\RHom^\bullet(\cE, \cQ(-H))=\RHomb(\cQ, \cE)^{\vee}[-3]=k[-2]$, which proves (3). Finally, (4) follows from applying $\Hom(-,\cE)$ to \eqref{tau-seq} and using Serre duality and $\RHomb(\cE,\cE)=k$.
\end{proof}

\subsection{The projection of \texorpdfstring{$\cE$}{the vector bundle} into \texorpdfstring{$\Ku(X)$}{the Kuznetsov component}} \label{Ku(X) gluing data subsection}

Let $\pi := \bR_{\oh_X(-H)}\bR_{\cE(-H)} \colon \D^b(X) \ra \Ku(X)$ be the right adjoint to the inclusion $\Ku(X) \hookrightarrow \D^b(X)$. Here $\Ku(X)=\langle \cE, \oh_X \rangle^{\perp}$ is the original Kuznetsov component.

\begin{lemma} \label{cohomology objects of pi lemma}
The projection object $\pi(\E)$ is the unique object that fits into a non-trivial exact triangle
\begin{equation}\label{piE-seq}
    \cQ(-H)[1]\to \pi(\cE)\to \cE.
\end{equation}
\end{lemma}

\begin{proof}
By Serre duality, we have $\RHom^\bullet(\cE, \cE(-H))=\RHom^\bullet(\cE, \cE)^{\vee}[-3]=k[-3]$. Then we have an exact triangle $\bR_{\cE(-H)}\cE\to \cE\to \cE(-H)[3]$. And by \eqref{tau-seq}, we see $\bR_{\oh_X(-H)}\cE(-H)=\cQ(-H)[-1]$. Thus, from $\bR_{\oh_X(-H)}\cE=\cE$ we obtain the triangle \eqref{piE-seq}. It is non-trivial since $\pi(\cE)\in \Ku(X)$, so $\cE$ cannot be a direct summand of $\pi(\cE)$. Finally, the uniqueness follows from Lemma \ref{cohmo-bdl} (4).
\end{proof}

\begin{lemma}\label{lemma_stability_projection}
Let $X$ be a GM threefold. Then we have

\begin{itemize}
    \item $\RHom^\bullet(\pi(\cE), \pi(\cE))=k\oplus k^2[-1]$ when $X$ is ordinary.
    
    \item $\RHom^\bullet(\pi(\cE), \pi(\cE))=k\oplus k^3[-1]\oplus k[-2]$ when $X$ is special.
\end{itemize}

Hence 
$\pi(\cE)$ is stable with respect to every Serre-invariant stability condition on $\Ku(X)$.
\end{lemma}

\begin{proof}
The first statement follows from applying $\Hom(-, \cE)$ to triangle (\ref{piE-seq}) and Lemma \ref{cohmo-bdl}, and also the fact that $\RHom^\bullet(\pi(\cE), \pi(\cE))=\RHom^\bullet(\pi(\cE), \cE)$ which is by adjunction. The last statement follows from Lemma \ref{ext23-stable}.
\end{proof}

\subsection{The analogous projection object for \texorpdfstring{$\cA_X$}{the alternative Kuznetsov component}}
In this subsection, we state and prove the analogous results as in Subsection \ref{Ku(X) gluing data subsection}, except for $\cA_X$ instead of $\Ku(X)$. Let $\pi' := \bR_{\cE}\bR_{\oh_X(-H)} \colon \D^b(X) \ra \cA_X$ be the right adjoint to the inclusion $\cA_X \hookrightarrow \D^b(X)$.

\begin{lemma} \label{pi' cohomology objects lemma}
The projection object $\pi'(\cQ^\vee)$ is the unique object fits into a non-trivial exact triangle
\begin{equation}\label{piQ-seq}
    \cE[1]\to \pi'(\cQ^{\vee})\to \cQ^{\vee}.
\end{equation}
\end{lemma}

\begin{proof}
The proof is completely analogous to the proof of Lemma \ref{cohomology objects of pi lemma}. By Serre duality, we have the vanishing $\RHomb(\cQ^{\vee}, \oh_X(-H))=\RHomb(\oh_X, \cQ^{\vee})^{\vee}=0$. Thus $\pi'(\cQ^{\vee})=\bR_{\cE}\cQ^{\vee}$. Then the result follows from Lemma \ref{cohmo-bdl} (4). Finally, the uniqueness also follows from Lemma \ref{cohmo-bdl} (4) and \eqref{piQ-seq} is non-trivial since $\RHomb(\cQ^{\vee}, \pi'(\cQ^{\vee}))=0$.
\end{proof}

\begin{remark}\leavevmode
Later in Section \ref{refined categorical Torelli section}, we will see that we have $\pi'(\cQ^\vee) \cong \pr(I_C)[1]$ where $C \sst X$ is a conic such that $I_C\not\in\cA_X$. 
\end{remark}

\begin{lemma} \label{AX_piQ_stable}
Let $X$ be a GM threefold. Then \begin{itemize}
    \item $\RHom^\bullet(\pi'(\cQ^{\vee}), \pi'(\cQ^{\vee}))=k\oplus k^2[-1]$ when $X$ is ordinary.
    
    \item $\RHom^\bullet(\pi'(\cQ^{\vee}), \pi'(\cQ^{\vee}))=k\oplus k^3[-1]\oplus k[-2]$ when $X$ is special.
\end{itemize}

Hence $\pi'(\cQ^{\vee})$ is stable with respect to every Serre-invariant stability condition on $\cA_X$.
\end{lemma}

\begin{proof}
It is not hard to check that $\Xi(\pi(\cE)) \cong \pi'(\cQ^\vee)[1]$, where $\Xi$ is the equivalence $\Ku(X) \simeq \cA_X$ in Lemma \ref{equivalence of original and alternative kuznetsov component}. Then the result follows from Lemma \ref{lemma_stability_projection}.
\end{proof}

\section{Conics on GM threefolds} \label{hilbert schemes section}

In this section, we collect some useful results regarding the birational geometry of GM threefolds and their Hilbert schemes of conics. The results in this section are all from \cite{debarre2012period},  \cite{logachev2012fano}, and \cite{iliev1994fano}.

Recall that a conic means a closed subscheme $C \sst X$ with Hilbert polynomial $p_C(t)=1+2t$, and a line means a closed subscheme $L \sst X$ with Hilbert polynomial $p_L(t)=1+t$. Denote their Hilbert schemes by $\cC(X)$ and $\Gamma(X)$, respectively. 

\subsection{Conics on ordinary GM threefolds}

Let $X$ be an ordinary GM threefold. Recall that it is a quadric section of a linear section of codimension $2$ of the Grassmannian $\Gr(2,5)=\Gr(2, V_5)$, where $V_5$ is a $5$-dimensional complex vector space. Let $V_i$ be an $i$-dimensional vector subspace of $V_5$. There are two types of $2$-planes in $\Gr(2,5)$; $\sigma$\emph{-planes} are given set-theoretically as $\{ [V_2] \mid V_1 \sst V_2 \sst V_4  \}$, and \emph{$\rho$-planes} are given by $\{ [V_2] \mid V_2 \sst V_3 \}$.

\begin{remark}
In \cite[Section 3.1]{debarre2012period}, the $\sigma$-planes and $\rho$-planes are called $\alpha$-planes and $\beta$-planes, respectively.
\end{remark}

By \cite[Section 3.1]{debarre2012period} and \cite[Section 3.1]{iliev:fano-manifold-degree-ten}, we have the following classification of conics on $X$.

\begin{definition}[{\cite[p. 5]{debarre2012period}}] \leavevmode
\begin{itemize}
    \item A conic $C \sst X$ is called a \emph{$\tau$-conic} if the $2$-plane $\langle C \rangle$ is not contained in $\Gr(2, V_5)$, there is a unique $V_4 \sst V_5$ such that $C \sst \Gr(2, V_4)$, the conic $C$ is reduced and if it is smooth, the union of corresponding lines in $\mathbb{P}(V_5)$ is a smooth quadric surface in $\mathbb{P}(V_4)$.
    \item A conic $C \sst X$ is called a \emph{$\sigma$-conic} if the $2$-plane $\langle C \rangle$ spanned by $C$ is an $\sigma$-plane, and if there is a unique hyperplane $V_4 \sst V_5$ such that $C \sst \Gr(2, V_4)$ and the union of the corresponding lines in $\PP(V_5)$ is a quadric cone in $\PP(V_4)$.
    \item A conic $C \sst X$ is called a \emph{$\rho$-conic} if the $2$-plane $\langle C \rangle$ spanned by $C$ is a $\rho$-plane, and the union of corresponding lines in $\PP(V_5)$ is this $2$-plane.
\end{itemize}
\end{definition}

The following lemma is very useful for computations:

\begin{lemma}\label{IC-class}
Let $X$ be an ordinary GM threefold and $C$ be a conic on $X$.

\begin{enumerate}
    \item If $C$ is a $\tau$-conic, then we have $\RHom^\bullet(\cE, I_C)=k$ and $\RHom^\bullet(\cE^{\vee}, I_C)=0$.
    
    \item If $C$ is a $\rho$-conic, then we have $\RHom^\bullet(\cE, I_C)=k^2\oplus k[-1]$ and $\RHom^\bullet(\cE^{\vee}, I_C)=0$.
    
    \item If $C$ is a $\sigma$-conic, then we have $\RHom^\bullet(\cE, I_C)=k$ and $\RHom^\bullet(\cE^{\vee}, I_C)=k[-1]\oplus k[-2]$.
\end{enumerate}

\end{lemma}

\begin{proof}
Note that if $\Hom(\cE, I_C)=k^a$, then $C\subset \mathrm{Gr}(2, 5-a)\cap X$. Since for any conic $C$, there is some $V_4$ such that $C\subset \mathrm{Gr}(2,V_4)$, then we have $\hom(\cE, I_C)\geq 1$ for any conic $C$. Now if $\hom(\cE, I_C)\geq 2$, we know that $C$ is contained in a $\rho$-plane $\mathrm{Gr}(2,3)$. Since $\langle C\rangle$ is not in $\mathrm{Gr}(2,5)$ for a $\tau$-conic $C$, and $\langle C \rangle$ is a $\sigma$-plane $\{V_2|V_1\subset V_2\subset V_4\}$ for a $\sigma$-conic, for these two types of conics we have $\Hom(\cE, I_C)=k$. Also, for a $\rho$-conic $C$, since $\langle C \rangle=\mathrm{Gr}(2,3)$, we have $\hom(\cE, I_C)\geq 2$. But if $\hom(\cE, I_C)\geq 3$, we know that $C\subset \mathrm{Gr}(2,2)$ which is impossible. Hence for a $\rho$-conic $C$ we have $\Hom(\cE, I_C)=k^2$. Now the result for Ext groups follows from applying $\Hom(\cE, -)$ to the short exact sequence $0\to I_C\to \oh_X\to \oh_C\to 0$ and $\chi(\cE, I_C)=1$.

First by stability and Serre duality, we have $\Hom(\cE^{\vee}, I_C)=\Ext^3(\cE^{\vee}, I_C)=0$. From $\chi(\cE^{\vee}, I_C)=0$, we only need to compute $\Ext^1(\cE^{\vee}, I_C)$. Since $\RHom^\bullet(\oh_X, I_C)=0$, applying $\Hom(-, I_C)$ to the tautological sequence, we have $\Hom(\cQ^{\vee}, I_C)=\Ext^1(\cE^{\vee}, I_C)$. Note that if $\Hom(\cQ^{\vee}, I_C)=k^a$, then $C\subset \mathrm{Gr}(2-a, 5-a)\cap X$. Thus we have $\hom(\cQ^{\vee}, I_C)\leq 1$ for any conic $C$. And since $\hom(\cQ^{\vee}, I_C)=1$ if and only if $C$ is contained in the zero locus of a global section of $\cQ$, which is a $\sigma$-3-plane in $\mathrm{Gr}(2, 5)$, we know that $\Hom(\cQ^{\vee}, I_C)=0$ for $C$ of type $\tau$ or $\rho$, and $\Hom(\cQ^{\vee}, I_C)=k$ for a $\sigma$-conic. Then the result follows.
\end{proof}

Now we recall some properties of the Fano surface of conics $\cC(X)$.

\begin{theorem}[{\cite{logachev2012fano}, \cite{debarre2012period}}] \label{ordinary-CX}
Let $X$ be an ordinary GM threefold. Then $\cC(X)$ is an irreducible projective surface. If $X$ is furthermore general, then $\cC(X)$ is smooth.
\end{theorem}

It is a fact that there is a unique $\rho$-conic on $X$, and there is a curve $L_{\sigma}\subset \cC(X)$ parameterise all $\sigma$-conics on $X$ (cf.~\cite[Section 5.1]{debarre2012period}), and we denote it by $c_X$. Furthermore, we have the following result which is a corollary of Logachev's Tangent Bundle Theorem (\cite[Section 4]{logachev2012fano}). 

\begin{lemma}[{\cite[p. 16]{debarre2012period}}] \label{minimal surface of S(X) lemma}
The only rational curve in $\cC(X)$ is $L_\sigma$. Furthermore, there exists a surface $\cC_m(X)$ and a map $\cC(X) \ra \cC_m(X)$ which contracts $L_\sigma$ to a point $[\pi]$. If $X$ is general, then $\cC_m(X)$ is the minimal surface of $\cC(X)$.
\end{lemma}

\begin{theorem}[{\cite[Section 5.2]{debarre2012period}}] \label{classical_invo_conic}
Let $X$ be a general ordinary GM threefold. Then there is a natural involution $\iota$ on $\cC_m(X)$, switching the points $[c_X]$ and $[\pi]$.
\end{theorem}

Another important result that we require is Logachev's Reconstruction Theorem. This was originally proved in \cite[Theorem 7.7]{logachev2012fano}, and then reproved later in \cite[Theorem 9.1]{debarre2012period}.

\begin{theorem}[Logachev's Reconstruction Theorem] \label{logachev's reconstruction theorem}
Let $X$ and $X'$ be general ordinary GM threefolds. If $\cC(X) \cong \cC(X')$, then $X \cong X'$.
\end{theorem}

\subsection{Conic and line transforms}

For this section, we follow \cite[Section 6.1]{debarre2012period}. Let $X$ be a general ordinary GM threefold, and let $C$ be a conic. Then in \cite[\S~6.1, Theorem 6.4]{debarre2012period}, the authors construct a new GM threefold $X_C$ and a birational map $\psi_C\colon X\dashrightarrow X_C$, called the \emph{conic transform}. Similarly, for any line $L\subset X$, a new GM threefold $X_L$ and a birational morphism $\psi_L\colon X\to X_L$ are constructed in \cite[Section 6.2]{debarre2012period}, called the \emph{line transform}.

Note that in \cite{debarre2015gushel}, such an $X_C$ is called the \emph{period partner} of $X$, and the line transforms are called the \emph{period duals}. We now list some important results about conic and line transforms below.

\begin{theorem}[{\cite[Theorem 6.4]{debarre2012period}}] \label{DIM theorem 6.4}
Let $X$ be a general ordinary GM threefold, and let $C \sst X$ be a conic. Then $\cC(X_C)$ is isomorphic to $\cC_m(X)$ blown up at the point $[C] \in \cC_m(X)$, where $\cC_m(X)$ is the minimal surface of $\cC(X)$.
\end{theorem}

\begin{proposition}[{\cite[Theorem 6.4, Remark 7.2]{debarre2012period}}] \label{DIM9.2}
Let $X$ be a general ordinary GM threefold. Then the isomorphism classes of conic transforms of $X$ are parametrized by the surface $\cC_m(X)/\iota$.
\end{proposition}

\begin{theorem}[{\cite[Theorem 1.6]{kuznetsov2019categorical}}] \label{dual_conj}
Let $X$ be a general ordinary GM threefold. Then the Kuznetsov components of all conic transforms and line transforms of $X$ are equivalent to $\cA_X$.
\end{theorem}

\subsection{Conics on special GM threefolds} Let $X$ be a special GM threefold. Recall that $X$ is a double cover $X\to Y$ of a degree five del Pezzo threefold $Y$ with branch locus a quadric hypersurface $\cB\subset Y$. When $X$ is general, $\cB$ is a smooth K3 surface of Picard number $1$ and degree $10$. Recall that $Y$ is a codimension $3$ linear section of $\mathrm{Gr}(2, 5)$. Let $\mathcal{V}$ be the tautological quotient bundle on $Y$. We recall some properties of $\cC(X)$ from \cite{iliev1994fano}.

\begin{theorem}[{\cite{iliev1994fano}}] \label{special_surface}
Let $X$ be a special GM threefold. Then $\cC(X)$ has two components $\cC_1$ and $\cC_2$. One of the components $\cC_2\cong \Sigma(Y)\cong \mathbb{P}^2$ parametrizes the preimage of lines on $Y$. Moreover, when $X$ is general, $\cC(X)$ is smooth away from $\cC_1\cap \cC_2$.
\end{theorem}

The following lemma will be useful in computations; it is similar to Lemma \ref{IC-class}.

\begin{lemma} \label{IC-class-sp}
Let $X$ be a special GM threefold and $C$ a conic on $X$. Then $\RHom^\bullet(\cE^{\vee}, I_C)\neq 0$ if and only if $C$ is the preimage of a line on $Y$. In this case $\RHom^\bullet(\cE^{\vee}, I_C)=k[-1]\oplus k[-2]$, and such a family of conics is parametrized by the Hilbert scheme of lines $\Sigma(Y)\cong \mathbb{P}^2$ on $Y$.
\end{lemma}

\begin{proof}
The proof is almost the same as the second part of the proof of Lemma~\ref{IC-class}. The same argument shows that $\RHom^\bullet(\cE^{\vee}, I_C)\neq 0$ if and only if $\Hom(\cQ^{\vee}, I_C)\neq 0$. The image of a non-trivial map $\cQ^{\vee}\to  I_C$ is the ideal sheaf of the zero locus of a section $s$ of $\cQ$, which is the preimage of the zero locus of a section of $\mathcal{V}$. By \cite[Lemma 2.18]{sanna2014rational}, the zero locus of a section of $\mathcal{V}$ is either a line or a point. Thus the zero locus of a section of $\cQ$ is either the preimage of a line on $Y$ which is a conic on $X$, or a zero-dimensional closed subscheme of length two. But this zero locus contains a conic $C\subset X$, so $C=Z(s)$ is the preimage of a line on $Y$ and the map $\cQ^{\vee}\to I_C$ is surjective. In particular,  such conics are exactly the preimages of lines on $Y$ and are parametrized by $\Sigma(Y)\cong\mathbb{P}^2$.
\end{proof}

\section{Conics and Bridgeland moduli spaces} \label{refined categorical Torelli section}

In this section, we study the moduli space of $\sigma$-stable objects of the $(-1)$-class $-x$ in the alternative Kuznetsov component $\mathcal{A}_X$ of a GM threefold $X$ and its relation to $\cC(X)$. Our main result in this section is Theorem \ref{theorem_irreduciblecomponent_conics}, which realizes the Bridgeland moduli space as a contraction of $\cC(X)$. 

First, we study those conics $C$ such that $I_C\notin \cA_X$.

\begin{proposition}\label{prop_contraction_P1_OGM}
Let $C\subset X$ be a conic on a GM threefold $X$. Then $I_C\not\in\mathcal{A}_X$ if and only if

\begin{enumerate}
    \item $C$ is a $\sigma$-conic when $X$ is ordinary. In particular, such a family of conics is parametrized by the line $L_{\sigma}$. 
    
    \item $C$ is the preimage of a line on $Y$ when $X$ is special. In particular, such a family of conics is parametrized by the Hilbert scheme of lines $\Sigma(Y)\cong \mathbb{P}^2$ on $Y$.
\end{enumerate}

Moreover, we have an exact sequence
\[0\to \cE\to \cQ^{\vee}\to I_C\to 0.\]

\end{proposition}

\begin{proof}
Note that $I_C\notin \cA_X$ if and only if $\RHomb(\cE^{\vee}, I_C)\neq 0$. When $X$ is ordinary, (1) follows from Lemma \ref{IC-class}. When $X$ is special, we deduce (2) from Lemma \ref{IC-class-sp}. Note that since $I_C\notin \cA_X$, we have $\Hom(\cQ^{\vee}, I_C)\neq 0$. The non-trivial map $\cQ^{\vee}\to I_C$ is surjective by the arguments in Lemma \ref{IC-class} and \ref{IC-class-sp}. Note that by the stability of $\cQ^{\vee}$, the kernel of $\cQ^{\vee}\twoheadrightarrow I_C$ is $\mu$-stable with the same Chern character as $\cE$, hence we have  $\ker(\cQ^{\vee}\twoheadrightarrow I_C)\cong \cE$ by \cite[Proposition 4.1]{debarre2012period}.
\end{proof}

\begin{proposition}\label{prop_conic_inKuz}
Let $X$ be a GM threefold and $C\subset X$ a conic on $X$. If $I_C\not\in\mathcal{A}_X$, then we have the exact triangle
$$\mathcal{E}[1]\rightarrow\mathrm{pr}(I_C)\rightarrow \cQ^{\vee}$$
and $\pr(I_C)\cong \pi'(\cQ^{\vee})$
\end{proposition}

\begin{proof}
By Proposition~\ref{prop_contraction_P1_OGM}, $I_C$ fits into the short exact sequence
$$0\rightarrow\cE\rightarrow\cQ^{\vee}\rightarrow I_C\rightarrow 0.$$
Applying the projection functor to this exact sequence, and note that applying the functor $\pr$ to the dual exact sequence of \eqref{tau-seq} gives
$\mathrm{pr}(\cQ^{\vee})=0$. Then we have $\mathrm{pr}(I_C)\cong\mathrm{pr}(\cE)[1]$. Now we compute the projection $\mathrm{pr}(\cE)$. Since $\mathrm{RHom}^\bullet(\cE^{\vee},\cE) \cong k[-3]$, we get an exact triangle $\cE^{\vee}[-3]\rightarrow\cE\rightarrow\bL_{\cE^{\vee}}\cE$. Now applying $\bL_{\mathcal{O}_X}$ to this triangle and using $\bL_{\oh_X}\cE^{\vee}=\cQ^{\vee}[1]$, we get $$\cQ^{\vee}[-2]\rightarrow\cE\rightarrow\mathrm{pr}(\cE).$$
Therefore we obtain the triangle
$$\cE[1]\rightarrow\mathrm{pr}(\cE)[1]\rightarrow\cQ^{\vee}$$
and the desired result follows from Lemma \ref{pi' cohomology objects lemma}. 
\end{proof}

Now the following two results follow from Proposition \ref{prop_conic_inKuz} and  Lemma \ref{AX_piQ_stable}.

\begin{lemma}
\label{lemma_ideal_projection_smooth} Let $X$ be a GM threefold. If $C\subset X$ is a conic such that $I_C\not\in\cA_X$, then
\begin{itemize}
    \item $\RHom^\bullet(\pr(I_C), \pr(I_C))=k\oplus k^2[-1]$ when $X$ is ordinary.
    
    \item $\RHom^\bullet(\pr(I_C), \pr(I_C))=k\oplus k^3[-1]\oplus k[-2]$ when $X$ is special.
\end{itemize}

\end{lemma}

\begin{lemma} \label{pr(IC)[1] when h0(EC)=1 stability}
Let $X$ be a GM threefold. If $I_C\not\in\mathcal{A}_X$, the projection $\pr(I_C)[1]$ is stable with respect to every Serre-invariant stability condition on $\cA_X$.
\end{lemma}

When $I_C\in \cA_X$, we cannot use Proposition \ref{ext23-stable} to prove the Bridgeland stability of $I_C$, since $\cC(X)$ can be singular and $\Ext^1(I_C, I_C)$ may have large dimension. Instead, we use a wall-crossing argument and the uniqueness of Serre-invariant stability conditions (Theorem \ref{all_in_one_orbit}).

\begin{lemma} \label{beta=0}
Let $X$ be a GM threefold. Let $F$ be an object with $\ch_{\leq 2}(F)=(1,0,-2L)$. Then there are no walls for $F$ in the range $-\frac{1}{2}\leq \beta<0$ and $\alpha>0$.
\end{lemma}

\begin{proof}
Recall that by \cite[Theorem 4.13]{bayer2020desingularization}, $\beta=0$ is the unique vertical wall of $F$. Any other wall is a semicircle centered along the $\beta$-axis, and its apex lies on the hyperbola $\mu_{\alpha, \beta}(F)=0$. Moreover, no two walls intersect. 

Note that when $\mu_{\alpha, \beta}(F)=0$ holds, we have $\beta<-\sqrt{\frac{2}{5}}<-\frac{1}{2}$, thus we know that there is no semicircular wall centered in the interval $-\frac{1}{2}\leq \beta<0$. Therefore, any semicircular wall in the range $-\frac{1}{2}\leq \beta<0$ will intersect $\beta=-\frac{1}{2}$. To prove the statement, we only need to show that there are no walls when $\beta=-\frac{1}{2}$. This follows from the fact that $\ch^{-\frac{1}{2}}_1(F)$ is minimal.
\end{proof}

\begin{lemma}
\label{proposition_stability_conic}
Let $C\subset X$ be a conic on a GM  threefold $X$ such that $I_C \in \cA_X$. Then $I_C[1]$ is stable with respect to every Serre-invariant stability condition on $\cA_X$.
\end{lemma}

\begin{proof}
By Lemma \ref{bms lemma 2.7} and Lemma \ref{beta=0}, we know that $I_C$ is $\sigma_{\alpha, \beta}$-semistable for every $(\alpha, \beta)\in V$. Since $I_C$ is torsion-free, $I_C[1]\in \Coh^0_{\alpha, \beta}(X)$ is $\sigma^0_{\alpha, \beta}$-semistable. Thus $I_C[1]\in \cA(\alpha, \beta)$ is $\sigma(\alpha, \beta)$-semistable. Then stability with respect to every Serre-invariant stability condition follows from Theorem \ref{Serre-invariant-GM} and Theorem \ref{all_in_one_orbit}.
\end{proof}

\subsection{The Bridgeland moduli space of class $-x$}

In this subsection, we are going to describe the Bridgeland moduli space $\mathcal{M}_{\sigma}(\mathcal{A}_X, -x)$ in Theorem \ref{theorem_irreduciblecomponent_conics}.

The proofs in this section seem technical. However, the only results in this section that will be used in other sections are Proposition \ref{moduli -x bijective} in the proof of Theorem \ref{theorem_irreduciblecomponent_conics}, so there is no harm for readers in skipping this whole section and assuming Proposition \ref{moduli -x bijective} and Theorem \ref{theorem_irreduciblecomponent_conics}.

We start with two lemmas.

\begin{lemma} \label{non existence of sheaf}
Let $X$ be a GM threefold and $E$ a $\mu$-semistable sheaf on $X$ with truncated Chern character $\mathrm{ch}_{\leq 2}(E)=(2,-H,aL)$. If $a\geq 1$ and $c_3(E)\geq 0$ and then we have $E\cong \mathcal{E}$.

\end{lemma}

\begin{proof}
By Lemma \ref{Li BG}, we have $a\leq 1$ which means $a=1$ by our assumption. Then $c_1(E)=-1$ and $c_2(E)=4$. Since  $c_3(E)\geq 0$, by \cite[Proposition 3.5]{brambilla2014vector} we have $\chi(E)=0$. This implies $c_3(E)=0$. Moreover, $E^{\vee \vee}$ also satisfies the assumptions above. Hence by the previous argument, we have $c_1(E^{\vee \vee})=-1, c_2(E^{\vee \vee})=4$ and $c_3(E^{\vee \vee})=0$ as well. In other words, $E=E^{\vee \vee}$. Since $c_3(E)=0$ and $E$ is reflexive of rank two, it is a vector bundle. Moreover, $E$ is a globally generated bundle by \cite[Proposition 3.5]{brambilla2014vector}. Thus $E\cong \cE$ by \cite[Proposition 4.1]{debarre2012period}.
\end{proof}

\begin{lemma} \label{unqiue_Q}
Let $X$ be a GM threefold and $E$ a $\mu$-semistable sheaf on $X$ with $\ch(E)=\ch(\cQ)$. Then we have $E\cong \cQ$.
\end{lemma}

\begin{proof}
First we show that $h^2(E)=0$; then from $\chi(E)=5$ we have $h^0(E)\geq 5$. Indeed, if $h^2(E)\neq 0$, then $\Hom(E, \oh_X(-H)[1])\neq 0$ by Serre duality. Therefore, we have a non-trivial extension
\[0\to \oh_X(-H)\to F\to E\to 0.\]
If $F$ is not $\mu$-semistable, then by the stability of $\oh_X(-H)$ and $E$, the minimal destabilizing quotient sheaf $F'$ of $F$ has $\ch_{\leq 1}(F')=(1,-H)$. Thus $F'^{\vee \vee}\cong \oh_X(-H)$. But if we apply $\Hom(-, \oh_X(-H))$ to the exact sequence above, we obtain $\Hom(F,\oh_X(-H))=0$ since this extension is non-trivial, which gives a contradiction. Then $F$ is $\mu$-semistable with $\ch_{\leq 2}(F)=(4,0,4L)$, which is impossible since $\Delta(F)<0$.

Now we can take five linearly independent elements in $H^0(E)$ and obtain a map $t \colon \oh^{\oplus 5}_X\to E$. From the stability of $\oh_X$ and $E$, we have $\mu(\Im(t))=0$ or $\mu(\Im(t))=\frac{1}{3}$. But the first case cannot happen, since then $\Im(t)$ is the direct sum of a number of copies of $\oh_X$, and this contradicts the construction of $t$. Thus $\mu(\Im(t))=\frac{1}{3}$ and $\ch_{\leq 1}(\Im(t))=(3,H)$. Also $\ch_{\leq 2}(\ker(t))=(2,-H, xL)$, where $x\geq 1$. Note that $\ker(t)$ is reflexive, thus we have $c_3(\ker(t))\geq 0$ since $\ker(t)$ has rank two. Then by stability of $\oh_X$ and $\Hom(\oh_X, \ker(t))=0$, it is not hard to see that $\ker(t)$ is $\mu$-semistable. Thus by Lemma \ref{non existence of sheaf} we have $\ker(t)\cong \cE$. Therefore $\ch(\Im(t))=\ch(E)$ and thus $t$ is surjective.

Now applying $\Hom(\cQ, -)$ to the exact sequence
\[0\to \cE\to  \oh^{\oplus 5}_X\to E\to 0,\]
from $\RHom^\bullet(\cQ, \oh_X)=0$ and $\Ext^1(\cQ, \cE)=k$ we have $\Hom(\cQ, E)=k$. Thus from the stability of $E$ and $\cQ$, we have $E\cong \cQ$ and the result follows.
\end{proof}

Now we introduce some notations. Let $\alpha>0$ and $\beta<0$. For an object $E\in \D^b(X)$, the limit central charge $Z^0_{0,0}(E)$ is defined as the limit of $Z^0_{\alpha, \beta}(E)$ when $(\alpha, \beta) \to (0, 0)$. Note that $Z^0_{\alpha, \beta}(E)$ is given by $\mathbb{Q}$-linear combinations of $\alpha, \beta, \alpha^2, \beta^2$, thus such a limit $Z^0_{0,0}(E)$ always exists. For $Z^0_{0,0}(E)\neq 0$, we can also define the limit slope $\mu^0_{0,0}(E)$ as follows:

\begin{itemize}
    \item If $\Im(Z^0_{0,0}(E))\neq 0$, then we define $\mu^0_{0,0}(E):=-\frac{\Re(Z^0_{0,0}(E))}{\Im(Z^0_{0,0}(E))}$.
    
    \item If $\Im(Z^0_{0,0}(E))=0$ and $\Re(Z^0_{0,0}(E))>0$, then we define $\mu^0_{0, 0}(E):=-\infty$.
    
    \item If $\Im(Z^0_{0,0}(E))=0$ and $\Re(Z^0_{0,0}(E))<0$, then we define $\mu^0_{0, 0}(E):=+\infty$.
\end{itemize}

Note that $Z^0_{0,0}(E)=0$ if and only if $\ch_{\leq 2}(E)$ is a multiple of $\ch_{\leq 2}(\oh_X)$. 

Let $E\in \Coh^0_{\alpha, \beta}(X)$.  By continuity, we can find a neighborhood $U_E$ of the origin such that for any $(\alpha, \beta)\in U_E$, the slopes $\mu^0_{\alpha, \beta}(E)$ and $\mu^0_{0,0}(E)$ are both negative or positive. 
Let $F\in \Coh^0_{\alpha, \beta}(X)$ be another object such that $E,F$ are both $\sigma^0_{\alpha, \beta}$-semistable in a neighborhood $U_{E,F}$ of the origin. If $\mu^0_{0,0}(E)>\mu^0_{0,0}(F)$, then by continuity, we can find a smaller neighborhood $U'_{E,F}$ such that $\mu^0_{\alpha, \beta}(E)>\mu^0_{\alpha, \beta}(F)$ holds for every $(\alpha, \beta)\in U'_{E,F}$. Thus we have $\Hom(E, F)=0$. We will use these two elementary facts repeatedly.


\begin{proposition} \label{F iso IC}
If $F\in \mathcal{A}(\alpha, \beta)$ is $\sigma(\alpha, \beta)$-stable such that $[F]=-x$ and $F$ is $\sigma^0_{\alpha, \beta}$-semistable for some $(\alpha, \beta)\in V$, then $F\cong I_C[1]$ for a conic $C$ on $X$.
\end{proposition}

\begin{proof}
Since $F$ is $\sigma^0_{\alpha, \beta}$-semistable and $\mu^0_{\alpha, \beta}(F)>0$, as in \cite[Proposition 4.6]{pertusi2020some} there is a triangle
\[F_1[1]\to F\to F_2\]
where $F_1\in \Coh^{\beta}(X)$ with $\mu^+_{\alpha, \beta}(F_1)<0$ and $F_2$ is supported on points. Thus $\ch(F_1)=(1,0,-2L, mP)$, where $m$ is the length of $F_2$. By Lemmas \ref{beta=0} and \ref{bms lemma 2.7}, $F_1$ is a rank one torsion-free sheaf, hence it is the ideal sheaf of a closed subscheme. Thus by \cite[Corollary 1.38]{sanna2014rational}, we have $m\leq 0$, which means $F_2=0$ and $F_1\cong F[-1]$. Thus by Lemma \ref{beta=0} again, $F[-1]$ is a $\mu$-semistable torsion free sheaf, which is of the form $F[-1]\cong I_C$ for a conic $C$ on $X$ since $\Pic(X)=\mathbb{Z}\cdot H$.
\end{proof}

When $F$ is not $\sigma^0_{\alpha, \beta}$-semistable for $(\alpha, \beta)\in V$, the argument is more complicated. Our main tools are the inequalities in \cite{petkovic2020note}, \cite[Proposition 4.1]{pertusi2020some}, Lemma \ref{Li BG} and Theorem \ref{Naoki BG}, which allow us to bound the rank and first two Chern characters $\ch_1,\ch_2$ of the destabilizing objects and their cohomology objects. Since $F\in \cA_X$, by using the Euler characteristics $\chi(\mathcal{O}_X,-)$ and $\chi(\cE^{\vee},-)$ we can obtain a bound on $\ch_3$. Finally, via a similar argument as in Lemma \ref{non existence of sheaf} we deduce that the Harder--Narasimhan factors of $F$ are the ones we expect.

\begin{proposition} \label{F iso prIC}
If $F\in \mathcal{A}(\alpha, \beta)$ is $\sigma(\alpha, \beta)$-stable such that $[F]=-x$ and $F$ is not $\sigma^0_{\alpha, \beta}$-semistable for every $(\alpha, \beta)\in V$, then $F$  fits into a triangle
\[\mathcal{E}[2]\to F\to \mathcal{Q}^{\vee}[1].\]
\end{proposition}

\begin{proof}
Since there are no walls for $F$ tangent to the wall $\beta=0$, by the local finiteness of walls and  \cite[Proposition 2.2.2]{bayer2011bridgeland} we can find an open neighborhood $U'$ of the origin such that the Harder--Narasimhan filtration with respect to $\sigma^0_{\alpha, \beta}$ is constant for every $(\alpha, \beta)\in U:=U'\cap V$. In the following we will only consider $\sigma^0_{\alpha, \beta}$ for $(\alpha, \beta)\in U$.

Let $B$ be the minimal destabilizing quotient object of $F$ and $0\to A\to F\to B\to 0$ be the destabilizing short exact sequence of $F$ in $\Coh^0_{\alpha, \beta}(X)$. Hence we know that $A, B\in \mathrm{Coh}^0_{\alpha, \beta}(X)$ and $B$ is $\sigma^0_{\alpha, \beta}$-semistable with $\mu^{0,-}_{\alpha, \beta}(A)>\mu^0_{\alpha, \beta}(F)> \mu^{0}_{\alpha, \beta}(B)$ for all $(\alpha, \beta)\in U$. By \cite[Remark 5.12]{bayer2017stability}, we have $\mu^0_{\alpha,\beta}(B)$ $\geq \min\{\mu^0_{\alpha, \beta}(F), \mu^0_{\alpha, \beta}(\mathcal{O}_X), \mu^0_{\alpha, \beta}(\mathcal{E}^{\vee})\}$. Hence the following relations hold for all $(\alpha, \beta)\in U$:

\begin{enumerate}[itemsep=0pt,parsep=0pt,label=(\alph*)]
    \item $\mu^0_{\alpha, \beta}(A)> \mu^0_{\alpha, \beta}(F) > \mu^0_{\alpha, \beta}(B)$,
    
    \item $\mathrm{Im}(Z^0_{\alpha, \beta}(A))\geq 0$, $\mathrm{Im}(Z^0_{\alpha, \beta}(B))> 0$,
    
    \item $\mu^0_{\alpha, \beta}(B)\geq \min\{\mu^0_{\alpha, \beta}(F), \mu^0_{\alpha, \beta}(\mathcal{O}_X), \mu^0_{\alpha, \beta}(\mathcal{E}^{\vee})\}$,
    
    \item $\Delta(B)\geq 0$.
    
\end{enumerate}

By continuity and taking $(\alpha, \beta)\to (0,0)$, we have:

\begin{enumerate} \label{1234}
    \item $\mu^0_{0,0}(A)\geq \mu^0_{0,0}(F)= 0 \geq  \mu^0_{0,0}(B)$,
    
    \item $\mathrm{Im}(Z^0_{0,0}(A))\geq 0$, $\mathrm{Im}(Z^0_{0,0}(B))\geq 0$,
    
    \item $\mu^0_{0,0}(B)\geq \min\{\mu^0_{0,0}(F), \mu^0_{0,0}(\mathcal{O}_X), \mu^0_{0,0}(\mathcal{E}^{\vee})\}$,
    
    \item $\Delta(B)\geq 0$.
    
\end{enumerate}

Assume that  $[A]=a[\mathcal{O}_X]+b[\mathcal{O}_H]+c[\mathcal{O}_L]+d[\mathcal{O}_P]$  for integers $a,b,c,d\in \ZZ$. Then we have $[B]=(-1-a)[\mathcal{O}_X]-b[\mathcal{O}_H]+(2-c)[\mathcal{O}_L]-(1+d)[\mathcal{O}_P]$. Then we see

\begin{itemize}
    \item $\mathrm{ch}(A)=(a, bH, \frac{c-5b}{10}H^2, \frac{\frac{5}{3}b+\frac{c}{2}+d}{10}H^3)$

    \item $Z^0_{0,0}(A)=bH^3+(\frac{c-5b}{10}H^3)\mathfrak{i}$, $Z^0_{0,0}(B)=-bH^3+(\frac{2-c+5b}{10}H^3)\mathfrak{i}$

    \item $\mu^0_{0,0}(A)=\frac{10b}{5b-c}, \mu^0_{0,0}(B)=\frac{-10b}{c-5b-2}$.
\end{itemize}

Note that $[F]=-[\mathcal{O}_X]+2[\mathcal{O}_L]-[\mathcal{O}_P]$. From $(2)$ we know $c-5b=0$, $1$ or $2$. But when $c-5b=2$, it is not hard to see that $(c)$ fails near the origin. Thus $c-5b=0$ or $1$.

We begin with two claims.

\medskip

\textbf{Claim 1:} \emph{We have $\RHom^\bullet(\oh_X, B)=\Hom(\oh_X, B)$ and $\RHom^\bullet(\oh_X, A)=\Ext^1(\oh_X, A)[-1]$.}

Since $F\in \cA_X$, we only need to prove that $\Ext^i(\oh_X, A)=0$ for $i\neq 1$. Indeed, since $\oh_X\in \Coh^0_{\alpha, \beta}(X)$ and $F\in \cA_X$, we have $\Ext^i(\oh_X, A)=0$ for all $i\leq 0$. Also, by Serre duality we have $\Ext^i(\oh_X, A)=\Hom(A, \oh_X(-H)[3-i])$. Thus from the fact that $\oh_X(-H)\in \Coh^0_{\alpha, \beta}(X)$, we obtain $\Hom(A, \oh_X(-H)[3-i])=0$ for $i\geq 2$. Therefore we have $\Ext^i(\oh_X, A)=0$ for $i\neq 1$.

\medskip

\textbf{Claim 2:} \emph{We have $\RHom^\bullet(\cE^{\vee}, B)=\Hom(\cE^{\vee}, B)$ and $\RHom^\bullet(\cE^{\vee}, A)=\Ext^1(\cE^{\vee}, A)[-1]$.}

Since $\cE^{\vee}$ and $\cE[2]\in \Coh^0_{\alpha, \beta}(X)$, the argument is the same as Claim 1.

\medskip

Now we deal with the cases $c-5b=0$ and $c-5b=1$ separately.

\textbf{Case 1 ($c-5b=0$):}

First, we assume that $c-5b=0$. By \ref{1234}, we have:

\begin{enumerate}
     \item $-2\leq b\leq 0$,

    \item $b^2+\frac{2a+2}{5}\geq 0$.
\end{enumerate}

\textbf{case 1.1 ($b=0$):} If $b=0$, then $c=0$ and $a\geq -1$. In this case we have $\ch_{\leq 2}(B)=(-1-a,0,2L)$. 
If $a=-1$, then $\ch_{\leq 2}(A)=\ch_{\leq 2}(\oh_X[1])=(-1,0,0)$, which is impossible since $\Im(Z^0_{\alpha, \beta}(A))<0$ for $(\alpha, \beta)\in V$. Thus $a\geq 0$, and $a\neq 0$ otherwise $\mu^0_{\alpha, \beta}(F)=\mu^0_{\alpha, \beta}(B)$ for any $(\alpha, \beta)\in V$. But then we have $\mu^0_{\alpha, \beta}(F)< \mu^0_{\alpha, \beta}(B)$ when $(\alpha, \beta)\in U$ is sufficiently close to the origin. This contradicts our assumption on $B$.

\textbf{case 1.2 ($b=-1$):} If $b=-1$, we have $c=-5$. In this case $\mathrm{ch}_{\leq 2}(A)=(a, -H, 0)$. Since $A\in \Coh^0_{\alpha, \beta}(X)$, we have $\mathrm{Im}(Z^0_{\alpha, \beta}(A))\geq 0$ for every $(\alpha, \beta)\in U$. Note that $\mathrm{Im}(Z^0_{\alpha, \beta}(A))=(\beta+\frac{a(\beta^2-\alpha^2)}{2})H^3$ and $0<\alpha<-\beta$, and we have $a\geq \frac{-2\beta}{\beta^2-\alpha^2}$. But note that when $\alpha=\frac{-\beta}{2}$ and $\beta\to 0$, we have  $\frac{-2\beta}{\beta^2-\alpha^2}\to +\infty$, thus we get a contradiction since $a$ is a finite number.

\textbf{case 1.3 ($b=-2$):} If $b=-2$, we have $c=-10$. In this case we have $\ch_{\leq 2}(A)=(a, -2H, 0)$. Similarly to case 1.2, we have $\mathrm{Im}(Z^0_{\alpha, \beta}(A))\geq 0$ for every $(\alpha, \beta)\in U$. Note that $\mathrm{Im}(Z^0_{\alpha, \beta}(A))=(2\beta+\frac{a(\beta^2-\alpha^2)}{2})H^3$ and $\alpha<-\beta$, and we have $a\geq \frac{-4\beta}{\beta^2-\alpha^2}$. Then as in case 1.2, we get a contradiction. 

\medskip

\textbf{Case 2 ($c-5b=1$):}  Now we assume that $c-5b=1$. Then by \ref{1234}, we have:

\begin{enumerate}
    \item $-1\leq b\leq 0$,
    
    
    \item $b^2+\frac{a+1}{5}\geq 0$.
\end{enumerate}

\textbf{case 2.1 ($b=0$):} If $b=0$, then $c=1$. Therefore $-1\leq a$. 
If $a=-1$, since $B$ is $\sigma^0_{\alpha, \beta}$-semistable, we know $\mathcal{H}^0_{\mathrm{Coh}^{\beta}(X)}(B)$ is either $0$ or supported on points. Thus $\mathrm{ch}_{\leq 2}(\mathcal{H}^{-1}_{\mathrm{Coh}^{\beta}(X)}(B))=(0,0, -L)$. But $\mathrm{Re}(Z_{\alpha, \beta}(\mathcal{H}^{-1}_{\mathrm{Coh}^{\beta}(X)}(B)))>0$ which is impossible since $\mathcal{H}^{-1}_{\mathrm{Coh}^{\beta}(X)}(B)\in \mathrm{Coh}^{\beta}(X)$ with $\mathrm{Im}(Z_{\alpha, \beta}(\mathcal{H}^{-1}_{\mathrm{Coh}^{\beta}(X)}(B)))=0$.

Therefore we have $a\geq 0$. Hence $\mathrm{ch}_{\leq 2}(B)=-(a+1,0, -L)$, where $a+1\geq 1$. This is also impossible since when $(\alpha, \beta)\in U$ is sufficiently close to the origin, we have $\mu^0_{\alpha, \beta}(B)>\mu^0_{\alpha, \beta}(F)$.

\textbf{case 2.2 ($b=-1$):}
We have $b=-1$ and $c=-4$. Hence $-6\leq a$. In this case $\ch_{\leq 2}(B)=(-1-a, H, L)$ and we have $\mu^0_{\alpha, \beta}(B)<0$ and $\ch_1^{\beta}(B)>0$ for when $(\alpha, \beta)\in U$ is sufficiently close to the origin. Thus $B\in \Coh^{\beta}(X)$ is $\sigma_{\alpha, \beta}$-semistable.  Applying Lemma \ref{Li BG} to $B$, we have $a\geq -3$. 

We first prove a claim.

\medskip

\textbf{Claim 3:} \emph{In the situation of case 2.2, we have $A$ is $\sigma^0_{\alpha, \beta}$-semistable. Hence $\RHom^\bullet(\oh_X, A)=0$, $\ch(A)=(a, -H, L, (\frac{7}{3}-a)P)$ and $\chi(\cE^{\vee}, A)=3-2a$.}

Assume $A$ is not $\sigma^0_{\alpha, \beta}$-semistable for some $(\alpha, \beta)\in U$. Then we can take a neighborhood $U'_A$ of the origin such that $A$ has constant Harder--Narasimhan factors for any $(\alpha, \beta)\in U_A:=U\cap U'_A\cap V$. Let $C$ be the minimal  destabilizing quotient object of $A$ with respect to $\sigma^0_{\alpha, \beta}$ for $(\alpha, \beta)\in U_A$.
In this case we have $\ch_{\leq 2}(A)=(a,-H,L)$. Since $\mathrm{Im}(Z^0_{0,0}(A))=\frac{1}{10}H^3$, we know that $\mathrm{Im}(Z^0_{0,0}(C))=0$ or $\frac{1}{10}H^3$. If $\mathrm{Im}(Z^0_{0,0}(C))=0$, then $\mu^0_{0,0}(C)=+\infty$ or $-\infty$. But the previous case contradicts $\mu^0_{\alpha, \beta}(A)>\mu^0_{\alpha, \beta}(C)$ and the latter case contradicts $\mu^0_{\alpha, \beta}(C)>\mu^0_{\alpha, \beta}(F)$. Therefore we have $\mathrm{Im}(Z^0_{0,0}(C))=\frac{1}{10}H^3$ and we can assume that $\ch_{\leq 2}(C)=(e,fH,L)$ where $e,f\in \mathbb{Z}$. Since $\mu^0_{0,0}(A)\geq \mu^0_{0,0}(C)\geq \mu^0_{0,0}(F)=0$, we have $10\geq -10f\geq 0$. If $f=0$, then $\ch_{\leq 2}(C)=(e,0,L)$ and $\ch_{\leq 2}(D)=(a-e,-H,0)$, where $D=\mathrm{cone}(A\to C)[-1]$. Then $\mu^{0-}_{\alpha, \beta}(D)>\mu^0_{\alpha, \beta}(A)$ for any $(\alpha, \beta)\in U_A$. Hence $$\mu^0_{\alpha, \beta}(D)=\frac{1+(a-e)\beta}{\beta+\frac{a-e}{2}(\beta^2-\alpha^2)}>\mu^0_{\alpha, \beta}(A)>\mu^0_{\alpha, \beta}(F).$$ 


But note that $\mu^0_{\alpha, \beta}(D)<0$ for $(\alpha, \beta)\in U_A$ that sufficiently closed to the origin, which gives a contradiction since $\mu^0_{\alpha, \beta}(D)>\mu^0_{\alpha_0, \beta_0}(F)$ holds for any $(\alpha, \beta)\in U_A$. 

Therefore the only possible case is $f=-1$, and hence $\mu^0_{0,0}(C)=10$. Since $\mu^0_{\alpha, \beta}(A)>\mu^0_{\alpha, \beta}(C)$ for $(\alpha, \beta)\in U_A$, we have $\rk C>a$. But this is impossible since $D, \oh_X\in \Coh^0_{\alpha, \beta}(X)$ but  $\ch_{\leq 2}(D)=(s, 0, 0)=s\cdot \ch_{\leq 2}(\oh_X)$ where $s=a-\rk C<0$. Now for the last statement, note that $\oh_X(-H)[2]\in \Coh^0_{\alpha, \beta}(X)$ is $\sigma^0_{\alpha, \beta}$-semistable with $\mu^0_{0,0}(\oh_X(-H)[2])=2$, hence we have $\Hom(A, \oh_X(-H)[2])=\Hom(\oh_X, A[1])=0$. Now combined with Claim 1, this proves our claim.

\medskip

Now we deal with the three cases $a=-3$, $-2\leq a\leq 1$ and $a\geq 2$ separately.

When $a=-3$, we have $\ch_{\leq 2}(B)=\ch_{\leq 2}(\cE^{\vee})$. Then since $\ch_{\leq 2}(B)$ is on the boundary of Lemma \ref{Li BG}, by a standard argument we know that $B$ is $\sigma_{\alpha, \beta}$-semistable for every $\alpha>0$ and $\beta<0$, as explained in \cite[Proposition 3.2]{Pertusi2021serreinv}. Thus by Lemma \ref{bms lemma 2.7}, $B$ is a $\mu$-semistable sheaf. From Claim 3 we have $\chi(\oh_X, B)=0$, hence $\ch(B)=\ch(\cE^{\vee})$ and by Lemma \ref{non existence of sheaf} we have $B\cong \cE^{\vee}$. But this implies $\Hom(\oh_X, A[1])=k^5$ since $F\in \cA_X$, which contradicts Claim 3.

When $-2\leq a\leq 1$, we have $\mu^{0}_{\alpha, \beta}(A)> \mu^0_{\alpha, \beta}(\cE[2])$. Since $A$ is $\sigma^0_{\alpha, \beta}$-semistable, we have $\Hom(A, \cE[2])=\Hom(\cE^{\vee}, A[1])=0$. Thus $\RHom^\bullet(\cE^{\vee}, A)=0$ by Claim 2. But this contradicts Claim 3 since $\chi(\cE^{\vee}, A)=3-2a$.

When $a\geq 2$, applying Theorem \ref{Naoki BG} to $B$, we have $a=2$. Thus $\ch_{\leq 2}(B)=\ch_{\leq 2}(\cQ^{\vee}[1])$. By Claim 3, we know that $\RHom^\bullet(\oh_X, B)=0$ and we get $\ch(B)=\ch(\cQ^{\vee}[1])$. Thus $\chi(\cE^{\vee}, B)=\hom(\cE^{\vee}, B)>0$. Therefore, if we apply $\Hom(-, B)$ to the exact sequence $0\to \cQ^{\vee}\to \oh_X^{\oplus 5}\to \cE^{\vee}\to 0$, we obtain $\hom(\cQ^{\vee}[1], B)>0$. Now by stability, we have $B\cong \cQ^{\vee}[1]$. Now $\ch(A)=\ch(\cE[2])$. By Claim 2 and Claim 3, we have $\ext^1(\cE^{\vee}, A)=\hom(A, \cE[2])=1$. Since $A$ is $\sigma^0_{\alpha, \beta}$-semistable and $\cE[2]$ is $\sigma^0_{\alpha, \beta}$-stable, we have $A\cong \cE[2]$.
\end{proof}

\begin{proposition} \label{moduli -x bijective}
Let $X$ be a GM threefold. Then every object in the moduli space $\cM_{\sigma}(\cA_X, -x)$ is of form $\pr(I_C)[1]$ for a conic $C\subset X$.
\end{proposition}

\begin{proof}

Note that $\mathrm{hom}(\cQ^{\vee}[1], \cE[2])=1$. Then the result follows from Proposition \ref{F iso IC} and Proposition \ref{F iso prIC}.
\end{proof}

Now we are ready to realize the Bridgeland moduli space $\mathcal{M}_{\sigma}(\mathcal{A}_{X},-x)$ as the contraction $\cC_m(X)$ of the Fano surface $\cC(X)$:

\begin{theorem}
\label{theorem_irreduciblecomponent_conics}
Let $X$ be a GM threefold and $\sigma$ a Serre-invariant stability condition on $\cA_X$. The projection functor $\pr \colon \D^b(X)\rightarrow\mathcal{A}_{X}$ induces a surjective morphism $p \colon \cC(X)\ra\mathcal{M}_{\sigma}(\mathcal{A}_{X},-x)$, where $p$ is

\begin{itemize}
    \item a blow-down morphism to a smooth point when $X$ is ordinary;
    
    \item a contraction of the component $\mathbb{P}^2$ to a singular point when $X$ is special.
\end{itemize}

In particular, when $X$ is general and ordinary,  $\mathcal{M}_{\sigma}(\mathcal{A}_{X},-x)$ is isomorphic to the minimal model $\mathcal{C}_m(X)$ of the Fano surface of conics on $X$.  When $X$ is general and special, the moduli space $\mathcal{M}_{\sigma}(\mathcal{A}_{X},-x)$ has only one singular point.
\end{theorem}

\begin{proof}
By Lemma \ref{pr(IC)[1] when h0(EC)=1 stability} and Lemma \ref{proposition_stability_conic}, $\pr(I_C)[1]$ is $\sigma$-stable for any conic $C\subset X$. Then we obtain a morphism $p\colon \cC(X)\to \mathcal{M}_{\sigma}(\mathcal{A}_{X},-x)$. Moreover, Proposition \ref{moduli -x bijective} implies that $p$ is surjective.

Now according to Proposition~\ref{prop_contraction_P1_OGM}, the family of conics $C\subset X$ with the property that $I_C\not\in\cA_X$ is parametrized by the line $L_{\sigma}$ when $X$ is ordinary, and the component $\PP^2$ when $X$ is special. Since $\pr(I_C)[1]\cong \pi'(\cQ^{\vee})[1]$ for $I_C\notin \cA_X$ by Proposition \ref{prop_conic_inKuz}, we know that $p(L_{\sigma})=[\pi'(\cQ^{\vee})[1]]$ when $X$ is ordinary, and $p(\PP^2)=[\pi'(\cQ^{\vee})[1]]$ when $X$ is special, where $[\pi'(\cQ^{\vee})[1]]\in \mathcal{M}_{\sigma}(\mathcal{A}_{X},-x)$ is the point represented by the object $\pi'(\cQ^{\vee})[1]$. Thus $p$ is a blow-down morphism to a smooth point when $X$ is ordinary and a contraction of the component $\mathbb{P}^2$ to a singular point when $X$ is special by Lemma \ref{AX_piQ_stable}.

When $X$ is general and ordinary, the Fano surface $\cC(X)$ is smooth by Theorem \ref{ordinary-CX}. Thus $\mathcal{M}_{\sigma}(\mathcal{A}_{X},-x)$ is a smooth surface obtained by blowing down a smooth rational curve $L_{\sigma}$ on the smooth irreducible projective surface $\cC(X)$. This implies that $\mathcal{M}_{\sigma}(\mathcal{A}_{X},-x)$ is also a smooth irreducible projective surface. On the other hand, it is known that there is a unique rational curve $L_{\sigma}\subset\mathcal{C}(X)$ and it is the unique exceptional curve by Lemma \ref{minimal surface of S(X) lemma}. Thus $\mathcal{M}_{\sigma}(\mathcal{A}_{X},-x)$ is isomorphic to the minimal model $\mathcal{C}_m(X)$ of Fano surface of conics on $X$.  

When $X$ is general and special, the last statement follows from Theorem \ref{special_surface} and Lemma \ref{lemma_ideal_projection_smooth}.
\end{proof}

\subsection{Involutions on $\cC_m(X)$}

In this section, we are going to describe the involution $\iota$ on $\cC_m(X)$ in  Theorem \ref{classical_invo_conic}, described in \cite[Section 5.2]{debarre2012period} using the involution on $\cA_X$. Recall that there is a natural involutive autoequivalence functor of $\cA_X$, denoted by $\tau_{\cA}$ (cf.~Remark \ref{rmk-tau}). When $X$ is special, it is induced by the natural involution $\tau$ on $X$, which comes from the double cover $X\to Y$. In this case it is easy to see that $\tau_{\cA}(\pr(I_C))\cong \pr(I_{\tau (C)})$.

When $X$ is ordinary, the situation is more subtle. In the following, we describe the action of $\tau_{\cA}$ on the projection into $\cA_X$ of an ideal sheaf of a conic $\mathrm{pr}(I_C)$ in this case, and compare with the involution $\iota$ on $\cC_m(X)$ described in \cite[Section 5.2]{debarre2012period}.

\begin{proposition}\label{involution-acting-conic} Let $X$ be an ordinary GM threefold and $C$ a conic on $X$.
\leavevmode\begin{enumerate}
    \item If $I_C\in\mathcal{A}_X$, then $\mathcal{\tau}_{\mathcal{A}}(I_C)$ is either \begin{enumerate}[label=(\alph*).]
        \item $I_{C'}$ such that $C\cup C'=Z(s)$ for $s\in H^0(\mathcal{E}^{\vee})$, where $Z(s)$ is the zero locus of the section $s$;
        \item or $\pi'(\cQ^{\vee})$, and in this case $C$ is the $\rho$-conic
    \end{enumerate} 
    \item If $I_C\notin\cA_X$, then $\tau_{\cA}(\mathrm{pr}(I_C))\cong I_{C''}$ for the $\rho$-conic $C''\subset X$. 
\end{enumerate}
Therefore, the involution induced by $\tau_{\cA}$ on $\cC_m(X)$ is the same as $\iota$ in Theorem \ref{classical_invo_conic}.
\end{proposition}

\begin{remark}
We can define a birational involution on $\cC(X)$ for any GM threefold $X$ as in Proposition \ref{involution-acting-conic} (1)(a), which is regular on the locus of conics $C$ with $\hom(\cE, I_C)=1$. 
\end{remark}

We first state some lemmas which we require for the proof of the proposition above.

\begin{lemma}\label{injectivity of map} 
Let $X$ be an ordinary GM threefold and $C$ be the $\rho$-conic on $X$. 
Then the natural morphism $s'\colon \cE^{\oplus 2}\to I_C$ is surjective and there is a short exact sequence $$0\rightarrow\cQ(-H)\rightarrow\cE^{\oplus 2}\xrightarrow{s'} I_C\rightarrow 0.$$
\end{lemma}

\begin{proof}
By Lemma \ref{IC-class}, we have $\Hom(\cE, I_C)=k^2$. Thus, taking two linearly independent elements in $\Hom(\cE, I_C)$,  we have a natural map $s'\colon\cE^{\oplus 2}\to I_C$. Moreover, since $\langle C\rangle =\mathrm{Gr}(2,3)$ and $\langle C\rangle\cap X=C$, we know that $s'$ is surjective. Let $K:=\ker(s')$. Then it is not hard to see that $\ch(K)=\ch(\cQ(-H))$. Note that $\Hom(\cE, K)=0$ and $K$ is reflexive. 

We claim that $K$ is $\mu$-semistable. Indeed, suppose $K$ is not $\mu$-semistable and let $K'$ be its maximal destabilizing subsheaf. Then $K'$ is also reflexive. Since $\Hom(\cE, K)=0$, we have $K'\neq \cE$. By the stability of $\cE$ and the fact that $K\subset \cE^{\oplus 2}$, we know that $\mu(K')=-\frac{1}{2}$. Since $\Hom(K', \cE)\neq 0$, by the stability of $K'$ and $\cE$ we have $K'\subset \cE$. Thus from $\ch_{\leq 1}(K')=\ch_{\leq 1}(\cE)$ we know that $\cE/K'$ is supported in codimension $\geq 2$, which gives a contradiction since $\cE$ and $K'$ are both reflexive.

Now the result follows from Lemma \ref{unqiue_Q}, since $K(H)$ is $\mu$-semistable with $\ch(K(H))=\ch(\cQ)$.
\end{proof}

\begin{lemma}\label{contraction}
Let $X$ be an ordinary GM threefold. Let $C$ be a conic on $X$. Then $$ \bL_{\cE}(I_C)=\begin{cases}
\mathbb{D}(I_{C'})\otimes\mathcal{O}_X(-H)[1], & \mathrm{RHom}^\bullet(\mathcal{E},I_C)=k\\
\pi(\mathcal{E}), & \mathrm{RHom}^\bullet(\mathcal{E},I_C)=k^2\oplus k[-1]
\end{cases}$$ such that $C\cup C'=Z(s)$ for $s\in H^0(\mathcal{E}^{\vee})$, where $Z(s)$ is the zero locus of the section $s$
\end{lemma}

\begin{proof}
By Lemma \ref{IC-class}, we have that $\mathrm{RHom}^\bullet(\mathcal{E},I_C)$ is either $k$ or $k^2\oplus k[-1]$.  If $\mathrm{RHom}^\bullet(\mathcal{E},I_C)=k$, then we have the triangle $$\mathcal{E}\rightarrow I_C\rightarrow\bL_{\cE}(I_C).$$
Taking cohomology with respect to the standard heart we get $$0\rightarrow\mathcal{H}^{-1}(\bL_{\cE}(I_C))\rightarrow\cE\xrightarrow{s} I_C\rightarrow\mathcal{H}^0(\bL_{\cE}(I_C))\rightarrow 0.$$
The image of the map $s$ is the ideal sheaf of an elliptic quartic $D=Z(s)$ for $s\in H^0(\cE^{\vee})$, thus we have following two short exact sequences:
$0\rightarrow\mathcal{H}^{-1}(\bL_{\cE}(I_C))\rightarrow\cE\rightarrow I_D\rightarrow 0$ and $0\rightarrow I_D\rightarrow I_C\rightarrow\mathcal{H}^0(\bL_{\cE}(I_C))\rightarrow 0$. Then $\mathcal{H}^{-1}(\bL_{\cE}(I_C))$ is a torsion-free sheaf of rank $1$ with the same Chern character as $\mathcal{O}_X(-H)$. It is easy to show that it must be $\mathcal{O}_X(-H)$. On the other hand $\mathcal{H}^0(\bL_{\cE}(I_C))$ is supported on the residual curve $C'$ of $C$ in $D$ and  $\mathcal{H}^0(\bL_{\cE}(I_C))\cong\mathcal{O}_{C'}(-H)$. Thus we have the triangle $$\mathcal{O}_X(-H)[1]\rightarrow\bL_{\cE}(I_C)\rightarrow\mathcal{O}_{C'}(-H)$$ and we observe that $\bL_{\cE}(I_C)$ is exactly the twisted derived dual of the ideal sheaf $I_{C'}$ of a conic $C'\subset X$, i.e. $\bL_{\cE}(I_C)\cong\mathbb{D}(I_{C'})\otimes\mathcal{O}_X(-H)[1]$. 

If $\mathrm{RHom}^\bullet(\mathcal{E},I_C)=k^2\oplus k[-1]$, then we have the triangle
$$\cE^2\oplus\cE[-1]\rightarrow I_C\rightarrow\bL_{\cE}(I_C) . $$
Taking the long exact sequence in cohomology with respect to the standard heart, we get
$$0\rightarrow\mathcal{H}^{-1}(\bL_{\cE}(I_C))\rightarrow\mathcal{E}^2\xrightarrow{s'} I_C\rightarrow\mathcal{H}^0(\bL_{\cE}(I_C))\rightarrow\cE\rightarrow 0.$$

Now by Lemma~\ref{injectivity of map}, $s'$ is surjective and the cohomology objects are given by $\mathcal{H}^{-1}(\bL_{\cE}(I_C))\cong\mathcal{Q}(-H)$ and $\mathcal{H}^{0}(\bL_{\cE}(I_C))\cong \cE$, which implies that $\bL_{\cE}(I_C)\cong\pi(\cE)$. 
\end{proof}

\begin{proof}[Proof of Proposition~\ref{involution-acting-conic}]
Since $\tau_\mathcal{A}\circ\tau_{\mathcal{A}}\cong \identity$, we have $\tau_{\mathcal{A}}\cong\tau^{-1}_{\mathcal{A}}$. By Proposition \ref{prop-serre-functor}, we have $\tau_{\mathcal{A}}\cong\tau_{\mathcal{A}}^{-1}\cong\bL_{\mathcal{O}_X}\circ\bL_{\mathcal{E}^{\vee}}(-\otimes\mathcal{O}_X(H))[-1]$. Then 
\begin{align*}
    \tau_{\mathcal{A}}(I_C) &\cong\bL_{\mathcal{O}_X}\circ\bL_{\mathcal{E}^{\vee}}(I_C\otimes\mathcal{O}_X(H))[-1] \\
    &\cong \bL_{\mathcal{O}_X}(\bL_{\mathcal{E}}(I_C)\otimes\mathcal{O}_X(H))[-1] .
\end{align*}
The left mutation $\bL_{\mathcal{E}}(I_C)$ is given by $$\mathrm{RHom}^\bullet(\mathcal{E},I_C)\otimes\mathcal{E}\rightarrow I_C\rightarrow\bL_{\mathcal{E}}(I_C).$$  Note that by Lemma \ref{IC-class},  $\mathrm{RHom}^\bullet(\mathcal{E},I_C)$ is either $k$ or $k^2\oplus k[-1]$, and in the latter case $C$ is the unique $\rho$-conic. Then by Lemma~\ref{contraction}, $$ \bL_{\mathcal{E}}(I_C)=\begin{cases}
\mathbb{D}(I_{C'})\otimes\mathcal{O}_X(-H)[1], & \mathrm{RHom}^\bullet(\mathcal{E},I_C)=k\\
\pi(\mathcal{E}), & \mathrm{RHom}^\bullet(\mathcal{E},I_C)=k^2\oplus k[-1]
\end{cases}$$

If $\mathrm{RHom}^\bullet(\mathcal{E},I_C)=k$, then $\tau_{\mathcal{A}}(I_C)\cong\bL_{\mathcal{O}_X}(\mathbb{D}(I_{C'}))$. We have the triangle $$\mathrm{RHom}^\bullet(\mathcal{O}_X,\mathbb{D}(I_{C'}))\otimes\mathcal{O}_X\rightarrow\mathbb{D}(I_{C'})\rightarrow\bL_{\mathcal{O}_X}(\mathbb{D}(I_{C'})) . $$
Note that $\mathrm{RHom}^\bullet(\mathcal{O}_X,\mathbb{D}(I_{C'}))\cong\mathrm{RHom}^\bullet(I_{C'},\mathcal{O}_X)=k\oplus k[-1]$. Then we have the triangle 
\begin{equation} \label{left mutation derived dual triangle}
\mathcal{O}_X\oplus\mathcal{O}_X[-1]\rightarrow\mathbb{D}(I_{C'})\rightarrow\bL_{\mathcal{O}_X}(\mathbb{D}(I_{C'})).
\end{equation}
The derived dual $\mathbb{D}(I_{C'})$ is given by the triangle $\mathcal{O}_X\rightarrow\mathbb{D}(I_{C'})\rightarrow\mathcal{O}_{C'}[-1]$. Then taking cohomology with respect to the standard heart of triangle (\ref{left mutation derived dual triangle}) we have the long exact sequence
$$0=\mathcal{H}^{-1}(\mathbb{D}(I_{C'}))\rightarrow\mathcal{H}^{-1}(\bL_{\mathcal{O}_X}(\mathbb{D}(I_{C'})))\rightarrow\mathcal{O}_X\rightarrow\mathcal{O}_X$$
$$\rightarrow\mathcal{H}^0(\bL_{\mathcal{O}_X}(\mathbb{D}(I_{C'})))\rightarrow\mathcal{O}_X\rightarrow\mathcal{O}_{C'}\rightarrow\mathcal{H}^1(\bL_{\mathcal{O}_X}(\mathbb{D}(I_{C'})))\rightarrow 0 . $$
Thus we have $\mathcal{H}^{-1}(\bL_{\mathcal{O}_X}(\mathbb{D}(I_{C'})))=0$, $\mathcal{H}^{1}(\bL_{\mathcal{O}_X}(\mathbb{D}(I_{C'})))=0$ and $\mathcal{H}^{0}(\bL_{\mathcal{O}_X}(\mathbb{D}(I_{C'})))\cong I_{C'}$. Hence $\tau_{\mathcal{A}}(I_C)\cong\bL_{\mathcal{O}_X}(\mathbb{D}(I_{C'}))\cong I_{C'}$. 

If $\mathrm{RHom}^\bullet(\mathcal{E},I_C)=k^2\oplus k[-1]$, then $\tau_{\mathcal{A}}(I_C)\cong \bL_{\mathcal{O}_X}\circ\bL_{\mathcal{E}^{\vee}}(I_C\otimes\mathcal{O}_X(H))[-1]  \cong\bL_{\mathcal{O}_X}(\pi(\mathcal{E})\otimes\mathcal{O}_X(H)[-1])$ by Lemma \ref{contraction}. Then using the triangle (\ref{piE-seq}), we have $\tau_{\cA}(I_C)\cong \pi'(\cQ^{\vee})$. Then $(2)$ follows from  $\tau_{\mathcal{A}}\cong \tau^{-1}_{\cA}$.  

Now since $\tau_{\cA}=S_{\cA_X}[-2]$ and $\tau_{\cA}$ acts trivially on $\cN(\cA_X)$, it induces an involution on the Bridgeland moduli space of any class with respect to any Serre-invariant stability condition. In particular, $\tau_{\cA}$ induces an involution on $\cC_m(X)\cong \cM_{\sigma}(\cA_X, -x)$ by Theorem \ref{theorem_irreduciblecomponent_conics}. By (1) and (2), this induced involution coincides with $\iota$ in Theorem \ref{classical_invo_conic}, described in \cite[Section 5.2]{debarre2012period}.
\end{proof}

\begin{remark}\label{rmk-open-locus}
Smooth $\tau$-conics form an open subscheme $U$ of $\cC(X)$. Therefore, the open subscheme $U\cap \iota(U)$ parameterizes smooth $\tau$-conics $C$ such that their involutive conics in Proposition \ref{involution-acting-conic} are smooth as well. The same also works for special GM threefolds, but replace $\tau$-conics with conics with $\hom(\cE, I_C)=1$ and $I_C\in \cA_X$, which are parametrized by $\cC(X)\setminus \PP^2$. In other words, for any GM threefold $X$, there is a two-dimensional open subscheme $\cC_1\subset \cC(X)$ parameterizing smooth conics $C$ with $\hom(\cE, I_C)=1$ such that their involutive conics are smooth.
\end{remark}

\section{The moduli space \texorpdfstring{$M_G(2,1,5)$ for GM threefolds}{corresponding to the other (-1)-class}} \label{M_G(2,1,5) section}

In this section, we investigate the moduli space of rank $2$ Gieseker-semistable sheaves on a GM threefold $X$ with Chern classes $c_1 =H, c_2=5L$ and $c_3=0$, denoted $M_G^X(2,1,5)$. We drop $X$ from the notation when it is clear from context on which threefold we work. Note that if $F \in M_G(2,1,5)$, then \[ \ch(F) = (2, H, 0, -\frac{5}{6}P) . \]

We are interested in $M_G(2,1,5)$ since it naturally appears in the description of the period fiber in \cite{debarre2012period}. Our main theorem in Section is Theorem \ref{irreducible component pr' theorem}, which realizes $M_G(2,1,5)$ as the Bridgeland moduli space $\cM_{\sigma}(\cA_X,y-2x)$.

First, we prove a classification result of sheaves in $M_G(2,1,5)$.

\begin{proposition} \label{objects in M_G(2,1,5) theorem}
Let $X$ be a GM threefold and  $F\in M_G(2,1,5)$. Then we have $\RHom^\bullet(\oh_X, F)=k^4$ and $\RHom^\bullet(\oh_X, F(-H))=0$. Moreover, $F$ is either a
\begin{enumerate}
    \item globally generated bundle which  fits into a short exact sequence
    \[ 0 \ra \oh_X \ra F \ra I_Z(H) \ra 0 \]
    where $Z$ is a projective normal smooth elliptic quintic curve; 
    \item non-locally free sheaf with a short exact sequence $$0\rightarrow F\rightarrow\cE^{\vee}\rightarrow\oh_L\rightarrow 0$$ where $L$ is a line on $X$.  Moreover, $F$ is uniquely determined by $L$.
\end{enumerate}
\end{proposition}

\begin{remark}
In \cite[Section 8]{debarre2012period}, they also did computations for non-globally generated bundles in $M^X_G(2,1,5)$. However, in the following proof, we will show such sheaves do not exist.
\end{remark}

\begin{proof}
The first statement follows from \cite[Proposition 3.5 (1)]{brambilla2014vector} and the fact $\chi(F)=4$. (1) and (2) also follow from \cite[Proposition 3.5]{brambilla2014vector} or the argument in \cite[Section 8]{debarre2012period}. So we only need to prove the non-existence of non-globally generated bundles in $M^X_G(2,1,5)$. If $F\in M^X_G(2,1,5)$ is a non-globally generated bundle, then as showed in \cite[Section 8]{debarre2012period}, we have an exact sequence
\[0\to F^{\vee}\to \oh_X^{\oplus 4}\to \cE^{\vee}\xra{a} \oh_L\to 0.\]
By (2), we know that $E:=\ker(a)$ is a non-locally free stable sheaf and $E\in M^X_G(2,1,5)$. Thus we have an exact sequence $0\to F^{\vee}\to \oh_X^{\oplus 4}\to E\to 0$. In particular, $E$ is generated by global sections.

However, we also have the following commutative diagram of exact sequences:
\[\begin{tikzcd}
	0 & {\oh_X^{\oplus4}} & {\oh_X^{\oplus 5}} & {\oh_X} & 0 \\
	0 & E & {\cE^{\vee}} & {\oh_L} & 0
	\arrow[from=2-1, to=2-2]
	\arrow[from=2-2, to=2-3]
	\arrow[from=2-3, to=2-4]
	\arrow[from=2-4, to=2-5]
	\arrow[from=1-1, to=1-2]
	\arrow[from=1-2, to=1-3]
	\arrow[from=1-3, to=1-4]
	\arrow[from=1-4, to=1-5]
	\arrow[from=1-4, to=2-4]
	\arrow[from=1-3, to=2-3]
	\arrow["{\mathrm{ev}}"', from=1-2, to=2-2]
\end{tikzcd}\]
where $\ev\colon \oh_X^{\oplus 4}\to E$ is the evaluation map. Then using Snake Lemma, we have an exact sequence
\[0\to \ker(\ev)\to \cQ^{\vee}\xra{s} I_L\to \mathrm{cok}(\ev)\to 0.\]
As shown in Lemma \ref{IC-class} and Proposition \ref{prop_contraction_P1_OGM}, the image of $s$ is the zero locus of a non-zero section of $\cQ$. It is a $\sigma$-conic when $X$ is ordinary, and a preimage of a line on $Y$ when $X$ is special. Hence in both cases, $\mathrm{im}(s)$ is an ideal sheaf of a conic, and $s$ is not surjective. Therefore, $\ev$ is not surjective as well and we get a contradiction.
\end{proof}

A natural question to ask is what Bridgeland moduli space we get after projecting a sheaf in $M_G(2,1,5)$ into the Kuznetsov component. Since it is easier in this setting, we will work with the alternative Kuznetsov component $\cA_X$ in this section. Our analysis of the projections of objects in $M_G(2,1,5)$ is based on the three cases listed in Proposition \ref{objects in M_G(2,1,5) theorem}. We begin with a Hom-vanishing result.

\begin{lemma} \label{RHom(E^vee, F) vanishing}
Let $X$ be a GM threefold and $F \in M^X_G(2,1,5)$. Then we have  $\RHom^\bullet(\cE^\vee, F) = 0$.
\end{lemma}

\begin{proof} 
By Serre duality and the stability of $\cE^{\vee}$ and $F$, we have $\Hom(\cE^\vee, F)=\Ext^3(\cE^\vee, F)=0$. Since $\chi(\cE^\vee, F)=0$, we only need to show that $\Ext^2(\cE^\vee, F)=0$ or $\Ext^2(\cE^\vee, F)=0$. By Serre duality, we have $\Ext^2(\cE^\vee, F)=\Hom(F, \cE[1])$. Since $\ch^0_{1}(F)=\ch^0_{1}(\cE[1])=1$, by Lemma \ref{bms lemma 2.7} we know that $F$ and $\cE[1]$ are both $\sigma_{\alpha, 0}$-stable for any $\alpha>0$. Then $\Hom(F, \cE[1])=0$ since $\mu_{\alpha, 0}(F)>\mu_{\alpha, 0}(\cE[1])$ when $0<\alpha$ is sufficiently small.
\end{proof}

We are now ready to give an explicit description of $\pr(F)$, for all objects $F \in M_G(2,1,5)$. Recall that for any line $L\subset X$, we have $\cQ|_L\cong \oh_L^{\oplus 2}\oplus \oh_L(1)$. Hence $L$ is contained in a unique $\sigma$-conic $C$. We define \emph{the residue line of $L$} to be the support of $\oh_C\twoheadrightarrow\oh_L$. Note that when $C$ is a double line, we have $L'=L$.

\begin{lemma} \label{pr'(F)=(iota F)^vee[1]}
Let $X$ be a GM threefold and $F \in M_G(2,1,5)$.

\begin{itemize}
    \item If $F$ is a globally generated bundle, then 
    \[\pr(F)\cong \ker(\ev)[1],\]
    where $\ev\colon \oh_X^{\oplus 4}\twoheadrightarrow F$ is the evaluation map.

    \item If $F$ is a non-locally free sheaf determined by a line $L\subset X$, then $\pr(F)$ is the unique object fits into a non-trivial exact triangle
    \[\cE[1]\to \pr(F)\to \oh_{L'}(-1),\]
    where $L'$ is the residue line of $L$.
\end{itemize}
\end{lemma}

\begin{proof}
As a result of Lemma \ref{RHom(E^vee, F) vanishing}, $\bL_{\cE^{\vee}}F = F$, so $\pr(F) = \bL_{\oh_X} F$. By Proposition \ref{objects in M_G(2,1,5) theorem} we have $\RHom^\bullet(\oh_X, F) = k^4$, and the triangle defining the left mutation is 
\begin{equation} \label{pr(F) triangle}
\oh^{\oplus 4}_X \xrightarrow{\ev} F \ra \pr(F). 
\end{equation}
In the cases where $F$ is globally generated, the evaluation map $\ev$ is surjective, so $\pr(F) = \ker(\ev)[1]$.

When $F$ is non-locally free, as in Proposition \ref{objects in M_G(2,1,5) theorem}, we have an exact sequence
\[0\to \ker(\ev)\to \cQ^{\vee}\xra{s} I_L\to \mathrm{cok}(\ev)\to 0.\]
As shown in Lemma \ref{IC-class} and Proposition \ref{prop_contraction_P1_OGM}, the image of $s$ is the zero locus of a non-zero section of $\cQ$, which is a $\sigma$-conic. Hence by Proposition  \ref{prop_contraction_P1_OGM}, we obtain $\ker(\ev)=\cE$ and $\mathrm{cok}(\ev)=\oh_{L'}(-1)$. Since $\pr(F)\in \cA_X$, by Serre duality we have $\RHomb(\cE^{\vee}, \pr(F))=\RHomb(\pr(F), \cE)^{\vee}[-3]=0$, which implies such triangle is non-trivial. And the uniqueness follows from $\Ext^2(\oh_{L'}(-1), \cE)=H^1(\cE(-1)|_L)=k$.
\end{proof}

\subsection{Stability of projection objects}

In the following, we prove the stability of $\pr(F)$ for any $F\in M_G^X(2,1,5)$.

\begin{lemma}
\label{prop_pr_embedding}
The functor $\mathrm{pr} \colon \D^b(X)\rightarrow\cA_X$ induces isomorphisms of $\mathrm{Ext}^k(\mathrm{pr}(F_1),\mathrm{pr}(F_2))$ and $\mathrm{Ext}^k(F_1,F_2)$ for all $k$ and for all $F_1, F_2\in M_G(2,1,5)$. 
\end{lemma}

\begin{proof}
We apply $\mathrm{Hom}(F_1,-)$ to the exact triangle $\mathcal{O}_X^{\oplus 4}\rightarrow F_2\rightarrow\mathrm{pr}(F_2)$. By adjunction of $\pr$ and the inclusion  $\cA_X\hookrightarrow \D^b(X)$, we have $\Ext^k(F_1, \pr(F_2))=\Ext^k(\pr(F_1), \pr(F_2))$ for all $k$. Thus we get a long exact sequence
 \begin{align*}
 \cdots &\rightarrow\mathrm{Ext}^k(F_1, \oh_X)^{\oplus 4}\rightarrow\mathrm{Ext}^k(F_1,F_2)\rightarrow\mathrm{Ext}^k(\mathrm{pr}(F_1),\mathrm{pr}(F_2))\ra \mathrm{Ext}^{k+1}(F_1, \oh_X)^{\oplus 4}\rightarrow \cdots .
 \end{align*} 
 Note that $\mathrm{Ext}^k(F_1, \oh_X)=\Ext^{3-k}(\oh_X, F_1(-H))=0$ for all $k$ by Proposition \ref{objects in M_G(2,1,5) theorem}. Thus the desired result follows. 
\end{proof}

Before we show the stability of projection objects, let us recall a classical result:

\begin{proposition}\label{prop-line-rhom}
Let $X$ be an ordinary GM threefold and $L\subset X$ be a line. Then $\RHomb(\oh_L, \oh_L)=k\oplus k[-1]$ or $k\oplus k^2[-1]\oplus k[-2]$. Moreover, when $X$ is general, we always have $\RHomb(\oh_L, \oh_L)=k\oplus k[-1]$.
\end{proposition}

\begin{proof}
The first statement follows from \cite[Lemma 4.2.1]{iskovskikh1999fano} and the second one follows from \cite[Theorem 4.2.7]{iskovskikh1999fano}.
\end{proof}

Now we are ready to prove the stability of $\pr(F)$.

\begin{proposition}
\label{prop_selfext_group_F}
Let $X$ be a GM threefold and $F\in M^X_G(2,1,5)$. Then we have $\RHom^\bullet(F,F) = k \oplus k^2[-1]$ or $\RHom^\bullet(F,F) = k \oplus k^3[-1]\oplus k[-2]$. Hence $\pr(F)$ is stable with respect to every Serre-invariant stability condition on $\cA_X$.

Moreover, 

\begin{enumerate}
    \item when $X$ is ordinary, if $\RHom^\bullet(F, F) = k \oplus k^3[-1]\oplus k[-2]$ then $F$ is a non-globally generated bundle or a non-locally free sheaf determined by a line $L$, and $[L]\in \Gamma(X)$ is a singular point. In particular, we always have $\RHom^\bullet(F,F) = k \oplus k^2[-1]$ when $X$is general;
    
    \item when $X$ is special, $\RHom^\bullet(F, F) = k \oplus k^3[-1]\oplus k[-2]$ if and only if $\tau^*F\cong F$, where $\tau$ is the natural involution on $X$.
\end{enumerate}

\end{proposition}

\begin{proof}
First, we assume that $X$ is ordinary. We have $\hom(F,F)=1$ and $\ext^3(F,F)=0$ by Serre duality and the stability of $F$. Since $\chi(F,F)=-1$, we need to  prove $\ext^2(F,F)=0$ or $1$.

When $F$ is a globally generated bundle, by the proof of \cite[Theorem 8.2]{debarre2012period}, we have  $\mathrm{ext}^1(F,F)=2$ and $\mathrm{ext}^2(F,F)=0$. When $F$ is non-locally free, there is a mistake made in the proof of \cite[Theorem 8.2]{debarre2012period} and we fix it here. From Proposition \ref{objects in M_G(2,1,5) theorem}, we have an exact sequence $0\to F\to \cE^{\vee}\to \oh_L\to 0$. Since $\RHomb(\cE^{\vee}, F)=0$ by Lemma \ref{RHom(E^vee, F) vanishing}, applying $\Hom(-, F)$ to this exact sequence, we get $\Ext^k(F,F)=\Ext^{k+1}(\oh_L, F)$ for any $k$. Now applying $\Hom(\oh_L, -)$ to this exact sequence, we get a long exact sequence
\[\cdots\to\Ext^2(\oh_L, \oh_L)\to \Ext^3(\oh_L, F)\to \Ext^3(\oh_L, \cE^{\vee})\to 0.\]
By Serre duality, we have $\Ext^3(\oh_L, \cE^{\vee})=H^0(\cE(-1)|_L)=0$. Then from Proposition \ref{prop-line-rhom}, we have $\ext^2(F,F)=\ext^3(\oh_L, F)\leq \ext^2(\oh_L, \oh_L)\leq 1$.  Moreover, if $\ext^2(F,F)=1$, then $\ext^2(\oh_L, \oh_L)=1$. In other words, $[L]\in \Sigma(X)$ is a singular point. This proves (1).

Now we assume that $X$ is special. Then by Lemma \ref{prop_pr_embedding} and Serre duality in $\Ku(X)$, we have 
\begin{align*}
    \mathrm{Ext}^2(F,F) &\cong \mathrm{Ext}^2(\mathrm{pr}(F),\mathrm{pr}(F)) \\
    &\cong \mathrm{Hom}(\mathrm{pr}(F),\tau_{\cA}(\mathrm{pr}(F))) \\
    &\cong \mathrm{Hom}(\mathrm{pr}(F),\mathrm{pr}(\tau^* F)) \cong \Hom(F, \tau^*F),
\end{align*}
where $\tau$ is the involution on $X$ induced by the double cover. Thus when $F\cong \tau^*F$, we have  $\mathrm{Ext}^2(F,F)=k$, and $\mathrm{Ext}^2(F,F)=0$ otherwise. Since $\Ext^3(F,F)=0$ and $\Hom(F,F)=k$, the result follows from $\chi(F,F)=-1$.

Finally, the stability of $\pr(F)$ follows from Lemma \ref{prop_pr_embedding} and Proposition \ref{ext23-stable}.
\end{proof}

\subsection{Involutions on \texorpdfstring{$M_G(2,1,5)$}{the Gieseker moduli space}} \label{involutions of MG(2,1,5) subsection}

In this subsection, we briefly recall the involutions that exist on $M_G(2,1,5)$ and compare it with the one induced by $\tau_{\cA}$. Let $F$ be a globally generated vector bundle, and consider the short exact sequence
\[ 0 \ra \ker(\ev) \ra H^0(X, F) \otimes \oh_X \xrightarrow{\ev} F \ra 0 . \]
Note that $\ker(\ev)$ is a rank $2$ vector bundle with $c_1=-H$ and $c_2=5L$ and no global sections, hence $\ker(\ev)^\vee \in M_G(2,1,5)$. Define $\iota F := \ker(\ev)^\vee$. This bundle $\iota F$ is globally generated, and we have $H^0(X, \iota F) \cong H^0(X, F)^\vee$ \cite[p. 29]{debarre2012period}. This defines a birational involution on $M_G^X(2,1,5)$. 

Note that there is no non-globally generated bundle in $M^X_G(2,1,5)$ by Proposition \ref{objects in M_G(2,1,5) theorem}, then the definition of $\iota$ on the non-locally free locus in \cite[Theorem 8.2]{debarre2012period} does not work. However, we can fix this issue as follows: for any non-locally free $F\in M^X_G(2,1,5)$ determined by a line $L$, we define $\iota(F):=F'$, where $F':=\ker(\cE^{\vee}\twoheadrightarrow \oh_{L'})$ is a non-locally free stable sheaf determined by the residue line $L'$ of $L$. This extends $\iota$ to be a regular involution on $M^X_G(2,1,5)$.

Note that for a special GM threefold, there is another involution on $M_G(2,1,5)$ induced by the involution $\tau$ on $X$,
\[\tau^*\colon M_G(2,1,5)\to M_G(2,1,5), \, \,  F\mapsto \tau^*F.\]
And it is clear that $\tau_{\cA}(\pr(F))\cong \pr(\tau^*F)$.

Now let $X$ be an ordinary GM threefold,  $\tau_{\mathcal{A}}$ be the involution of $\mathcal{A}_X$, and $\iota$ be the geometric involution of $M_G(2,1,5)$ defined above. Then $\tau_{\mathcal{A}}$ induces involutions of the Bridgeland moduli spaces of $\sigma$-stable objects $\mathcal{M}_{\sigma}(\mathcal{A}_X,-x)$ and $\mathcal{M}_{\sigma}(\mathcal{A}_X, y-2x)$. In Proposition~\ref{involution-acting-conic}, we already showed that the action of $\tau_{\cA}$ on $\mathcal{M}_{\sigma}(\mathcal{A}_X,-x)$ induces a geometric involution on $\mathcal{C}_m(X)$. 
In this section, we show that the involution induced by $\tau_{\mathcal{A}}$ is also compatible with $\iota$ on $M_G(2,1,5)$. 

\begin{proposition} \label{MG(2,1,5) involution compatibility}
Let $X$ be an ordinary GM threefold and $F\in M^X_G(2,1,5)$. Then $\tau_{\mathcal{A}}\mathrm{pr}(F)\cong\mathrm{pr}(\iota(F))$.
\end{proposition}

\begin{proof} \leavevmode
\begin{enumerate}
    \item If $F$ is a non-locally free sheaf determined by a line $L$, then by Lemma~\ref{pr'(F)=(iota F)^vee[1]} we have the triangle $$\mathcal{E}[1]\rightarrow\mathrm{pr}(F)\rightarrow\mathcal{O}_{L'}(-1).$$ Then since $\tau_{\mathcal{A}}\cong\tau_{\mathcal{A}}^{-1}\cong\bL_{\mathcal{O}_X}\circ\bL_{\mathcal{E}^{\vee}}(-\otimes\mathcal{O}_X(H))[-1]$, $\tau_{\mathcal{A}}(\mathrm{pr}(F))$ is given by a triangle $$\bL_{\mathcal{O}_X}\bL_{\mathcal{E}^{\vee}}(\mathcal{E}^{\vee})\rightarrow\tau_{\mathcal{A}}(\mathrm{pr}(F))\rightarrow \bL_{\mathcal{O}_X}\bL_{\mathcal{E}^{\vee}}(\mathcal{O}_{L'})[-1].$$ Note that $\bL_{\mathcal{E}^{\vee}}(\mathcal{E}^{\vee})=0$, hence $\tau_{A}(\mathrm{pr}(F))\cong\bL_{\mathcal{O}_X}\bL_{\mathcal{E}^{\vee}}(\mathcal{O}_{L'})[-1]$. It is easy to see $\mathrm{RHom}^\bullet(\mathcal{E}^{\vee},\mathcal{O}_{L'})=k$, therefore we have $\mathcal{E}^{\vee}\rightarrow\mathcal{O}_{L'}\rightarrow\bL_{\mathcal{E}^{\vee}}\mathcal{O}_{L'}$. Also, since $\mathcal{E}^{\vee}\rightarrow\mathcal{O}_{L'}$ is surjective, we have $\bL_{\mathcal{E}^{\vee}}\mathcal{O}_{L'}\cong F'[1]$, where $F':=\ker (\mathcal{E}^{\vee}\rightarrow\mathcal{O}_{L'})$ is a non-locally free sheaf in $M_G(2,1,5)$ determined by $L'$ as in Proposition~\ref{objects in M_G(2,1,5) theorem}. Thus $\tau_{\mathcal{A}}(\mathrm{pr}(F))\cong \bL_{\oh_X}F'\cong \pr(F')=\pr(\iota(F))$.
     
    \item If $F$ is a globally generated vector bundle, consider the standard short exact exact sequence
    $$0\rightarrow\ker(\ev)\rightarrow H^0(X,F)\otimes\mathcal{O}_X\xrightarrow{\ev} F\rightarrow 0 . $$ Dualizing the sequence and applying $\pr$, we get the triangle $$\mathrm{pr}(F^{\vee})\rightarrow\mathrm{pr}(\mathcal{O}_X^{\oplus 4})\rightarrow\mathrm{pr}(\ker(\ev)^{\vee})\cong\mathrm{pr}(\iota F).$$
    Note that $F^{\vee}\in\mathcal{A}_X$ and $\mathrm{pr}(\mathcal{O}_X)=0$, thus we get $\mathrm{pr}(\iota F)\cong F^{\vee}[1]$. 
Since $F\in M_G(2,1,5)$ is a globally generated vector bundle, we have $F\cong\iota E$ for some globally generated vector bundle $E$. Then $\mathrm{pr}(F)=\mathrm{pr}(\iota E)\cong E^{\vee}[1]\cong E\otimes\mathcal{O}_X(-H)[1]$, hence $\tau_{\mathcal{A}}(\mathrm{pr}(F))\cong\tau_{\mathcal{A}}(E\otimes\mathcal{O}_X(-H))[1]\cong\mathrm{pr}(E)\cong\mathrm{pr}(\iota F)$.
\end{enumerate}
\end{proof}

\subsection{The Bridgeland moduli space of class $y-2x$}

In this subsection, we show that $M_G(2,1,5)\cong \cM_{\sigma}(\mathcal{A}_X, y-2x)$. 

\begin{theorem} \label{irreducible component pr' theorem}
Let $X$ be a GM threefold and $\sigma$ be a Serre-invariant stability condition on $\cA_X$. Then the projection functor $\pr \colon \D^b(X) \ra \cA_X$ induces an isomorphism $M^X_G(2,1,5)\cong \cM_{\sigma}(\cA_X, y-2x)$.
\end{theorem}


We split the proof of this theorem into a series of lemmas and propositions. Recall that in \ref{Theorem_stabilitycondition_exists} we defined
\[V:=\{(\alpha, \beta)\colon -\frac{1}{10}<\beta<0, 0<\alpha<-\beta\}.\]

\begin{proposition}\label{mod215_connected}
Let $F\in \mathcal{A}(\alpha, \beta)$ be a $\sigma(\alpha, \beta)$-stable object with numerical class $y-2x$ for every $(\alpha, \beta)\in V$. Then $F=\mathrm{pr}(E)$ for some $E\in M_G(2,1,5)$.
\end{proposition}

\begin{proof}
First, we argue as in \cite[Proposition 4.6]{pertusi2020some}. When $(\alpha_0, \beta_0)=(0,0)$, we have $\mu^0_{\alpha_0, \beta_0}(F)=-\infty$. Since there are no walls intersecting with $\beta=0$ as in \cite[Proposition 4.6]{pertusi2020some}, we know that $F$ is $\sigma^0_{\alpha, 0}$-semistable for all $\alpha>0$. By the definition of the double-tilted heart, we have a triangle
\[A[1]\to F\to B\]
such that $A$ (respectively $B$) is in $\mathrm{Coh}^0(X)$ with its $\sigma_{\alpha, 0}$-semistable factors having slope $\mu_{\alpha, 0}\leq 0$ (respectively $\mu_{\alpha, 0}>0$). Since $F$ is $\sigma^0_{\alpha, 0}$-semistable and $\mu^0_{\alpha, 0}(F)<0$, we have that $A[1]=0$ and $B\cong F$. Since $\ch^0_1(F)$ is minimal, there are no walls on $\beta=0$, and we know that $F$ is $\sigma_{\alpha, 0}$-semistable for every $\alpha>0$. Thus by Lemma \ref{bms lemma 2.7}, $\mathcal{H}^{-1}(F)$ is a $\mu$-semistable reflexive sheaf and $\mathcal{H}^0(F)$ is $0$ or supported in dimension $\leq 1$. 

If $\mathcal{H}^0(F)$ is supported in dimension 0, then $\mathrm{ch}(\mathcal{H}^0(F))=bP$ for $b\geq 1$. But this is impossible since then $c_3(\cH^{-1}(F))>0$ and by \cite[Proposition 3.5]{brambilla2014vector} we have $\chi(\mathcal{H}^{-1}(F))=0$, which implies $b=0$. 

If $\mathcal{H}^0(F)$ is supported in dimension 1, we can assume $\mathrm{ch}(\mathcal{H}^0(F))=aL+\frac{b}{2}P$ where $a\geq 1$ and $b$ are integers. Thus $\mathrm{ch}(\mathcal{H}^{-1}(F))=2-H+aL+(\frac{5}{6}+\frac{b}{2})P$.  Now from Lemma \ref{non existence of sheaf}, we know $\mathcal{H}^{-1}(F)\cong \mathcal{E}$ and $\mathrm{ch}(\mathcal{H}^0(F))=L-\frac{P}{2}$. Thus $\mathcal{H}^0(F)\cong \mathcal{O}_L(-1)$ for some line $L$ on $X$. Therefore we have a triangle
\[\mathcal{E}[1]\to F\to \mathcal{O}_L(-1) . \]
In this case we have $\mathrm{Hom}(\mathcal{O}_L(-1), \mathcal{E}[2])=\mathrm{Hom}(\mathcal{E}^{\vee}(1), \mathcal{O}_L[1])=H^1(L, \mathcal{E}(-1)
|_L)=H^1(L, \mathcal{O}_L(-1)\oplus \mathcal{O}_L(-2))=k$.  Hence by Lemma \ref{pr'(F)=(iota F)^vee[1]}, $F\cong \mathrm{pr}(E)$ for some $E\in M_G(2,1,5)$ such that $E$ is locally free but not globally generated.

If $\mathcal{H}^0(F)=0$, we have $F[-1]\cong \mathcal{H}^{-1}(F)$. Then $F[-1]$ is a $\mu$-semistable sheaf. Since $F[-1]$ is reflexive and $c_3(F[-1])=0$, $F[-1]\in M_G(2,-1,5)$ is a stable vector bundle. Thus by Lemma \ref{pr'(F)=(iota F)^vee[1]}, we know  $F[-1]=\mathrm{pr}(E)$ for some $E\in M_G(2,1,5)$ such that $E$ is a globally generated vector bundle.
\end{proof}


\begin{lemma}
\label{prop_pr_injective}
The functor $\pr \colon \D^b(X) \ra \cA_X$ is injective on all objects in $M_G(2,1,5)$, i.e. if $\pr(F_1) \cong \pr(F_2)$, then $F_1 \cong F_2$.
\end{lemma}

\begin{proof}
For the case of globally generated vector bundles, by Corollary \ref{pr'(F)=(iota F)^vee[1]}, $\pr(F_1) \cong \pr(F_2)$ implies that  
\begin{equation*} \label{iota(F_1)^vee=iota(F_2)^vee}
(\iota F_1)^\vee \cong (\iota F_2)^\vee . 
\end{equation*} 
Note that $(\iota F_i)^{\vee}\cong\iota F_i\otimes\mathcal{O}_X(-H)$ for $i=1,2$. Then we get $\iota F_1\cong\iota F_2$. Finally, we apply $\iota$ to both sides. Since it is an involution $\iota^2= \mathrm{id}$, so $F_1 \cong F_2$ as required. 

For the case of non-locally free sheaves $F$, recall that from Lemma \ref{pr'(F)=(iota F)^vee[1]} we have $\cH^{-1}(\pr(F))=\cE$ and  $\cH^0(\pr(F)) = \oh_L(-H)$. Since  $F$ is uniquely determined by the line $L$, and $\Hom(\oh_L(-H), \cE[2])=k$, the object $\pr(F)$ is also uniquely determined by the line $L$. Thus $\pr(F_1) \cong \pr(F_2)$ implies $F_1 \cong F_2$, as required.
\end{proof}

\begin{proof}[Proof of Theorem~\ref{irreducible component pr' theorem}]
Using Proposition \ref{prop_selfext_group_F}, we know that the projection functor $\pr$ induces a morphism  \[ p \colon M_G(2,1,5)\to \mathcal{M}_\sigma(\cA_X, y-2x) \] which is bijective on points by Proposition~\ref{mod215_connected} and Lemma \ref{prop_pr_injective}, and bijective on tangent spaces by Lemma \ref{prop_pr_embedding}. Hence it is an isomorphism.
\end{proof}

\section{Refined and birational categorical Torelli theorems for GM threefolds} \label{torelli_sec}

In this section, we will prove several refined/birational categorical Torelli theorems for GM threefolds, using results from the previous sections.

\subsection{The universal family for \texorpdfstring{$\mathcal{C}_m(X)$}{the minimal surface of the Hilbert scheme of conics}}

In this subsection, we show that $\cC_m(X)\cong\mathcal{M}_{\sigma}(\cA_X,-x)$ admits a universal family, which thus gives a fine moduli space. Let $\mathcal{I}$ be the universal ideal sheaf of conics on $X\times\mathcal{C}(X)$ and $\mathcal{I}_{ L_{\sigma}}$ be the universal ideal sheaf of conics restricted to $X\times L_{\sigma}$. Let $q \colon X\times\mathcal{C}(X)\rightarrow X$ and $\pi \colon X\times\mathcal{C}(X)\rightarrow\mathcal{C}(X)$ be the projection maps on the first and second factors, respectively. Let $\mathcal{G}':=\mathrm{pr}(\mathcal{I}_{ L_{\sigma}})$ be the projected family in $\mathcal{A}_{X\times L_{\sigma}}$. Let $t\in L_{\sigma}\cong\mathbb{P}^1$ be any point. Then $j_t^*\mathrm{pr}(\mathcal{I}_{ L_{\sigma}})\cong A$, where $j_t \colon X_t\rightarrow X_t\times L_{\sigma}$ and $A\in\cA_X$ is $A\cong \pr(I_C)$ for $I_C\notin \cA_X$ by Proposition~\ref{prop_conic_inKuz}. Then $\mathcal{G}'\cong q^*(A)\otimes\pi^*\mathcal{O}_{ L_{\sigma}}(k)$ for some $k\in\mathbb{Z}$. Now let $\mathcal{G}:=\mathrm{pr}(\mathcal{I})\otimes\pi^*\mathcal{O}_{\mathcal{C}(X)}(kE)$, where $E\cong L_{\sigma}\cong\mathbb{P}^1$ is the unique exceptional curve on $\mathcal{C}(X)$. 

\begin{proposition}\label{universal-family}
The object $(p_X)_*\mathcal{G}$ is the universal family of $\mathcal{C}_m(X)$, where $p_X=\mathrm{id}_X \times p \colon X \times \cC(X) \rightarrow  X \times \cC_m(X)$.
\end{proposition}

\begin{proof}\leavevmode
\begin{enumerate}
    \item If $s=[A]=\pi\in\mathcal{C}_m(X)$, $s$ is contracted from the unique rational curve $ L_{\sigma}\cong\mathbb{P}^1\subset\mathcal{C}(X)$. Note that in this case $p_X|_{ L_{\sigma}}=q$. Then 
    \begin{align*}
        i_s^*(p_X)_*\mathcal{G} &\cong i_s^*(p_X)_*(\mathcal{G}'\otimes\pi^*\mathcal{O}_{\mathcal{C}(X)}(kE)) \\
        &\cong i_s^*q_*(q^*(A)\otimes\pi^*\mathcal{O}_{ L_{\sigma}}(k)\otimes\pi^*\mathcal{O}_{\mathcal{C}(X)}(kE)) \\
        &\cong i_s^*q_*(q^*(A)\otimes(\pi^*\mathcal{O}_{ L_{\sigma}}(k)\otimes\mathcal{O}_{ L_{\sigma}}(kE))) \\
        &\cong i_s^*q_*(q^*(A)\otimes\pi^*(\mathcal{O}_{ L_{\sigma}}(k)\otimes\mathcal{O}_{ L_{\sigma}}(-k))) \\
        &\cong i_s^*q_*(q^*(A)) \cong i_s^*(A) \cong A .
    \end{align*}
    \item If $s=[I_C]$, then $\mathcal{C}_m(X)$ and $\mathcal{C}(X)$ are isomorphic outside $ L_{\sigma}$. Note that $p$ restricts to $\mathrm{id}$ on $\mathcal{C}(X)\smallsetminus L_{\sigma}$.  Then 
    \begin{align*}
        i_s^*(p_X)_*\mathcal{G} &\cong i_s^*(p_X)_*(\mathrm{pr}(\mathcal{I})\otimes\pi^*\mathcal{O}_{\mathcal{C}(X)}(kE)) \\
        &\cong j_s^*(\mathrm{pr}(\mathcal{I}))\otimes j_s^*\pi^*\mathcal{O}_{\mathcal{C}(X)}(kE) \\
        &\cong I_C\otimes (\pi\circ j_s)^*\mathcal{O}_{\mathcal{C}(X)}(kE) \\
        &\cong I_C\otimes (i_s\circ\pi_s)^*\mathcal{O}_{\mathcal{C}(X)}(kE) \cong I_C .
    \end{align*}
\end{enumerate}
See below for the commutative diagrams which summarise the maps in the proof:
\[ \begin{tikzcd}[arrows={-Stealth}]
  X_s\rar["\cong"]\arrow[hookrightarrow]{d}{j_s}& X_s\rar["\pi_s"]\arrow[hookrightarrow]{d}{i_s} & \{s\}\dar \\%
X\times\mathcal{C}(X)\rar[swap, "p_X"] & X\times\mathcal{C}_m(X)\rar[swap, " "] & \mathcal{C}_m(X)
\end{tikzcd}
\]

\[ \begin{tikzcd}[arrows={-Stealth}]
  X_s\rar["j_s"]\arrow{d}{\pi_s}& X\times\mathcal{C}(X)\arrow{d}{\pi} \dar \\%
\{s\}\rar[swap, "i_s"] & \mathcal{C}(X)
\end{tikzcd}
\]

\end{proof}

\subsection{A refined categorical Torelli theorem for ordinary GM threefolds}

We now prove a refined categorical Torelli theorem for ordinary GM threefolds. 
\begin{theorem}\label{refined categorical torelli theorem}
Let $X$ and $X'$ be general ordinary GM threefolds such that $\Phi \colon \mathcal{K}u(X)\simeq \mathcal{K}u(X')$ is an equivalence and $\Phi(\pi(\cE))\cong\pi(\cE')$. Then $X\cong X'$.
\end{theorem}

\begin{proof}
Since $\Phi$ commutes with Serre functors, it preserves the stability of an object with respect to any Serre-invariant stability condition. Then the existence of the universal family on $\mathcal{C}_m(X)\cong \cM_{\sigma}(\cA_X, -x)$ guarantees a morphism from $\mathcal{C}_m(X)$ to $\mathcal{C}_m(X')$, denoted by $\psi$, which is induced by $\Psi$ (for more details on the construction of the morphism $\psi$, see \cite{bernardara2012categorical,altavilla2019moduli}). Since $\Phi$ is an equivalence, $\psi$ is an isomorphism. On the other hand, we have $\psi([\pi_X])=[\pi_{X'}]$ by the assumption, where $\pi_X:=\pi'(\cQ^{\vee})$ and $\pi_{X'}:=\pi'(\cQ^{\vee})$. Then $\psi$ induces an isomorphism $\phi \colon \mathcal{C}(X)\cong\mathcal{C}(X')$ by blowing up $[\pi_X]\in\mathcal{C}_m(X)$ and $[\pi_{X'}]\in\mathcal{C}_m(X')$, respectively. Then we have $X\cong X'$ by Logachev's Reconstruction Theorem \ref{logachev's reconstruction theorem}. 
\end{proof}

\subsection{Birational categorical Torelli theorem for ordinary GM threefolds} \label{birational categorical Torelli section}

In this subsection, we show a birational categorical Torelli theorem for ordinary GM threefolds, i.e. assuming the Kuznetsov components are equivalent leads to a birational equivalence of the ordinary GM threefolds.

\begin{theorem} \label{birational OGM torelli theorem}
Let $X$ and $X'$ be general ordinary GM threefolds such that $\mathcal{A}_X \simeq \mathcal{A}_{X'} $. Then $X'$ is a conic transform or a conic transform of a line transform of $X$. In particular, we have $X\simeq X'$.
\end{theorem}

\begin{proof}
The equivalence $\Phi\colon \mathcal{A}_X \xrightarrow{\sim} \mathcal{A}_{X'}$ sends $-x$ to either itself or $y-2x$ in $\cN(\mathcal{A}_{X'})$ up to sign, since they are only $(-1)$-class. By the same argument as in Theorem \ref{refined categorical torelli theorem} and \cite{bernardara2012categorical,altavilla2019moduli}, we thus get two possible induced isomorphisms between Bridgeland moduli spaces
\[
\begin{tikzcd}
{\cC_m(X)\cong\cM_\sigma(\mathcal{A}_X, -x)} \arrow[r, "\gamma"] \arrow[rd, "\gamma'"] & {\cC_m(X')\cong\cM_{\Phi(\sigma)}(\mathcal{A}_{X'}, -x)}  \\
                                                                  & {M_G^{X'}(2,1,5)\cong\cM_{\Phi(\sigma)}(\mathcal{A}_{X'}, y-2x)}
\end{tikzcd}
\]

If we have the isomorphism $\gamma$, then we blow up $\mathcal{C}_m(X)$ at the distinguished point $[\pi_X] := [\Xi(\pi(\cE))]$, and blow up $\mathcal{C}_m(X')$ at the point $[C]:=[\Phi(\pi_X)]=\gamma([\pi_X])$. We have $$\mathcal{C}(X) \cong \mathrm{Bl}_{[ \pi_X ]}\mathcal{C}_m(X)\cong \mathrm{Bl}_{[C]} \mathcal{C}_m(X'),$$ 
and $\mathrm{Bl}_{[C]} \mathcal{C}_m(X') \cong \mathcal{C}(X'_C)$ by Theorem \ref{DIM theorem 6.4}, so $\mathcal{C}(X) \cong \mathcal{C}(X_C')$. Therefore by Logachev's Reconstruction Theorem \ref{logachev's reconstruction theorem} we have $X \cong X'_C$.

For the second case, we get $\mathcal{C}_m(X)\cong M_G^{X'}(2,1,5)$. And by \cite[Proposition 8.1]{debarre2012period} we have a birational equivalence $M_G^{X'}(2,1,5) \simeq \mathcal{C}(X_L')$ of surfaces, where $L\subset X'$ is a line. Then we see $\mathcal{C}_m(X)$ is birationally equivalent to $\mathcal{C}(X'_L)$. Let $\mathcal{C}_m(X'_L)$ be the minimal surface of $\mathcal{C}(X'_L)$. Note that the surfaces here are all smooth surfaces of general type. By the uniqueness of minimal models of surfaces of general type, we get $\mathcal{C}_m(X)\cong\mathcal{C}_m(X'_L)$, which implies $X\cong (X'_L)_C\simeq X'$ for a conic $C\subset X'_L$ as in the first case. 
\end{proof}

\begin{remark}
Theorem \ref{birational OGM torelli theorem} proves a conjecture \cite[Conjecture 1.7]{kuznetsov2019categorical} of Kuznetsov--Perry for \emph{general ordinary} GM varieties of dimension $3$.
\end{remark}

In \cite{debarre2012period}, the authors proved that $\cC_m(X_L)$ is birational to $M^X_G(2,1,5)$. The following corollary shows that they are indeed isomorphic.


\begin{corollary} \label{line_trans}
Let $X$ be a general ordinary GM threefold, and $X_L$ be a line transform of $X$. Then we have $\cC_m(X_L)\cong M^X_G(2,1,5)$.
Moreover, this isomorphism commutes with involutions $\iota$ and $\iota'$ on both sides, thus giving an isomorphism $\cC_m(X_L)/\iota\cong M^X_G(2,1,5)/\iota'$.
\end{corollary}

\begin{proof}
By the same argument as in the proof of Theorem \ref{birational OGM torelli theorem}, we have $\cC_m(X_L)\cong \cC_m(X)$ or $\cC_m(X_L)\cong M^X_G(2,1,5)$. Note that $\cC_m(X_L)\cong \cC_m(X)$ implies that $X_L\cong X_C$ for some conic $C\subset X$ as in Theorem \ref{birational OGM torelli theorem}. But this is impossible by \cite[Remark 7.3]{debarre2012period}. Thus we always have $\cC_m(X_L)\cong M^X_G(2,1,5)$. The last statement follows from the fact that any equivalence between Kuznetsov components commutes with Serre functors, and the involutions on $\cC_m(X_L)$ and $M^X_G(2,1,5)$ can be induced by Serre functors up to shift by Propositions \ref{involution-acting-conic} and Proposition \ref{MG(2,1,5) involution compatibility}.
\end{proof}

Since the intermediate Jacobian $J(X)$ is invariant under conic and line transforms, as a corollary we have

\begin{corollary} \label{K_J}
Let $X$ and $X'$ be general ordinary GM threefolds. If $\Ku(X)\simeq \Ku(X')$, then we have $J(X)\cong J(X')$.
\end{corollary}

In fact, we can relax the assumptions on $X$ by looking at the singularities of Bridgeland moduli spaces.

\begin{theorem} \label{theorem_KP_conjecture}
Let $X$ and $X'$ be general GM threefolds (they can be either general ordinary or general special) and suppose their Kuznetsov components $\cA_{X}\simeq\cA_{X'}$ are equivalent. Then $X$ is birationally equivalent to $X'$.
\end{theorem}

\begin{proof}
First, we claim that if $X$ and $X'$ are general GM threefolds such that $\Phi\colon \cA_X\simeq\cA_{X'}$, then both $X$ and $X'$ are ordinary or special simultaneously. Indeed, we may assume $X'$ is ordinary and $X$ is special. Then the equivalence would identify the moduli space $\cC_m(X)\cong \mathcal{M}_{\sigma}(\cA_X,-x)$ of stable objects of class $-x$ in $\cA_X$ with either the moduli space $\cC_m(X')\cong \mathcal{M}_{\sigma'}(\cA_{X'},-x)$ or $M_G^{X'}(2,1,5)\cong\mathcal{M}_{\sigma'}(\cA_{X'}, y-2x)$. But $\mathcal{C}_m(X)$ has a unique singular point by Theorem \ref{theorem_irreduciblecomponent_conics}, and both $\mathcal{C}_m(X')$ and $M_G^{X'}(2,1,5)$ are smooth for $X'$ general by Theorem \ref{theorem_irreduciblecomponent_conics} and Theorem \ref{irreducible component pr' theorem}. This means that neither identification is possible, so the claim follows.

Now $X$ and $X'$ are both general ordinary or general special, hence the result follows from Theorem \ref{birational OGM torelli theorem} and \ref{SGM categorical torelli}.
\end{proof}

\begin{corollary} \label{prop_general_preserves}
Let $X$ and $X'$ be general GM threefolds such that one of them is ordinary and their Kuznetsov components $\cA_{X}\simeq\cA_{X'}$ are equivalent. Then they are both general ordinary and $X$ is birationally equivalent to $X'$.
\end{corollary}

\subsection{A categorical Torelli theorem for special GM threefolds} \label{SGM categorical torelli}
In this subsection, we show that the Kuznetsov component of a general special GM threefold $X$ determines the isomorphism class of $X$.

Recall from Section \ref{GM threefolds and their derived categories section} that every special GM threefold $X$ is a double cover of a degree $5$ index $2$ prime Fano threefold $Y$ branched over a quadric hypersurface $\mathcal{B}$ in $Y$. Since $X$ is smooth and general, $(\cB,h)$ is a smooth degree $h^2=10$ K3 surface with Picard number $1$. There is a natural geometric involution $\tau$ on $X$ induced by the double cover. The Serre functor on $\Ku(X)$ is given by $S_{\Ku(X)}=\tau\circ [2]$.

\begin{theorem}
\label{theorem_torelli_sGM}
Let $X$ and $X'$ be general special GM threefolds with $\Phi \colon \Ku(X)\simeq\Ku(X')$. Then $X\cong X'$. 
\end{theorem}

\begin{proof}
By \cite[Theorem 1.1, Section 8.2]{KP2017}, the equivariant triangulated category $\Ku(X)^{\mu_2}$ is equivalent to $\D^b(\cB)$, where $\mu_2$ is the group of square roots of $1$ generated by the involution $\tau$ acting on $\Ku(X)$. Assume there is an equivalence $\Phi \colon \Ku(X)\simeq  \Ku(X')$. Since $S_{\Ku(X)} \cong \tau[2]$ and $S_{\Ku(X')}\cong \tau'[2]$, $\Phi$ commutes with the involutions $\tau$ and $\tau'$ on $\Ku(X)$ and $\Ku(X')$, respectively. Then we get an induced equivalence
$$\Psi:\Ku(X)^{\mu_2}\simeq\Ku(X')^{\mu_2'}$$ where $\mu_2=\langle \tau\rangle$, $\mu_2'=\langle\Phi\circ\tau\circ\Phi^{-1}=\tau'\rangle$ and $\mu_2\cong\mu_2'$. 
Thus we have $\Psi \colon \D^b(\cB)\simeq \D^b(\cB')$. We know that $\cB$ and $\cB'$ are smooth projective surfaces with polarizations $h$ and $h'$, respectively, so $\Psi$ is a Fourier--Mukai functor by Orlov's Representability Theorem \cite[Theorem 2.2]{orlov1997equivalences}. Moreover, $(\cB,h)$ and $(\cB',h')$ are both Picard number 1 smooth projective K3 surfaces of degree $h^2=h'^2=10=2\times 5$. Then by \cite[Theorem 1.10]{Oguiso2002} and \cite[Corollary 1.7]{SHOY2002}, there is an isomorphism $\phi \colon \cB\cong \cB'$. Since they both have Picard number one, we obtain $\phi^*(h')=h$. On the other hand $Y_5$ is rigid \cite[\S~4.1]{kuznetsov2009derived}, which implies $X\cong X'$.
\end{proof}

\begin{remark}\label{another-proof-SGM}
Theorem~\ref{theorem_torelli_sGM} can also be proved via Bridgeland moduli spaces with respect to the Kuznetsov component $\cA_X$. The details are contained in another paper of authors \cite{jz2021brillnoether}.
\end{remark}

\section{The Debarre--Iliev--Manivel conjecture} \label{DIM conjecture section}

Let $\cX_{10}$ be the moduli space of smooth ordinary GM threefolds and $\cA_{10}$ be the moduli space of $10$-dimensional principal polarised abelian varieties. In \cite[pp. 3-4]{debarre2012period}, the authors make the following conjecture regarding the general fiber of the period map:

\begin{conjecture}[{\cite[pp. 3-4]{debarre2012period}}] \label{DIM conjecture}
A general fiber $\cP^{-1}([J(X)])$ of the period map $\cP \colon \cX_{10} \ra \cA_{10}$ at the intermediate Jacobian $J(X)$ of an ordinary GM threefold $X$ is the union of $\cC_m(X)/ \iota $ and a surface birationally equivalent to $M^X_G(2,1,5)/ \iota'$, where $\iota, \iota'$ are geometrically meaningful involutions.
\end{conjecture}

\begin{remark} \label{rmk_DIM_conj}
Note that by Corollary \ref{line_trans}, the surface birationally equivalent to $M_G(2,1,5)/ \iota'$ in \cite{debarre2012period}, parametrizing conic transforms of a line transform of $X$, is actually isomorphic to $M_G(2,1,5)/ \iota'$. Thus this conjecture predicts that a general fiber $\cP^{-1}([J(X)])$ is actually the disjoint union of $\cC_m(X)/ \iota $ and $M_G(2,1,5)/ \iota'$.
\end{remark}

We will prove a categorical analogue of this conjecture.
Consider the ``categorical period map"
\[ \cP_{\mathrm{cat}} \colon \cX_{10} \ra \{ \cA_X \} / \simeq, \, \, X\mapsto\cA_X\]
where $\cX_{10}$ is the moduli space of isomorphism classes of GM threefolds and $\{ \cA_X \} / \simeq$ is the set of equivalence classes of Kuznetsov components of GM threefolds.
Note that a global description of a ``moduli of Kuznetsov components" $\{ \cA_X \} / \simeq$ is not known, however, local deformations are controlled by the second Hochschild cohomology $\HH^2(\cA_X)$. The fiber of the ``categorical period map" $\mathcal{P}_{\mathrm{cat}}$ over $\cA_X$ for an ordinary GM threefold is defined as the isomorphism classes of all ordinary GM threefolds $X'$ such that $\cA_{X'}\simeq\cA_{X}$.


\begin{theorem}\label{categorical period map}
 The general fiber $\cP^{-1}_{\mathrm{cat}}([\cA_X])$ of the categorical period map over the alternative Kuznetsov component of an ordinary GM threefold $X$ is the union of $\mathcal{C}_m(X)/\iota$ and $M_G^X(2,1,5)/\iota'$ where $\iota, \iota'$ are geometrically meaningful involutions.
\end{theorem}

\begin{proof}
The general fiber $\cP^{-1}_{\mathrm{cat}}([\cA_X])$ of the categorical period map consists of GM threefolds $X'$ such that there is an equivalence of Kuznetsov components $\cA_{X'}\simeq\cA_{X}$. Then by Theorem~\ref{theorem_KP_conjecture}, $X'$ is also a general ordinary GM threefold. Thus by Theorem \ref{birational OGM torelli theorem} and Theorem \ref{dual_conj}, we know that  $\cA_{X'}\simeq \cA_X$ if and only if $X'$ is a conic transform of $X$, or a conic transform of a line transform of $X$. Then the result follows from Proposition \ref{DIM9.2} and Corollary \ref{line_trans}.
\end{proof}


The Kuznetsov components of prime Fano threefolds of index $1$ and $2$ are often regarded as categorical analogues of the intermediate Jacobians of these threefolds, and it is known that if there is a Fourier--Mukai type equivalence  $\mathcal{K}u(X)\simeq\mathcal{K}u(X')$ (or $\mathcal{A}_X\simeq\mathcal{A}_{X'}$), then $J(X)\cong J(X')$ by \cite{perry2020integral}. For the converse, we have the following result.

\begin{theorem}\label{categorical analogue of Jacobian-theorem}
For smooth prime Fano threefolds $X$, if $X$ is one of the following:\begin{itemize}
    \item $Y_d, \quad \quad \, \, \, \, 2\leq d\leq 5$
    \item $X_{2g-2}, \quad g=5, 7, 8, 9, 10, 12$,
\end{itemize}
then the intermediate Jacobian $J(X)$ uniquely determines the Kuznetsov component $\mathcal{K}u(X)$, i.e.~for another prime Fano threefold $X'$ of the same degree, if  $J(X)\cong J(X')$, then $\mathcal{K}u(X)\simeq\mathcal{K}u(X')$. 
\end{theorem}

\begin{proof}
If $X$ is an index $2$ prime Fano threefold $Y_d$ of degree $2\leq d\leq 5$, then the statement follows from the Torelli theorems for $Y_d$. Now let $X_d$ be a degree $d$ index one prime Fano threefold. If $X=X_8$, the statement follows from its Torelli theorem. If $X=X_{12}, X_{18}, X_{16}$, their intermediate Jacobians are Jacobians of curves: $J(X_{12})\cong J(C_7)$, $J(X_{16})\cong J(C_3)$, and $J(X_{18})\cong J(C_2)$. But $\mathcal{K}u(X_{12}) \simeq \D^b(C_7)$,  $\mathcal{K}u(X_{16})\simeq \D^b(C_3)$ and $\mathcal{K}u(X_{18})\simeq \D^b(C_2)$. Thus the statement follows from the classical Torelli theorem for curves. If $X=X_{14}$, the statement follows from the Kuznetsov conjecture for the pair $(Y_3, X_{14})$ \cite{kuznetsov2003derived} and the Torelli theorem for cubic threefolds. If $X=X_{22}$, the statement is trivial since $\Ku(X_{22})\cong \Ku(Y_5)$ (\cite{KPS2018}) and $Y_5$ is rigid, so  $\mathcal{K}u(X)\simeq\mathcal{K}u(X')$ is always true.  
\end{proof}


Therefore, it is natural to make the following conjecture:
\begin{conjecture}\label{categorical analogue of Jacobian}
Let $X$ be a prime Fano threefold of index one or two. Then the intermediate Jacobian $J(X)$ uniquely determines the Kuznetsov component $\mathcal{K}u(X)$, i.e.~for another prime Fano threefold $X'$ of the same degree, if  $J(X)\cong J(X')$, then $\mathcal{K}u(X)\simeq\mathcal{K}u(X')$. 
\end{conjecture}


Surprisingly, in the case of general ordinary GM threefolds, we can restate the Debarre--Iliev--Manivel Conjecture \ref{DIM conjecture} as Conjecture \ref{categorical analogue of Jacobian}:

\begin{proposition} \label{DIM conjecture equivalent prop}
The Debarre--Iliev--Manivel Conjecture \ref{DIM conjecture} is equivalent to Conjecture~\ref{categorical analogue of Jacobian} for general ordinary GM threefolds.
\end{proposition}

\begin{proof}
First, we assume that  Conjecture~\ref{categorical analogue of Jacobian} holds. Then by Corollary \ref{K_J} and Theorem~\ref{categorical period map}, the Debarre--Iliev--Manivel Conjecture \ref{DIM conjecture} holds. 

On the other hand, we assume the Debarre--Iliev--Manivel Conjecture \ref{DIM conjecture} holds. Then for any $X$ and $X'$ such that $J(X)\cong J(X')$, $X$ is either a conic transform of $X'$, or $X$ is a conic transform of a line transform of $X'$. In both cases, we have $\mathcal{K}u(X)\simeq\mathcal{K}u(X')$ by the Duality Conjecture Theorem \ref{dual_conj}. Thus Conjecture~\ref{categorical analogue of Jacobian} holds. 
\end{proof}

\begin{appendix}

\section{Uniqueness of Serre-invariant stability conditions}\label{appendix}

In this appendix, we aim to prove the uniqueness of Serre-invariant stability conditions on $\Ku(X)$ for several prime Fano threefolds $X$ (Theorem \ref{all_in_one_orbit}). We start with a general criterion for when two numerical stability conditions with the same central charge are equal. We always assume that any triangulated category is $k$-linear and of finite type, i.e.~$\sum_i \ext^i(A, B)<+\infty$ for any two objects $A, B$. Therefore, the Euler form and the numerical Grothendieck group are well-defined.

\begin{theorem}\label{theorem-general-criterion}
Let $\cD$ be a $k$-linear triangulated category of finite type. Assume that

\begin{enumerate}[label=(\Alph*).]

    \item $\chi(x,x)\leq 1-n$ for a positive integer $n$ and any non-zero $x\in \cN(\cD)$,

    \item there exists an object object $D$ satisfies
    \[n\leq \ext^1(D,D)<2n,\]


\end{enumerate}

Let $\sigma_1=(\cA_1, Z)$ be a numerical stability condition on $\cD$ and $D_1, D_2\in \cA_1$ be two $\sigma_1$-stable objects satisfying:

\begin{enumerate}[label=(C).]
     \item for any two objects $A,B\in \cD$, if $\phi^+_{\sigma_1}(B)<\phi^-_{\sigma_1}(A)$, then $\Hom(B,A[2])=0$,
\end{enumerate}
\begin{enumerate}[label=(D).]
    \item if $E$ is a $\sigma_1$-semistable object with $\chi(E,D_1)\geq 0$ and $\chi(E, D_2)\geq 0$, then there exist $k\in \{1,2\}$ such that $\chi(E,D_k)>0$ and $\mu_{\sigma_1}(E)<\mu_{\sigma_1}(D_k)$, and
\end{enumerate}
\begin{enumerate}[label=(E).]
    \item if $E$ is a $\sigma_1$-semistable object with $\chi(D_1,E)\geq 0$ and $\chi(D_2,E)\geq 0$, then there exist $k\in \{1,2\}$ such that $\chi(D_k, E)>0$ and  $\mu_{\sigma_1}(E)>\mu_{\sigma_1}(D_k)$.
\end{enumerate}

If $\sigma_2=(\cA_2, Z)$ is a numerical stability condition on $\cD$ satisfies (C), (D) and (E) such that $D_1$ and $D_2$ are $\sigma_2$-stable with $\phi_{\sigma_2}(D_1)=\phi_{\sigma_1}(D_1)$ and $\phi_{\sigma_2}(D_2)=\phi_{\sigma_1}(D_2)$, then $\sigma_1=\sigma_2$. 

\end{theorem}

We first prove several lemmas. By the same proof as in \cite[Lemma 2.5]{bayer:derived-auto-K3}, we have the following generalized version of Weak Mukai Lemma:

\begin{lemma}\label{lem-mukai}
Let $\cD$ be a $k$-linear triangulated category with finite-dimensional $\Hom$-space. Then for any exact triangle $A\to E\to B$ with $\Hom(A,B)=\Hom(B,A[2])=0$, we have
\[\ext^1(A,A)+\ext^1(B,B)\leq \ext^1(E,E).\]
\end{lemma}

\begin{lemma}\label{lem-general-property}
Let $\cD$ be a $k$-linear triangulated category of finite type satisfies (A). Assume that there is a stability condition $\sigma=(\cA, Z)$ on $\cD$ satisfies (C).

\begin{enumerate}
    \item The homological dimension of $\cA$ is at most $2$.

    \item For any exact triangle $A\to E\to B$ with $\phi^-_{\sigma}(A)>\phi^+_{\sigma}(B)$, we have
    \[\ext^1(A,A)+\ext^1(B,B)\leq \ext^1(E,E).\]

    \item For any non-zero object $A\in \cD$, we have $\ext^1(A,A)\geq n$.

    \item If a non-zero object $E$ is not $\sigma$-semistable, then any Harder--Narasimhan factor $A$ of $E$ satisfies
    \[\ext^1(A, A)<\ext^1(E, E).\]

    \item Any object $E$ with 
    \[n\leq \ext^1(E,E)<2n\]
is $\sigma$-semistable.
\end{enumerate}

\end{lemma}

\begin{proof}
Let $A,B\in \cA$. Then we have $\phi^+_{\sigma}(A)\leq 1<\phi^-_{\sigma}(B[k])$ for any $k\geq 1$. Therefore, by (C) we get $\Hom(A,B[k+2])=\Hom(A,B[k][2])=0$ for any $k\geq 0$. This proves (1).

Now for (2), note that $\Hom(A,B)=0$ and by (C) we have $\Hom(B,A[2])$. Then the result follows from Lemma \ref{lem-mukai}.

Next, we prove (3). If $A\neq 0\in \cA$, then from (1), we get $\chi(A,A)=\hom(A,A)-\ext^1(A,A)+\ext^2(A,A)$. Since $\chi(A,A)\leq 1-n$ by (A), we know that $\ext^1(A,A)\geq n$ in this case. Now for a general non-zero object $A\in \cD$, if $A$ is $\sigma$-semistable, then it is in $\cA$ up to shift and the result follows from the previous argument. So we assume that $A$ is not $\sigma$-semistable. Let $A'$ be the first Harder--Narasimhan factor of $A$ with respect to $\sigma$, and $A'':=\cone(A'\to A)$. We have $\phi_{\sigma}(A')>\phi^+_{\sigma}(A'')$. Using (2) and $\sigma$-semistability of $A'$, we obtain $n\leq \ext^1(A',A')+\ext^1(A'',A'')\leq \ext^1(A,A)$ and hence (3) follows. And (4) follows from the induction on the number of Harder--Narasimhan factors of $E$ and using (2) and (3).

Finally, if such $E$ in (5) is not $\sigma$-semistable, then by the existence of Harder--Narasimhan filtration, we can find a triangle $A\to E\to B$ with $\phi^-_{\sigma}(A)>\phi^+_{\sigma}(B)$. By (2) and (3), this contradicts our assumption on $\ext^1(E, E)$. Thus $E$ is $\sigma$-semistable. 
\end{proof}

Now we are ready to prove our criterion.

\begin{proof}[{Proof of Theorem \ref{theorem-general-criterion}}]
Since $\sigma_1$ and $\sigma_2$ have the same central charge, it remains to show $\cA_1=\cA_2$. By our assumptions, $D_1, D_2\in \cA_1\cap \cA_2$ are both $\sigma_1$-stable and $\sigma_2$-stable with phases in $(0,1]$.

\textbf{Step 1.} First, we show that if $E$ is a $\sigma_i$-semistable object which is also $\sigma_j$-semistable, then $\phi_{\sigma_1}(E)=\phi_{\sigma_2}(E)$, where $\{i, j\}=\{1, 2\}$. Since $\sigma_1$ and $\sigma_2$ satisfy the same assumptions, in the following, we will take $i=2$ and $j=1$. The other case can be deduced from the same argument but exchanges the role of $\sigma_1$ and $\sigma_2$.

Up to shift, we can assume that $E\in \cA_2$. Since $\sigma_1$ and $\sigma_2$ have the same central charge, we have $\phi_{\sigma_1}(E)-\phi_{\sigma_2}(E)=2m$ for an integer $m$. Then we see
\begin{equation}\label{eq-bound-m}
    2m<\phi_{\sigma_1}(E)\leq 2m+1.
\end{equation}

\begin{itemize}
    \item Assume that there exist $k,l\in\{1,2\}$ such that $\chi(D_k, E)<0$ and $\chi(E, D_l)<0$. Then by Lemma \ref{lem-general-property} (1) and the fact that $E,D_1,D_2\in \cA_2$, we have
    \[\ext^1(D_k, E)\neq 0, \quad \ext^1(E, D_l)\neq 0,\]
    which imply
    \[-1<\phi_{\sigma_1}(D_k)-1\leq \phi_{\sigma_1}(E)\leq \phi_{\sigma_1}(D_l)+1\leq 2.\]
    Hence, by \eqref{eq-bound-m} we get $2m<2$ and $-1<2m+1$, which means $m=0$ and we obtain $\phi_{\sigma_1}(E)=\phi_{\sigma_2}(E)$ as required.

    \item Assume that $\chi(E, D_1)\geq 0$ and $\chi(E, D_2)\geq 0$. By (D), there is an integer $s\in \{1,2\}$ such that $\chi(E, D_s)>0$ and $\mu_{\sigma_1}(E)<\mu_{\sigma_1}(D_s)$. Thus we have $\mu_{\sigma_2}(E)<\mu_{\sigma_2}(D_s)$ as well. Since $E, D_s\in \cA_2$, we get $\phi_{\sigma_2}(E)<\phi_{\sigma_2}(D_s)$, which implies $\Hom(E, D_s[2])=0$ by (C). Then from $\chi(E, D_s)>0$ and Lemma \ref{lem-general-property} (1), we obtain $\Hom(E, D_s)\neq 0$, and hence
    \[\phi_{\sigma_1}(E)\leq\phi_{\sigma_1}(D_s)\leq 1.\]
    
    Now if one of $\chi(D_1, E)$ and $\chi(D_2, E)$ is negative, the same argument as in the first case shows that $-1< \phi_{\sigma_1}(E)$. 
    
    If $\chi(D_1, E)\geq 0$ and $\chi(D_2, E)\geq 0$, then by (E), there is an integer $t\in \{1,2\}$ such that $\chi(D_t, E)>0$ and $\mu_{\sigma_1}(E)=\mu_{\sigma_2}(E)>\mu_{\sigma_1}(D_t)=\mu_{\sigma_2}(D_t)$. Since $E, D_t\in \cA_2$,  we get $\phi_{\sigma_2}(E)>\phi_{\sigma_2}(D_t)$, which by (C) implies $\Hom(D_t, E[2])=0$. Then together with $\chi(D_t, E)>0$ and Lemma \ref{lem-general-property} (1), we see $\Hom(D_t, E)\neq 0$. Therefore, we have $0<\phi_{\sigma_1}(D_t)\leq\phi_{\sigma_1}(E)$. In both cases, we always have $\phi_{\sigma_1}(E)\in (-1, 2]$. By \eqref{eq-bound-m}, we get $m=0$ and $\phi_{\sigma_1}(E)=\phi_{\sigma_2}(E)$ as required.

    \item Assume that $\chi(D_1, E)\geq 0$ and $\chi(D_2, E)\geq 0$. Then using (E), by a similar argument as the second case, we obtain $\phi_{\sigma_1}(E)=\phi_{\sigma_2}(E)$. This completes the first step.
\end{itemize}

\medskip

\textbf{Step 2.} Next, we prove that 
 an object $E$ is $\sigma_1$-semistable if and only if $\sigma_2$-semistable. We show this by induction on $\ext^1(E, E)$. If $\ext^1(E, E)<2n$, then from (B) we know such $E$ exists. By Lemma \ref{lem-general-property} (5), $E$ is both $\sigma_1$-semistable and $\sigma_2$-semistable.

Now assume that the statement holds for any object $F$ with $\ext^1(F, F)<N$ for an integer $N>2n$. Let $E$ be an object with $\ext^1(E,E)=N$. If $E$ is $\sigma_i$-semistable but not $\sigma_j$-semistable for $\{i, j\}=\{1, 2\}$, let $A$ be the first Harder--Narasimhan factor of $E$ with respect to $\sigma_j$ and $B$ be the last one. Therefore, we see
\begin{equation}\label{eq-app1}
    \phi_{\sigma_j}(A)>\phi_{\sigma_j}(B).
\end{equation}
And from Lemma \ref{lem-general-property} (4), we have 
\[\ext^1(A, A)< \ext^1(E, E), \quad \ext^1(B, B)<\ext^1(E, E).\]
Therefore, by the induction hypothesis, $A$ and $B$ are $\sigma_i$-semistable as well and the first step implies 
\begin{equation}\label{eq-app2}
    \phi_{\sigma_1}(A)=\phi_{\sigma_2}(A), \quad\phi_{\sigma_1}(B)=\phi_{\sigma_2}(B).
\end{equation}
But then from $\Hom(A, E)\neq 0$ and $\Hom(E, B)\neq 0$, we have $\phi_{\sigma_i}(A)\leq \phi_{\sigma_i}(E)\leq \phi_{\sigma_i}(B)$, which contradicts \eqref{eq-app1} and \eqref{eq-app2}. Hence $E$ is $\sigma_j$-semistable. This completes our induction argument.

\medskip

\textbf{Step 3.} Finally, by the previous two steps, we know that an object $E$ is $\sigma_1$-semistable if and only if $\sigma_2$-semistable with $\phi_{\sigma_1}(E)=\phi_{\sigma_2}(E)$. Since every non-zero object in the heart is obtained by extensions of semistable objects with phases in $(0,1]$, we know that $\cA_1=\cA_2$. This ends the proof of our theorem.
\end{proof}

\subsection{Applications to Kuznetsov components of Fano threefolds}

Let $Y_d$ be smooth index $2$ degree $d\geq 2$ prime Fano threefold and $X_{4d+2}$ an index $1$ degree $4d+2$ prime Fano threefold. In this section, we apply Theorem \ref{theorem-general-criterion} to show that all Serre-invariant stability conditions on $\mathcal{K}u(Y_d)$ and $\mathcal{K}u(X_{4d+2})$ (or $\mathcal{A}_{X_{4d+2}}$) are in the same $\widetilde{\mathrm{GL}}^+(2,\mathbb{R})$-orbit for each $d\geq 2$ (Theorem \ref{all_in_one_orbit}). 

Recall that the Kuznetsov component of an index two prime Fano threefold $Y_d$ of degree $d$ is defined by $\Ku(Y_d):=\langle \oh_{Y_d}, \oh_{Y_d}(H)\rangle^{\perp}$. The numerical Grothendieck group $\cN(\Ku(Y_d))$ is a rank two lattice generated by two classes
\[v=1-\frac{1}{d}H^2,\quad w=H-\frac{1}{2}H^2+(\frac{1}{6}-\frac{1}{d})H^2.\]
Moreover, under this basis, the Euler form is given by the matrix
\[ \begin{pmatrix}  
-1 & -1 \\
1-d & -d
\end{pmatrix} . \]


For index one cases, we assume that $d\geq 2$. Then the Kuznetsov component is defined by $\Ku(X_{4d+2}):=\langle \cE_{X_{4d+2}}, \oh_{X_{4d+2}}\rangle^{\perp}$, where $\cE_{X_{4d+2}}$ is a certain exceptional bundle pulled back from a Grassmannian (cf.~\cite{kuznetsov2003derived}).

By \cite{bayer2017stability}, $\sigma(\alpha, \beta)$ is a stability condition on $\Ku(Y_d)$ and $\Ku(X_{4d+2})$ for suitable $(\alpha, \beta)$. Moreover, according to \cite{pertusi2020some,Pertusi2021serreinv}, these stability conditions are all Serre-invariant.

Since for every index one prime Fano threefold $X_{4d+2}$ with $d\geq 3$, there is an index two prime Fano threefold $Y_d$ with $\Ku(Y_d)\simeq \Ku(X_{4d+2})$ by \cite{kuznetsov2003derived}, hence we only need to consider Kuznetsov components of $Y_d$ for $d\geq 2$ and $X_{10}$. Moreover, $\Ku(Y_4)$ is equivalent to the derived category of a smooth curve, and $\Ku(Y_5)$ is equivalent to the derived category of the $3$-Kronecker quiver. In these two cases, the result is known by \cite{macri2007stability} and \cite{dimitrov2019bridgeland}. So in the following, we mainly focus on $\cD=\Ku(Y_d)$ for $2\leq d\leq 3$ or $\Ku(X_{10})$. We first prove some properties of Serre-invariant stability conditions.

\begin{lemma}\label{lem-phase}
Let $\cD=\Ku(Y_d)$ for $2\leq d\leq 3$ or $\Ku(X_{10})$ and $\sigma=(\cA, Z)$ be a Serre-invariant stability condition on $\cD$. Then $\cD$ satisfies (A) and (B) in Theorem \ref{theorem-general-criterion} and $\sigma$ satisfies (C). Moreover, for any $\sigma$-semistable object $E\in \cD$, we have :

\begin{enumerate}
    \item if $\cD=\Ku(Y_3)$, then
    \[\phi_{\sigma}(E)+1\leq \phi_{\sigma}(S_{\cD}(E))<\phi_{\sigma}(E)+2,\]

    \item if $\cD=\Ku(Y_2)$ or $\Ku(X_{10})$, then
    \[\phi_{\sigma}(S_{\cD}(E))=\phi_{\sigma}(E)+2.\]
\end{enumerate}
\end{lemma}

\begin{proof}
It is clear that $\cD$ satisfies (A). And by \cite[Lemma 5.16]{pertusi2020some} and Proposition \ref{prop_selfext_group_F}, $\cD$ also satisfies (B). When $\cD=\Ku(Y_3)$, by \cite[Lemma 5.9]{pertusi2020some} we have $\phi_{\sigma}(S_{\cD}(E))<\phi_{\sigma}(E)+2$. When $\cD=\Ku(Y_2)$ or $\Ku(X_{10})$, recall that $S^2_{\cD}\cong [4]$. Then the same argument as in \cite[Lemma 5.9]{pertusi2020some} shows that $\phi_{\sigma}(S_{\cD}(E))\leq \phi_{\sigma}(E)+2$. Then for any two objects $A,B\in \cD$ with $\phi^+_{\sigma}(B)<\phi^-_{\sigma}(A)$, using (1) and (2) we see 
\[\phi^-_{\sigma}(A[2])=\phi^-_{\sigma}(A)+2>\phi^+_{\sigma}(B)+2\geq \phi^+_{\sigma}(S_{\cD}(B)).\]
Hence $\Hom(B, A[2])=\Hom(A[2], S_{\cD}(B))=0$ and the condition (C) is satisfied.

When $\cD=\Ku(Y_3)$, from \cite[Lemma 5.11]{pertusi2020some}, we get $\ext^1(E,E)\neq 0$, which implies $\phi_{\sigma}(E)+1\leq \phi_{\sigma}(S_{\cD}(E))$. This proves (1).

When $\cD=\Ku(Y_2)$ or $\Ku(X_{10})$, since $\cD$ satisfies (A) and (C), by Lemma \ref{lem-general-property} we have $\ext^1(E, E)\neq 0$, which implies $\phi_{\sigma}(E)+1\leq \phi_{\sigma}(S_{\cD}(E))$. Now since $[E]=[S_{\cD}(E)]\in \cN(\cD)$, we have $\phi_{\sigma}(E)-\phi_{\sigma}(S_{\cD}(E))\in 2\ZZ$. Hence we get $\phi_{\sigma}(S_{\cD}(E))= \phi_{\sigma}(E)+2$.

\end{proof}

Before verifying (D) and (E), we need several lemmas.

\begin{lemma}\label{lem-incidence}
Let $X$ be a GM threefold.

\begin{enumerate}
    \item $\mathrm{Hilb}_X^{3t+m}=\varnothing$ for $m<1$. Thus for any conic $C\subset X$ and line $L\subset X$, we have $\Hom(I_C, \oh_L(-k))=0$ for any $k>1$.

    \item If a line $L$ and a conic $C$ satisfies $L\cap C\neq \varnothing$, then $L\cap C$ is of length one and $L\cup C$ is a twisted cubic.

    \item Let $\mathcal{I}\subset \Gamma(X)\times \cC(X)$ be the incidence variety, i.e.
\[\cI=\{(L, C)~\colon~ L\cap C\neq \varnothing\}.\]
Then the projection maps $\cI\to \cC(X)$ and $\cI\to \Gamma(X)$ are surjective.
\end{enumerate}
\end{lemma}

\begin{proof}
By \cite[Corollary 1.38]{sanna2014rational}, we have $\mathrm{Hilb}_X^{3t+m}=\varnothing$ for $m<0$. Thus to prove (1), we only need to show $\mathrm{Hilb}_X^{3t}=\varnothing$. From \cite[Corollary 1.38]{sanna2014rational}, $\langle C\rangle\cong \PP^2$ for any $[C]\in \mathrm{Hilb}_X^{3t}$. Since $X$ is an intersection of quadrics, such a $C$ cannot exist on $X$. Hence $\mathrm{Hilb}_X^{3t}=\varnothing$. Now note that the kernel of any non-zero map $I_C\to \oh_L(-k)$ is the ideal sheaf of a closed subscheme with the Hilbert polynomial $3t+m$ for $m\leq 2-k$. Therefore, $\Hom(I_C, \oh_L(-k))=0$ when $k>1$. This proves (1). For (2), note that $\chi(\oh_{L\cup C})=2-\mathrm{length}(L\cap C)$, then the result follows from (1).

Finally, we prove (3). Since $\dim \Gamma(X)=1$, all lines on $X$ sweep out a surface $S$ in $X$. By $\Pic(X)=\ZZ\oh_X(H)$, we see $S\in |mH|$ for $m>0$. Thus $C.S\geq C.H>0$ for any conic $C$. In other words, $C\cap S\neq \varnothing$, hence any conic on $X$ intersects with a line. Thus $\cI\to \cC(X)$ is surjective. Similarly, since $X$ is covered by conics, any line intersects with a conic. Then $\cI\to \Gamma(X)$ is surjective.
\end{proof}

\begin{lemma}
\label{lemma_existence_line_cubic}
Let $X$ be a GM threefold. Then there exists a line $L$ and twisted cubics $C$ and $D$ on $X$ such that

\begin{enumerate}
    \item $[L]\in \Gamma(X)$ is a smooth point,

    \item $I_C\notin \Ku(X)$ and $L\cup C=Z(s)$ for a section $s\in H^0(\cE^{\vee})$,

    \item $L\subset D$, $I_D\in \Ku(X)$ and $\ext^1(I_D, I_D)=3$.
\end{enumerate} 
\end{lemma}

\begin{proof}
Let $\mathcal{I}\subset \Gamma(X)\times \cC(X)$ be the incidence variety. We denote by $\cC_1$ the sublocus of $\cC(X)$ parametrizing smooth conics $Z$ such that their involutive conics are also smooth and $\hom(\cE, I_Z)=1$. By Remark \ref{rmk-open-locus}, $\cC_1$ is an open subscheme of $\cC(X)$. Let $\cI_1:=\cI|_{\Gamma(X)\times \cC_1}$. From \cite[Theorem 3.4 (iii)]{iskovskih:fano-3fold-II} and \cite[Section 3.1]{iliev:lines-on-gushel-threefold}, $\Gamma(X)$ is generic smooth. This implies that the image of $p\colon\cI_1\to \Gamma(X)$ contains a smooth point. 

Let $L\subset X$ be a line such that $[L]\in \Gamma(X)$ is smooth and contained in the image of $p$. Then $p^{-1}([L])$ is non-empty and there is a conic $[Z]\in \cC_1$ such that $L\cap Z\neq \varnothing$. We set $D:=L\cup Z$. And since $\Hom(\cE, I_L)\neq 0$, there is a section $s\in H^0(\cE)$ such that $L\subset Z(s)$. We define $C$ to be the residue curve of $L$ in $Z(s)$. It is clear that $C$ and $D$ are twisted cubics by Lemma \ref{lem-incidence}. Moreover, $L$ and $Z$ intersect transversely at a single point. Then it remains to check $I_C\notin \Ku(X)$, $\ext^1(I_D, I_D)=3$ and $I_D\in \Ku(X)$.

Since $C\subset Z(s)$, it is clear that $\Hom(\cE, I_C)\neq 0$, i.e.~$I_C\notin \Ku(X)$. Moreover, by the construction, we have an exact sequence
\begin{equation}\label{eq-ID}
    0\to I_D\to I_Z\to \oh_L(-1)\to 0.
\end{equation}
Note that all conics are connected, hence the smoothness implies irreducibility. Since $\Hom(\cE, I_Z)=k$ and the involutive conic $Z'$ is smooth, we know that $Z\cup Z'$ only has two irreducible components which are both of degree $2$, hence does not contain $D$. This means the unique non-zero map in $\Hom(\cE, I_Z)=k$ does not factor through $I_D$. Hence the induced map $\Hom(\cE, I_Z)\to \Hom(\cE, \oh_L(-1))$ is injective. By $\RHomb(\cE, \oh_L(-1))=k$, this map is actually an isomorphism. Therefore, applying $\Hom(\cE, -)$ to \eqref{eq-ID}, we obtain $\RHomb(\cE, I_D)=0$, which implies $I_D\in \Ku(X)$.

To show $\ext^1(I_D, I_D)=3$, since $\hom(I_D, I_D)=1$ and $\ext^3(I_D, I_D)=0$, by $\chi(I_D, I_D)=-2$ we only need to prove $\ext^2(I_D, I_D)=0$. From the construction above, we see $\ext^2(I_Z, I_Z)=0$. Moreover, $\ext^2(\oh_L(-1), \oh_L(-1))=0$ since $[L]\in \Gamma(X)$ is a smooth point. And by the transversality of the intersection of $L$ and $Z$, we see the derived restriction $\oh_Z|_L\cong \oh_{L\cap Z}\in \D^b(L)$. Hence $\ext^2(\oh_Z, \oh_L(-1))=\ext^1(I_Z, \oh_L(-1))=0$. Finally, by Lemma \ref{lem-incidence} we have $\hom(I_Z, \oh_L(-2))=\ext^3(\oh_L(-1), I_Z)=0$. Then $\ext^2(I_D, I_D)=0$ follows from applying \cite[Lemma 2.27]{pirozhkov2020admissible} to \eqref{eq-ID}.
\end{proof}

\begin{lemma}\label{lem-F-stability}
Let $X$ be a GM threefold and $L,C,D$ as in Lemma \ref{lemma_existence_line_cubic}. We define $F_1:=\pr'(I_C), F_1':=I_D$ and $F_2:=\pr'(I_L)$. Then the objects $F_1, F_1', F_2$ are stable with respect to any Serre-invariant stability condition on $\Ku(X)$. Moreover, $F_1$ and $F_1'$ have the same phase.
\end{lemma}

\begin{proof}
From the construction, we see $\ext^1(F_1', F_1')=3$. By the same argument as in \cite[Corollary 5.4]{Zhang2020Kuzconjecture}, we have $\ext^1(F_1, F_1)=3$. Finally, applying \cite[Lemma 2.27]{pirozhkov2020admissible} to $\cE^{\oplus 2}\to I_L\to F_2$ and using $\RHomb(\oh_L, \oh_L)=k\oplus k[-1]$ implies $\ext^1(F_2, F_2)=3$. Then the stability of $F_1, F_1'$ and $F_2$ follows from Proposition \ref{ext23-stable}. 

As $[F_1]=[F_1']\in \cN(\Ku(X))$, we have $\chi(F_1, F_1')<0$. Since $\ext^i(F_1, F_1')=\ext^i(I_C, I_D)=0$ for $i\notin \{1, 2\}$, we get $\Hom(F_1, F_1'[1])=\Hom(F_1'[1], S_{\Ku(X)}(F_1))\neq 0$. Using Lemma \ref{lem-phase}, we obtain $\phi_{\sigma}(F_1')-1<\phi_{\sigma}(F_1)<\phi_{\sigma}(F_1')+1$ for any Serre-invariant stability condition $\sigma$. Thus $\phi_{\sigma}(F_1)=\phi_{\sigma}(F_1')$ since $[F_1]=[F_1']\in \cN(\Ku(X))$.
\end{proof}

\begin{lemma}\label{lem-central-charge}
Let $\cD=\Ku(Y_d)$ for $2\leq d\leq 3$ or $\Ku(X_{10})$. Then there exist two objects $F_1, F_2\in \cD$ such that for any Serre-invariant stability condition $\sigma$ on $\cD$,  $F_1$ and $F_2$ are $\sigma$-stable with
\[\phi_{\sigma}(F_2)-1<\phi_{\sigma}(F_1)<\phi_{\sigma}(F_2).\]

In particular, the image of the central charge is not contained in a line for any Serre-invariant stability condition on $\cD$.
\end{lemma}

\begin{proof}
When $\cD=\Ku(Y_d)$, we define $F_2:=R\mathcal{H}om(I_L, \oh_{Y_d}(-H))[1]$ and $F_1=I_L$, where $L\subset Y_d$ is a line. Then by \cite[Lemma 5.13]{pertusi2020some} and \cite[Remark 4.8]{pertusi2020some}, $F_1$ and $F_2$ are $\sigma$-stable for any Serre-invariant stability condition $\sigma$ on $\cD$ with $\phi_{\sigma}(F_2)-1<\phi_{\sigma}(F_1)<\phi_{\sigma}(F_2).$

Now assume that $\cD=\Ku(X_{10})$. We take $F_1, F_1'$ and $F_2$ as in Lemma \ref{lem-F-stability}. By \cite[Proposition 3.3, 5.3]{Zhang2020Kuzconjecture}, we have $\pr'(I_C)\cong \pr'(G)$, where $G$ fits into an exact triangle
\[\oh_X(-H)[1]\to G\to \oh_L(-2)\]
and is the twisted derived dual of the line $L$.

First, we prove that $\Hom(F_2, F_1[1])\neq 0$. By adjunction, we have $\Hom(F_2, F_1[1])=\Hom(I_L, \pr'(G)[1])$. And by \cite[Proposition 5.3]{Zhang2020Kuzconjecture},  $\pr'(G)$ fits into an exact triangle
\begin{equation}\label{eq-prG}
    G\to \pr'(G)\to \cE.
\end{equation}
Since $\cE|_L\cong \oh_L\oplus \oh_L(-1)$, it is easy to see $\Ext^i(I_L, \cE)$ for $i\neq 2$. So applying $\Hom(I_L, -)$ to \eqref{eq-prG}, we get $\Hom(I_L, \pr'(G)[1])=\Hom(I_L, G[1])=\Hom(I_L, \oh_L(-2)[1])$. As the normal bundle $N_{L/X_{10}}$ is either $\oh_L\oplus \oh_L(-1)$ or $\oh_L(1)\oplus \oh_L(-2)$ by \cite[Lemma 4.2.1]{iskovskikh1999fano}, we see the derived restriction $I_L|_L\cong N_{L/X}^{\vee}\oplus \oh_L(1)[1]$ from \cite[Proposition 11.8]{huybrechts2006fourier}. Then $\Hom(I_L, \oh_L(-2)[1])\neq 0$ follows from a direct computation.

Next, we show that $\Hom(F_1', F_2)\neq 0$. By the definition of $\pr'$, $\pr'(I_L)$ fits into an exact triangle $\cE^{\oplus 2}\to I_L\to \pr'(I_L)$. Then applying $\Hom(I_D, -)$ to this triangle, the result follows from $\Hom(I_D, \cE)=0$ and $\Hom(I_D, I_L)\neq 0$ since $L\subset D$.

By Lemma \ref{lem-F-stability}, $F_1, F_1'$ and $F_2$ are all stable with respect to any Serre-invariant stability condition on $\cD$. Therefore, combined with above results we get $\phi_{\sigma}(F_1)=\phi_{\sigma}(F_1')<\phi_{\sigma}(F_2)<\phi_{\sigma}(F_1)+1$ as desired.
\end{proof}

Now we are ready to verify conditions (D) and (E) in Theorem \ref{theorem-general-criterion}.

\begin{lemma}\label{lem-cond-DE}
Let $\cD=\Ku(Y_d)$ for $2\leq d\leq 3$ or $\Ku(X_{10})$. Then there exists a Serre-invariant stability condition $\sigma_1$ on $\cD$ with two $\sigma_1$-stable objects $D_1, D_2$ satisfying (D) and (E). Moreover,

\begin{itemize}
    \item we can assume that for any Serre-invariant stability condition $\sigma$ on $\cD$, $D_1$ and $D_2$ are $\sigma$-stable with
\[\phi_{\sigma}(D_1)-1<\phi_{\sigma}(D_2)<\phi_{\sigma}(D_1),~\text{or}~~\phi_{\sigma}(D_1)<\phi_{\sigma}(D_2)<\phi_{\sigma}(D_1)+1,\]
\item any Serre-invariant stability condition on $\cD$ with the same central charge as $\sigma_1$ satisfies (D) and (E).
\end{itemize}

\end{lemma}

\begin{proof}
When $\cD=\Ku(Y_d)$, we define $D_1:=R\mathcal{H}om(I_L, \oh_{Y_d}(-H))[1]$ and $D_2=I_L[1]$, where $L\subset Y_d$ is a line. Then by \cite[Lemma 5.13]{pertusi2020some} and \cite[Remark 4.8]{pertusi2020some}, $D_1$ and $D_2$ are $\sigma$-stable for any Serre-invariant stability condition $\sigma$ on $\cD$ with $\phi_{\sigma}(D_1)<\phi_{\sigma}(D_2)<\phi_{\sigma}(D_1)+1.$ In this case, we take $\sigma_1:=\sigma(\alpha, -\frac{1}{2})$ for $\alpha>0$ sufficiently small. Then by \cite[Section 4]{pertusi2020some}, $D_1, D_2\in \cA(\alpha, -\frac{1}{2})$. Now a direct computation shows that, for any object $E$ with $[E]=av+bw$, we have
\begin{itemize}
    \item $\chi(E,D_2)=a+(d-1)b$, $\chi(D_2,E)=a+b$; and $\mu^0_{\alpha, -\frac{1}{2}}(E)>\mu^0_{\alpha, -\frac{1}{2}}(D_2) \iff b<0$
    
    \item $\chi(E,D_1)=-b$, $\chi(D_1,E)=-[(d-2)a+(d-1)b]$; and  $\mu^0_{\alpha, -\frac{1}{2}}(E)>\mu^0_{\alpha, -\frac{1}{2}}(D_1) \iff a+b<0$.
\end{itemize}
Then it is straightforward to check (D) and (E) for $\sigma_1$.

When $\cD=\Ku(X_{10})$, we use the equivalence $\Xi$ in Lemma \ref{equivalence of original and alternative kuznetsov component} and prove every thing on $\cA_{X_{10}}$. Let $\sigma_1:=\sigma(\alpha, \beta)$, where $\beta<0$ and $\alpha>0$ with $-\beta$ and $\alpha$ are sufficiently small. We set $D_1=I_C[1]$ and $D_2=\pr(F)[1]$, where $C\subset X$ is a smooth conic with $I_C\in \cA_{X_{10}}$ and $F\in M_G(2,1,5)$ is non-locally free. It is clear that $D_1, D_2\in \cA(\alpha, \beta)$ and are stable with respect to any Serre-invariant stability condition on $\cA_{X_{10}}$ by Lemma \ref{beta=0} and Proposition \ref{prop_selfext_group_F}. As in the previous case, it is straightforward to check (D) and (E) for $\sigma_1$. Now we show that for any Serre-invariant stability condition $\sigma$ on $\cA_{X_{10}}$, we have $\phi_{\sigma}(D_1)-1<\phi_{\sigma}(D_2)<\phi_{\sigma}(D_1).$ Indeed, if $\sigma=\sigma_1$, then this follows from a direct computation of the slope function of $\sigma_1$. When $\sigma\neq\sigma_1$, by Lemma \ref{lem-central-charge}, up to $\GL$-action, we can assume that $\sigma$ and $\sigma_1$ has the same central charge and $\phi_{\sigma}(D_1)=\phi_{\sigma_1}(D_1)$. Thus $\phi_{\sigma}(D_2)-\phi_{\sigma_1}(D_2)\in 2\ZZ$. We claim that $\phi_{\sigma}(D_1)-2< \phi_{\sigma}(D_2)<\phi_{\sigma}(D_1)$, which implies $\phi_{\sigma}(D_2)=\phi_{\sigma_1}(D_2)$ and the result follows. Indeed, by Proposition \ref{objects in M_G(2,1,5) theorem} we have an exact sequence $0\to F\to \cE\to \oh_L(-1)\to 0$ for a line $L\subset X_{10}$. Hence applying $\Hom(-,D_1)$ to this exact sequence and use $\Hom(\cE, D_1[-1])\neq 0$ (Lemma \ref{IC-class}) and adjunction of $\pr$, we have $\Hom(F, I_C)=\Hom(D_2, D_1)=\Hom(D_1, S_{\cA_{X_{10}}}(D_2))\neq 0$. Then by Lemma \ref{lem-phase} we obtain $\phi_{\sigma}(D_1)-2< \phi_{\sigma}(D_2)<\phi_{\sigma}(D_1)$ as desired.

The final statement follows from the fact that (D) and (E) in this case only depend on the central charge and numerical classes $[D_1]$ and $[D_2]$, as we have seen above.
\end{proof}

Applying Theorem \ref{theorem-general-criterion}, we obtain the uniqueness of Serre-invariant stability conditions.

\begin{theorem}\label{all_in_one_orbit}
Let $\cD=\Ku(Y_d)$ for $2\leq d\leq 3$ or $\Ku(X_{10})$. Then all Serre-invariant stability conditions on $\cD$ are in the same $\GL$-orbit.
\end{theorem}

\begin{proof}
Let $\sigma_1$, $D_1$ and $D_2$ as in Lemma \ref{lem-cond-DE}. Let $\sigma_2$ be another Serre-invariant stability condition on $\cD$. By Lemma \ref{lem-central-charge}, up to $\GL$-action, we can assume that $\sigma_1$ and $\sigma_2$ have the same central charge. Moreover, up to shift we can assume that $\phi_{\sigma_1}(D_1)=\phi_{\sigma_2}(D_1)$. Thus $\phi_{\sigma_1}(D_2)-\phi_{\sigma_2}(D_2)\in 2\ZZ$. And by Lemma \ref{lem-cond-DE}, $\sigma_2$ also satisfies (D) and (E).

Now from Lemma \ref{lem-cond-DE}, we have $\phi_{\sigma_k}(D_1)-1<\phi_{\sigma_k}(D_2)<\phi_{\sigma_k}(D_1)$ or $\phi_{\sigma_k}(D_1)<\phi_{\sigma_k}(D_2)<\phi_{\sigma_k}(D_1)+1$ for any $k\in \{1, 2\}$. This implies $\phi_{\sigma_1}(D_2)=\phi_{\sigma_2}(D_2)$. Therefore, by Lemma \ref{lem-phase} and Lemma \ref{lem-cond-DE}, we can apply Theorem \ref{theorem-general-criterion} and get $\sigma_1=\sigma_2$.
\end{proof}

\begin{remark}\label{remark_history_uniqueness}
The idea of the proof of Theorem~\ref{all_in_one_orbit} was first explained to us by Arend Bayer. In \cite[Proposition 4.21]{Zhang2020Kuzconjecture}, one of the authors made an attempt to prove this statement but the argument is incomplete. Here, we fill the gaps and give a more general argument. Later, in \cite[Theorem 3.1]{FeyzbakhshPertusi2021stab}, the authors also prove the uniqueness of Serre-invariant stability conditions for a general triangulated category satisfying a list of assumptions and include Kuznetsov components of cubic threefolds and very general cubic fourfolds. The assumptions used in  \cite[Theorem 3.1]{FeyzbakhshPertusi2021stab} are (A), (B), and the Serre functor of $\cD$ satisfies $S^r_{\cD}=[k]$ with $0<k/r<2$ or $r=2$ and $k=4$, while our Theorem \ref{theorem-general-criterion} also works for general triangulated categories which are not fractional Calabi--Yau but with extra assumptions (C), (D) and (E). Indeed, if we take $D_1=D$ and $D_2=S_{\cD}(D)[-2]$  in (D) and (E) where $D$ is an object in (B), then one can show that when $k/r<2$, Theorem \ref{theorem-general-criterion} implies \cite[Theorem 3.1]{FeyzbakhshPertusi2021stab}. Moreover, Theorem \ref{theorem-general-criterion} can be applied to the derived category of a smooth projective curve or a generalized Kronecker quiver as well.
\end{remark}

\end{appendix}

\bibliographystyle{alpha}
{\small{\bibliography{hochschild}}}

\end{document}